\documentclass[12pt]{amsart}
\usepackage[dvips]{graphicx}

   {\setlength{\parsep}{0pt}%
    \begin{list}{}{\setlength{\topsep}{0pt}\setlength{\partopsep}{0pt}%
                   \setlength{\leftmargin}{2\parindent}%
                   \item }} {\end{list}}

\newenvironment{hangpar}%
   {\setlength{\parsep}{0pt}%
    \begin{list}{}{\setlength{\topsep}{0pt}\setlength{\partopsep}{0pt}%
                   \setlength{\leftmargin}{4\parindent}%
                   \setlength{\itemindent}{-\leftmargin}%
                   \item }} {\end{list}}

\newcommand{\myfigure}[1]{\goodbreak\begin{figure}[!ht]#1\end{figure}}
\newcommand{\mytable}[1]{\goodbreak\begin{table}[!ht]#1\end{table}}
\newtheorem{thm}{Theorem}[section]
\newtheorem{prop}[thm]{Proposition}
\newtheorem{cor}[thm]{Corollary}
\newtheorem{lem}[thm]{Lemma}
\newtheorem{conj}[thm]{Conjecture}
\newtheorem{exa}[thm]{Example}
\newtheorem{defn}[thm]{Definition}

\newtheorem{rem}[thm]{Remark}
\newtheorem{note}[thm]{Notation}
\newtheorem{alg}[thm]{Algorithm}

\newcommand{\ben}{\begin{enumerate}}
\newcommand{\een}{\end{enumerate}}
\newcommand{\ble}{\begin{lem}}
\newcommand{\ele}{\end{lem}}
\newcommand{\bth}{\begin{thm}}
\newcommand{\eth}{\end{thm}}
\newcommand{\bpr}{\begin{prop}}
\newcommand{\epr}{\end{prop}}
\newcommand{\bco}{\begin{cor}}
\newcommand{\eco}{\end{cor}}
\newcommand{\bcon}{\begin{conj}}
\newcommand{\econ}{\end{conj}}
\newcommand{\bde}{\begin{defn}}
\newcommand{\ede}{\end{defn}}
\newcommand{\bex}{\begin{exa}}
\newcommand{\eex}{\end{exa}}
\newcommand{\brem}{\begin{rem}}
\newcommand{\erem}{\end{rem}}
\newcommand{\bnot}{\begin{note}}
\newcommand{\enot}{\end{note}}
\newcommand{\balg}{\begin{alg}}
\newcommand{\ealg}{\end{alg}}
\newcommand{\bform}{\noindent \begin{hangpar} $}
\newcommand{\eform}{ $ \end{hangpar}}

\newcommand{\<}{\left <}
\renewcommand{\>}{\right > }

\newcommand{\C}{\mathbb{C}}
\newcommand{\PP}{\mathbb{P}}
\newcommand{\RR}{\mathbb{R}}

\title[BM of Real Singular Curves]{Braid Monodromy Computation of Real Singular Curves}

\author[Kaplan]{Kaplan Shmuel}
\address{Department of Mathematics and Statistics, Bar-Ilan University, Ramat-Gan 52900, Israel}
\email{kaplansh@macs.biu.ac.il}

\author[Liberman]{Liberman Eran}
\thanks{This paper is part of the First and Second authors PhD. Thesis in Bar-Ilan university}
\thanks{Partially supported by the Emmy Noether Research Institute for
Mathematics, Bar-Ilan University and the Minerva Foundation, Germany, and by the Excellency Center ``Group %%@
theoretic methods in the study of algebraic varieties'' of the National Science Foundation of Israel, and by Eager, European Algebraic Geometry Research.}
\address{Department of Mathematics and Statistics, Bar-Ilan University, Ramat-Gan 52900, Israel}
\email{liberme1@macs.biu.ac.il}

\author[Teicher]{Teicher Mina}
\address{Department of Mathematics and Statistics, Bar-Ilan University, Ramat-Gan 52900, Israel}
\email{teicher@macs.biu.ac.il}

\date{\today}

\begin{document}

\begin{abstract}
We generalize the Moishezon Teicher algorithm that was suggested for the computation of the braid monodromy of an almost real curve. The new algorithm suits a larger family of curves, and enables the computation of braid monodromy not only of caspidal curves, but of general algebraic curves, with some non simple singularities. Moreover, it works also when in the fiber the curve admits any number of imaginary points. We also provide two examples of how to use the generalized algorithm.
\end{abstract}

\maketitle

\section*{Introduction}
Classification of surfaces is a major area of study in algebraic geometry. Hence, different methods for achieving this goal have grown to be very important. It was discovered that branch curves of hypersurfaces are invariant in some sense for those hypersurfaces, and the braid monodromy was suggested as an invariant for the branch curves.

The braid monodromy is an important tool in order to compute invariants of curves and surfaces, for example it is used as one way to compute the fundamental group of the complement of a curve. This invariant was used by Chisini \cite{chisini}, Kulikov \cite{Kul} and Kulikov-Teicher \cite{KuTe} in order to distinguish between connected components of the moduli space of surfaces of general type. 

Many previous works were done using the braid monodromy in order to classify surfaces, branch curves and in order to compute the fundamental groups of the complement of curves (for example see \cite{Suci,Garber1,Garber2} among many others). However, there is still a growing need of computed examples. Until now, most examples were computed on caspidal curves with only simple singularities and limited number of imaginary points. 

In this paper we improve the the algorithm given in \cite{BGTII} for the computation of braid monodromy of an almost real curve. The improved algorithm works on general algebraic curves, it enables to consider other singularities which were not considered in the past, where some of the singularities are not simple, and it also enables the computation with any number of complex points in any fiber above the real line.

The paper is organized as follows: In Section $1$ we give preliminaries and definitions. Section $2$ is devoted to the computations of
the local braid monodromy of the new types of singular points. In Section $3$ we present a new model for the fiber of the curve above a given point. We introduce the algorithm, which calculates the braid monodromy, and prove it's steps in Section $4$. Finally, in Section $5$ we give two examples 
for computations according to the algorithm.

\section{Preliminaries} \label{Preliminaries}

We assume that the reader is familiar with the braid group and the braid monodromy subjects (see for example \cite{Artin, Birman, Dehornoy, BGTII}); Therefore, we give here only notations and some of the basic definitions involved.

Let $S$ be a curve in $\C^2$. Denote by $\pi_1:S \to \C$ and by $\pi_2:S \to \C$ the projection to the first
and second coordinate respectively, defined in the obvious way. Let $n$ to be the degree of $S$ with respect to $\pi _1$. For $x \in \C$ we denote $K(x)$ the second coordinate of the points in $S$ which are in the fiber above $x$ (i.e., $K(x)=\pi_2(\pi_1^{-1}(x)$)).

Let $N \subset \C$ be the set $N=N(S)=\{x \in\C \ | \ |K(x)|<n\}=\{x_1, \cdots ,x_p\}$. Take $E$ to be a closed disc in $\C$ for which $N \subset E \setminus \partial E$.
In addition take $D$ to be a closed disc in $\C$ for which $D$ contains all the points $\{K(x) \ | \ x\in E\}$.
That means that when restricted to $E$, we have $S \subset E \times D$.

With these definitions we define the \emph{braid monodromy}

\bde
Let $M \in \partial E$ be the base point of $\pi _1(E \setminus N)$, and let $\sigma$ be an element of $\pi _1(E \setminus N)$. To $\sigma$ there are $n$ lifts in $S$, each begins and ends at the points of $M \times K(M)$. Projecting these lifts using $\pi _2:S \to \C$ we get $n$ paths in $D$ which begin and end at the points of $K(M)$. These paths induce a diffeomorphism of $\pi _1(D \setminus K(M))$ which is an element of the $\mathcal{MCG}$ of $\pi _1(D \setminus K(M))$ known to be the braid group $B_n$ \cite{Birman}. We call the homomorphism $\varphi :\pi _1(E \setminus N) \to B_n$ the \emph{braid monodromy} of $S$ with respect to $E \times D,\pi _1,$ and $M$.
\ede

\subsection{Real singular curves}

We call the curve $S$ a \emph{real singular curve} if it is defined by a polynomial with real coefficients, and if the set $N(S) \subset \RR$.
We restrict ourselves to work with subset of real singular curves as follows:

\bde
Let $B=B(S)=\{x \in \RR \ | \ \exists z_1,z_2 \in K(x) \setminus \RR \text{ s.t. } \Im(z_1)=\Im(z_2) \}$.
\ede

Note that if $z \in K(x) \setminus \RR$ so does it's complex conjugate $z'$ (because, $S$ is defined over $\RR$) so when $x \in B$ we have at least two pairs of points $z_1,z_2$ and $z'_1,z'_2$ such that $\Im(z_1)=\Im(z_2)$, $\Im(z'_1)=\Im(z'_2)$ and $\Im(z_1)=-\Im(z'_1)$.

\noindent
We work with curve $S$ which satisfies the following conditions:
\ben
\item
$S$ is real singular.
\item
For all $x \in N$ there is only one singular points of $S$ with $x$ as it's first coordinate.
\item
The set $B$ is finite.
\item
For all $x \in B$ the number of pairs $z_1,z_2 \in K(x) \setminus \RR$ s.t. $\Im(z_1)=\Im(z_2)$ is $2$.
\item
$B \cap N = \emptyset$.
\een

Our algorithm will apply to curves with the following types of singularities (looking at these singularities was motivated by the fact that we are interested mostly, but not only, in computing braid monodromy for branch curves of hypersurfaces).

\ben
\item
$a_1$: a branch point, topologically equivalent to $y^2-x=0$.
\item
$a_2$: a branch point, topologically equivalent to $y^2+x=0$.
\item
$b$: a tangent point, topologically equivalent to $y(y-x^2)=0$.
\item
$c$: intersection of $m$ non singular branches of $S$, transversal to each other.
\item
$d_1$: a cusp point, topologically equivalent to $y^2-x^\nu=0$, ($3 \leq \nu$).
\item
$d_2$ a cusp point, topologically equivalent to $y^2+x^\nu=0$, ($3 \leq \nu$).
\item
$d_3$: a cusp point, topologically equivalent to $y^3-x^2=0$.
\item
$e_1$: two cusps intersecting at one point, topologically equivalent to $(y^2+x^3)(y^2-x^3)=0$.
\item
$e_2$: two cusps intersecting at one point, topologically equivalent to $(y^3+ax^2)(y^3-bx^2)=0$ ($a\neq b$),($0<a,b$).
\item
$f_1$: a cusp point intersecting a line, topologically equivalent to $y(y^2-x^\nu)=0$.
\item
$f_2$: a cusp point intersecting a line, topologically equivalent to $y(y^2+x^\nu)=0$.
\item
$g_1$: a cusp point intersecting a line, topologically equivalent to $(y-ax)(y^2-x^\nu)=0$.
\item
$g_2$: a cusp point intersecting a line, topologically equivalent to $(y-ax)(y^2+x^\nu)=0$.
\item
$g_3$: a cusp point intersecting a line, topologically equivalent to $(y)(y^3-x^2)=0$.
\een

\brem $ $
\ben
\item
In our case topologically equivalent means that there is a deformation which transforms the given singular point in the curve $S$ into singularity types defined above, without changing the multiplicity of any point with respect to the projection $\pi_1$.
\item
Without loss of generality we assume that $a>0$ in the cases of singularity of type $g_1$ and $g_2$ since when $a<0$ all constructions and proofs are similar.
\item
The points of type $d_1$ and $d_2$ generalize the results previously obtained in \cite{BGTII} for the points $a_1$ and $a_2$. We keep the notations of the points $a_1$ and $a_2$ in order to comply with previous works.
\een
\erem

The computation of the braid monodromy of the curve is done using computations of local braid monodromy at each singular point and combining the results. So our aim now is to describe the local braid monodromy at each type of singularity defined above.

\section{Computation of local Braid Monodromy}

In this section we compute the local braid monodromy induced by each of the singularities defined above. We use the connection between the geometrical braid and it's description using Artin generators \cite{Artin,Birman}, by looking at the projection of the geometrical $3$-dimensional braid to the $xz$ plane keeping in mind which strings come above the other at intersection points of the projection. When two strings start or end at the same position we always rotate the braid clockwise by a very small angle in order to get a general position of the beginning and ending points of the strings. 

Local braid monodromy of singular points of types $a_1,a_2,b,c,d_1$ were discussed in \cite{BGTI}, \cite{BGTII} we will not repeat these computations here.

\subsection{Local braid monodromy of singularity type $d_2$}
The computation here is similar to the one presented in \cite{BGTI}[Proposition-Example VI.1.1]

\bpr \label{Prop d_2}
Let $S$ be the curve $y^2+x^\nu$. Then, the local braid monodromy at the origin is generated by the braid $\sigma _1^\nu$.
\epr

\begin{proof}
Take $E$ to be the unit disk, $D$ be a closed disk sufficiently large such that $S | _E \subset E \times D$.
$N=(0+0i) \in E$. We take $\Gamma=\{x=e^{2\pi i t} \ | \ t \in [0,1]\}$ to be the element of $\pi _1(E \setminus N,1)$. Thus, when lifting $\Gamma$ into $S$ we get two paths:
$$\delta_1(t)=(e^{2\pi i t},e^{\frac{\nu(2\pi i t+\pi i)}{2}})$$
$$\delta_2(t)=(e^{2\pi i t},-e^{\frac{\nu(2\pi i t+\pi i)}{2}})$$
Where, $t \in[0,1]$

Projecting them to $D$ we get two paths:
$$a_1(t)=e^{\frac{\nu(2\pi i t+\pi i)}{2}}$$
$$a_2(t)=-e^{\frac{\nu(2\pi i t+\pi i)}{2}}$$
Which induce the braid $\sigma _1^\nu$ which transpose the position of two strings $\nu$ times, as described in Figure \ref{braid d_2}.

\myfigure{
\begin{center}
\includegraphics[angle=90,scale=0.20]{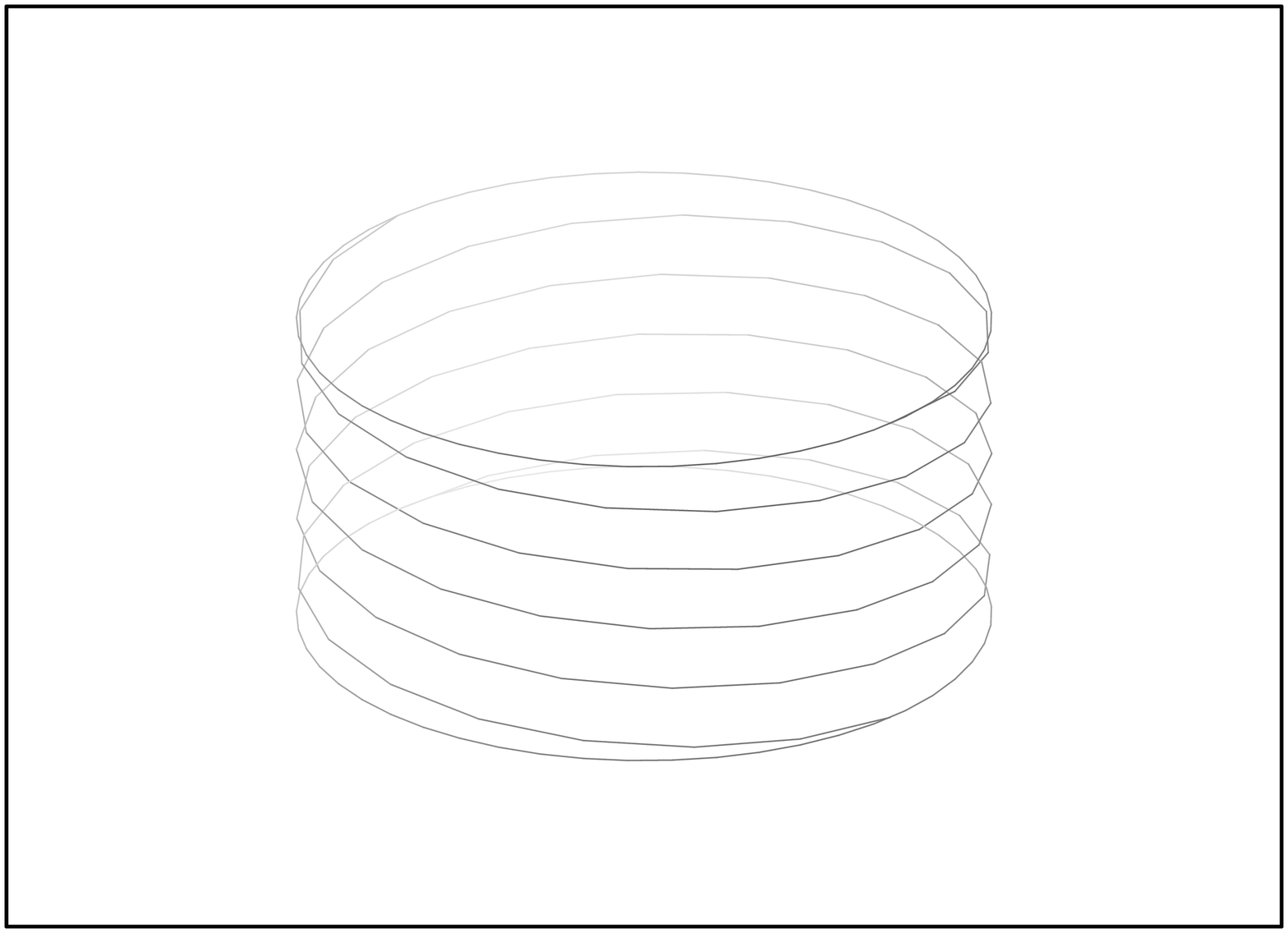}
\caption{The curve of type $d_2$ where $\nu=5$.}\label{braid d_2}
\end{center}
}

\end{proof}

\subsection{Local braid monodromy of singularity type $d_3$}

\bpr \label{Prop d_3}
Let $S$ be the curve $y^3-x^2$. Then, the local braid monodromy at the origin is generated by the braid $(\sigma _1 \sigma _2)^2$.
\epr

\begin{proof}
Take $E$ to be the unit disk, $D$ be a closed disk sufficiently large such that $S | _E \subset E \times D$.
$N=(0+0i) \in E$. We take $\Gamma=\{x=e^{2\pi i t} \ | \ t \in [0,1]\}$ to be the element of $\pi _1(E \setminus N,1)$. Thus, when lifting $\Gamma$ into $S$ we get $3$ paths:
$$\delta_k(t)=(e^{2\pi i t},e^{\frac{4}{3}\pi i t+\frac{2}{3}\pi ik})$$
Where, $t \in[0,1]$ and $k=0,1,2$.

Projecting them to $D$ we get $3$ paths:
$$a_k(t)=e^{\frac{4}{3}\pi i t+\frac{2}{3}\pi ik}$$
\noindent
Where $k=0,1,2$.

Which induce the braid $(\sigma _1 \sigma _2)^2$, see Figure \ref{braid d_3}.

\myfigure{
\begin{center}
\includegraphics[angle=90,scale=0.20]{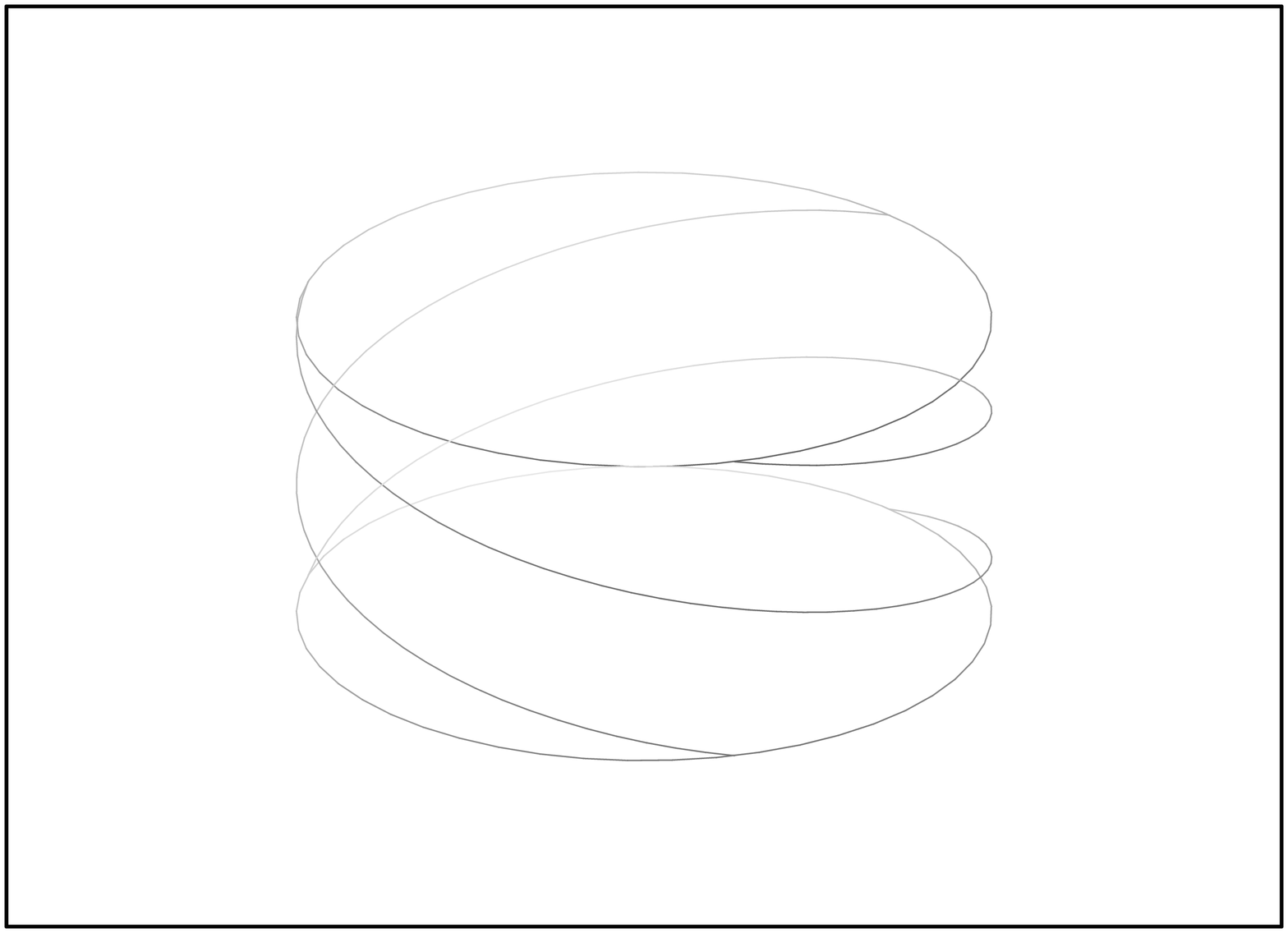}
\caption{The curve of type $d_3$.}\label{braid d_3}
\end{center}
}

\end{proof}

\subsection{Local braid monodromy of singularity type $e_1$}

\bpr \label{Prop e_1}
Let $S$ be the curve $(y^2+x^3)(y^2-x^3)$. Then, the local braid monodromy at the origin is generated the braid $(\sigma _2 \sigma _1 \sigma _3)^6$.
\epr

\begin{proof}
Take $E$ to be the unit disk, $D$ be a closed disk sufficiently large such that $S | _E \subset E \times D$.
$N=(0+0i) \in E$. We take $\Gamma=\{x=e^{2\pi i t} \ | \ t \in [0,1]\}$ to be the element of $\pi _1(E \setminus N,1)$. Thus, when lifting $\Gamma$ into $S$ we get $4$ paths:
$$\delta_1(t)=(e^{2\pi i t},e^{3\pi i t})$$
$$\delta_2(t)=(e^{2\pi i t},-e^{3\pi i t})$$
$$\delta_2(t)=(e^{2\pi i t},e^{3\pi i t+\frac{\pi}{2} i})$$
$$\delta_4(t)=(e^{2\pi i t},-e^{3\pi i t+\frac{\pi}{2} i})$$
Where, $t \in[0,1]$

Projecting them to $D$ we get two paths:
$$a_1(t)=e^{3\pi i t}$$
$$a_2(t)=-e^{3\pi i t}$$
$$a_2(t)=e^{3\pi i t+\frac{\pi}{2} i}$$
$$a_4(t)=-e^{3\pi i t+\frac{\pi}{2} i}$$
Which induce the braid $(\sigma _2 \sigma _1 \sigma _3)^6$, as can be seen in Figure \ref{braid e_1}.

\myfigure{
\begin{center}
\includegraphics[angle=90,scale=0.20]{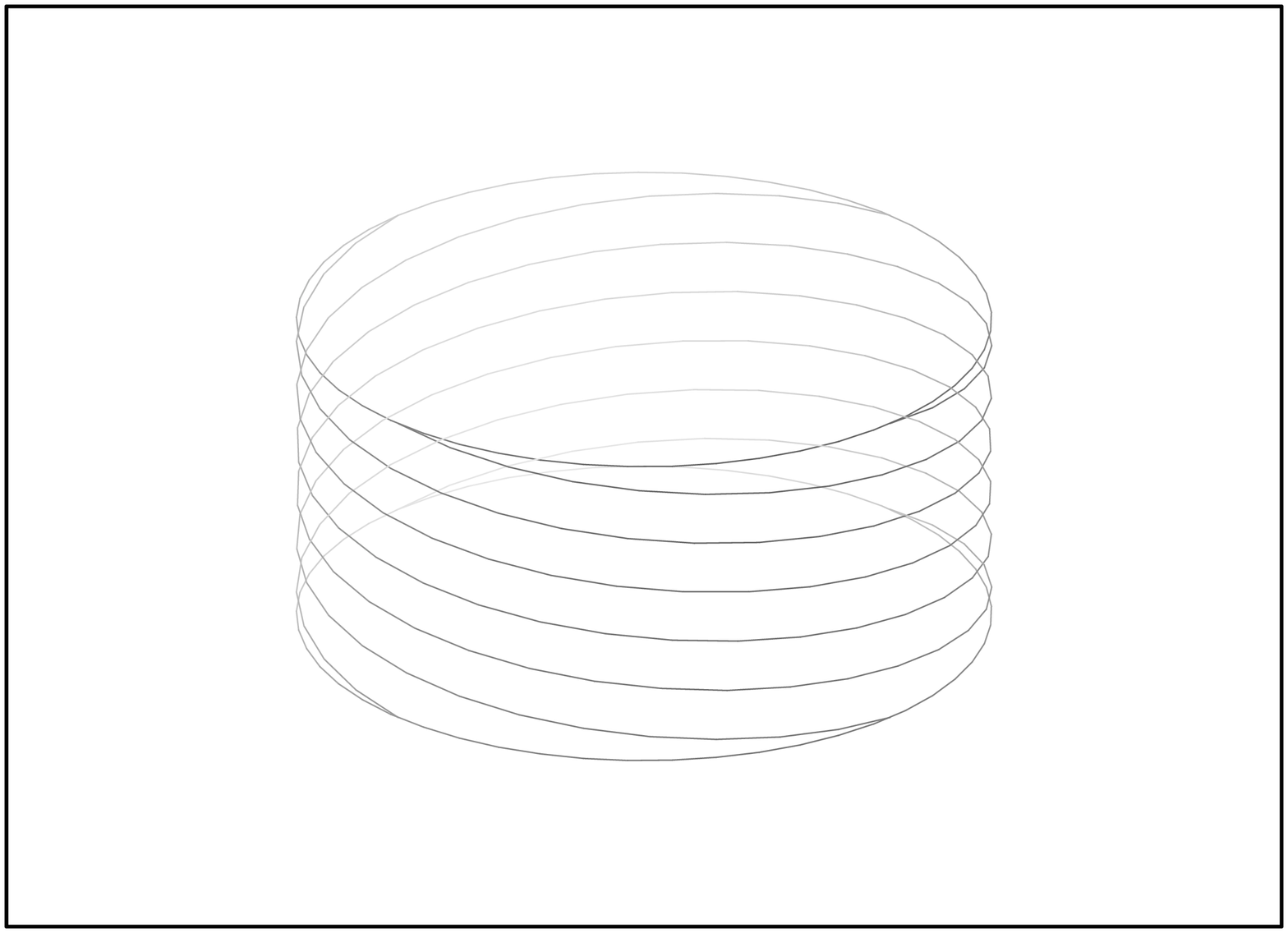}
\caption{The curve of type $e_1$.}\label{braid e_1}
\end{center}
}

\end{proof}

\subsection{Local braid monodromy of singularity type $e_2$}
\bpr \label{Prop e_2}
Let $S$ be the curve $(y^3+x^2)(y^3-x^2)$. Then, the local braid monodromy at the origin is generated by the braid $(\sigma _2 \sigma _4 \sigma _1 \sigma _3 \sigma _5)^4$.
\epr

\begin{proof}
Take $E$ to be the unit disk, $D$ be a closed disk sufficiently large such that $S | _E \subset E \times D$.
$N=(0+0i) \in E$. We take $\Gamma=\{x=e^{2\pi i t} \ | \ t \in [0,1]\}$ to be the element of $\pi _1(E \setminus N,1)$. Thus, when lifting $\Gamma$ into $S$ we get $6$ paths:
$$\delta_{1,k}(t)=(e^{2\pi i t},e^{\frac{4 \pi i t}{3}+\frac{2 \pi i k}{3}})$$
$$\delta_{2,k}(t)=(e^{2\pi i t},e^{\frac{4 \pi i t}{3}+\frac{\pi}{3}+\frac{2 \pi i k}{3}})$$
Where, $t \in[0,1]$, and $k=1, \cdots ,3$.

Projecting them to $D$ we get two paths:
$$a_{1,k}(t)=e^{\frac{4 \pi i t}{3}+\frac{2 \pi i k}{3}}$$
$$a_{1,k}(t)=e^{\frac{4 \pi i t}{3}+\frac{\pi}{3}+\frac{2 \pi i k}{3}}$$
Which induce the braid $(\sigma _2 \sigma _4 \sigma _1 \sigma _3 \sigma _5)^4$, as can be seen in Figure \ref{braid e_2}.

\myfigure{
\begin{center}
\includegraphics[angle=90,scale=0.20]{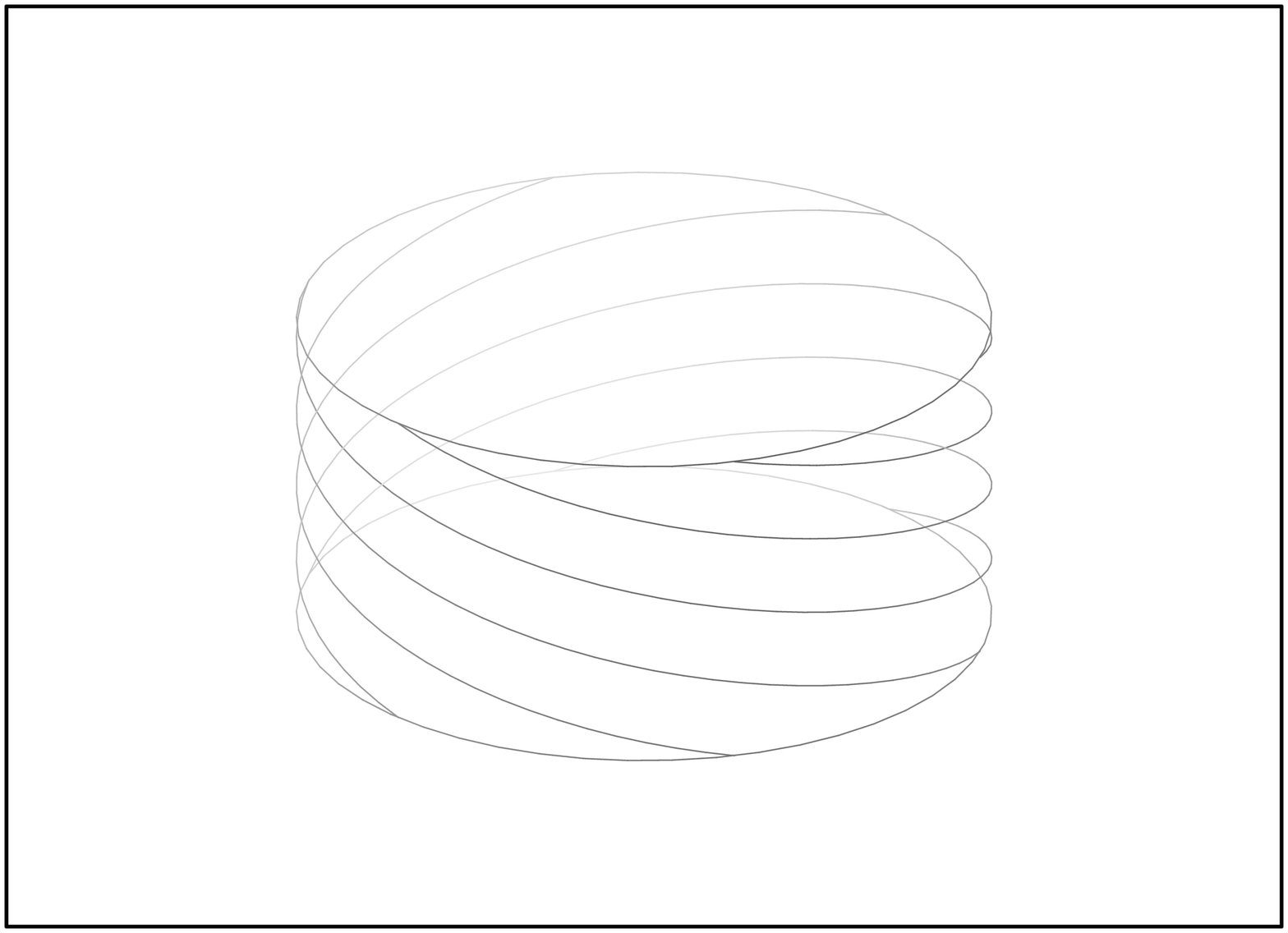}
\caption{The curve of type $e_2$.}\label{braid e_2}
\end{center}
}

\end{proof}

\subsection{Local braid monodromy of singularity type $f_1$}

\bpr \label{Prop f_1}
Let $S$ be the curve $(y^2-x^\nu)y$. Then, the local braid monodromy at the origin is generated by the braid $(\sigma _2 \sigma _1 \sigma _2)^\nu$.
\epr

\begin{proof}
Take $E$ to be the unit disk, $D$ be a closed disk sufficiently large such that $S | _E \subset E \times D$.
$N=(0+0i) \in E$. We take $\Gamma=\{x=e^{2\pi i t} \ | \ t \in [0,1]\}$ to be the element of $\pi _1(E \setminus N,1)$. Thus, when lifting $\Gamma$ into $S$ we get Three paths:
$$\delta_1(t)=(e^{2\pi i t},e^{\nu \pi i t})$$
$$\delta_2(t)=(e^{2\pi i t},-e^{\nu \pi i t})$$
$$\delta_3(t)=(e^{2 \pi i t},0)$$
Where, $t \in[0,1]$.

Projecting them to $D$ we get three paths:
$$a_1(t)=e^{\nu \pi i t}$$
$$a_2(t)=-e^{\nu \pi i t}$$
$$a_3(t)=0$$
Which induce the braid $(\sigma _2 \sigma _1 \sigma _2)^\nu$, as can be seen in Figure \ref{braid f_1}.

\myfigure{
\begin{center}
\includegraphics[angle=90,scale=0.20]{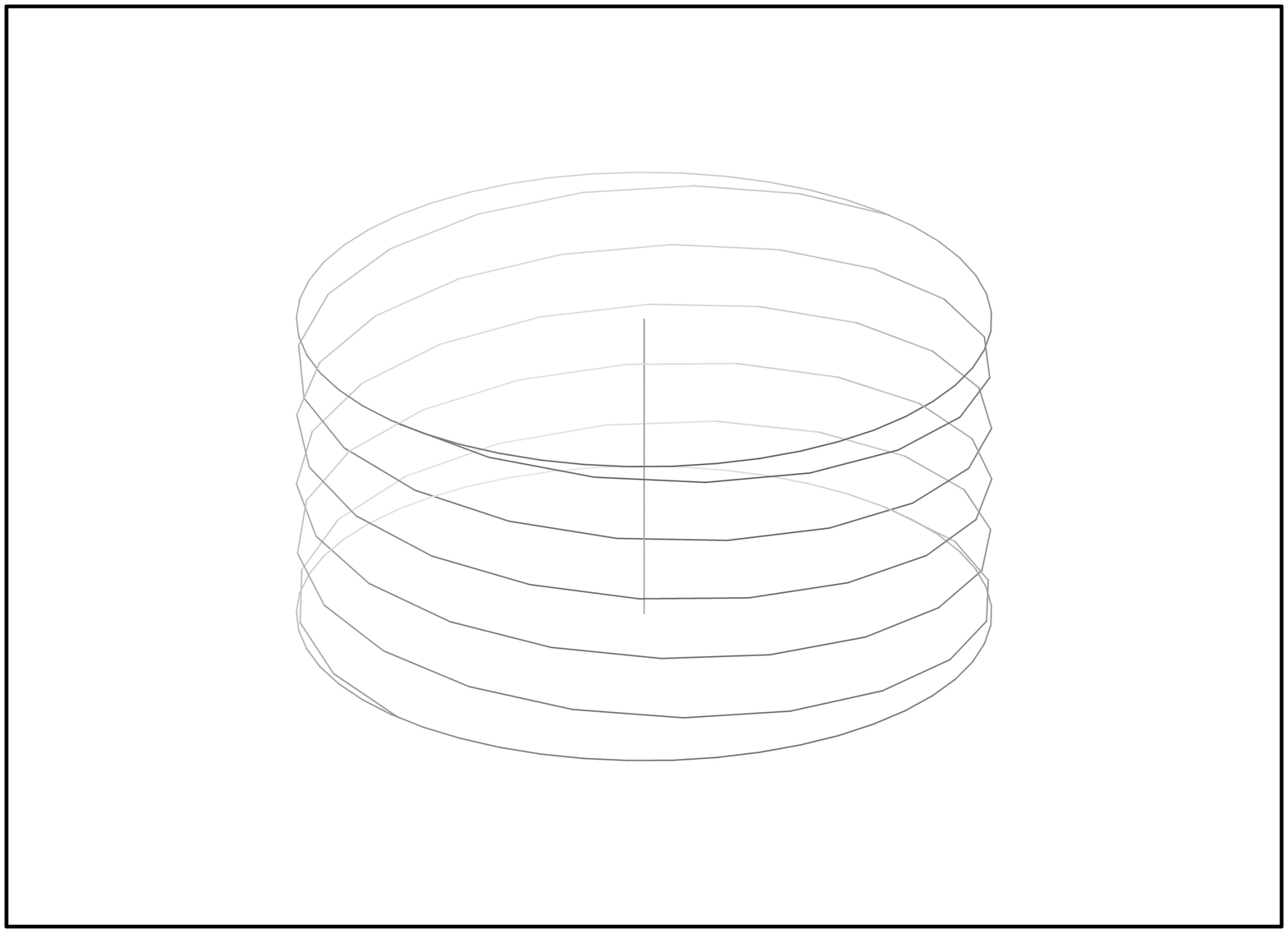}
\caption{The curve of type $f_1$ where $\nu=5$.}\label{braid f_1}
\end{center}
}

\end{proof}

\subsection{Local braid monodromy of singularity type $f_2$}

\bpr
Let $S$ be the curve $(y^2+x^\nu)y$. Then, the local braid monodromy at the origin is generated by the braid $(\sigma _2 \sigma _1 \sigma _2)^\nu$.
\epr

\begin{proof}
The proof is similar to the proof of Proposition \ref{Prop f_1}.
\end{proof}

\subsection{Local braid monodromy of singularity type $g_1$}

\bpr \label{Prop g_1}
Let $S$ be the curve $(y^2-x^\nu)(y-ax)$. Then, the local braid monodromy at the origin is generated by the braid $(\sigma _2 \sigma _1 \sigma _2^\nu \sigma _1 \sigma _2)$.
\epr

\begin{proof}
Without loos of generality we may choose $a=2$, this will help us to avoid multiple intersections between the line and the cusp inside the unit disk. Take $E$ to be the unit disk, $D$ be a closed disk sufficiently large such that $S | _E \subset E \times D$ (in case where $a \neq 2$ we may need to take $E$ to be a smaller disk).
$N=(0+0i) \in E$. We take $\Gamma=\{x=e^{2\pi i t} \ | \ t \in [0,1]\}$ to be the element of $\pi _1(E \setminus N,1)$. Thus, when lifting $\Gamma$ into $S$ we get Three paths:
$$\delta_1(t)=(e^{2 \pi i t},e^{\nu \pi i t})$$
$$\delta_2(t)=(e^{2 \pi i t},-e^{\nu \pi i t})$$
$$\delta_3(t)=(e^{2 \pi i t},2e^{2 \pi i t})$$
Where, $t \in[0,1]$.

Projecting them to $D$ we get three paths:
$$a_1(t)=e^{\nu \pi i t}$$
$$a_2(t)=-e^{\nu \pi i t}$$
$$a_3(t)=2e^{2 \pi i t}$$
Which induce the braid $(\sigma _2 \sigma _1 \sigma _2^\nu \sigma _1 \sigma _2)$, as can be seen in Figure \ref{braid g_1}.

\myfigure{
\begin{center}
\includegraphics[angle=90,scale=0.20]{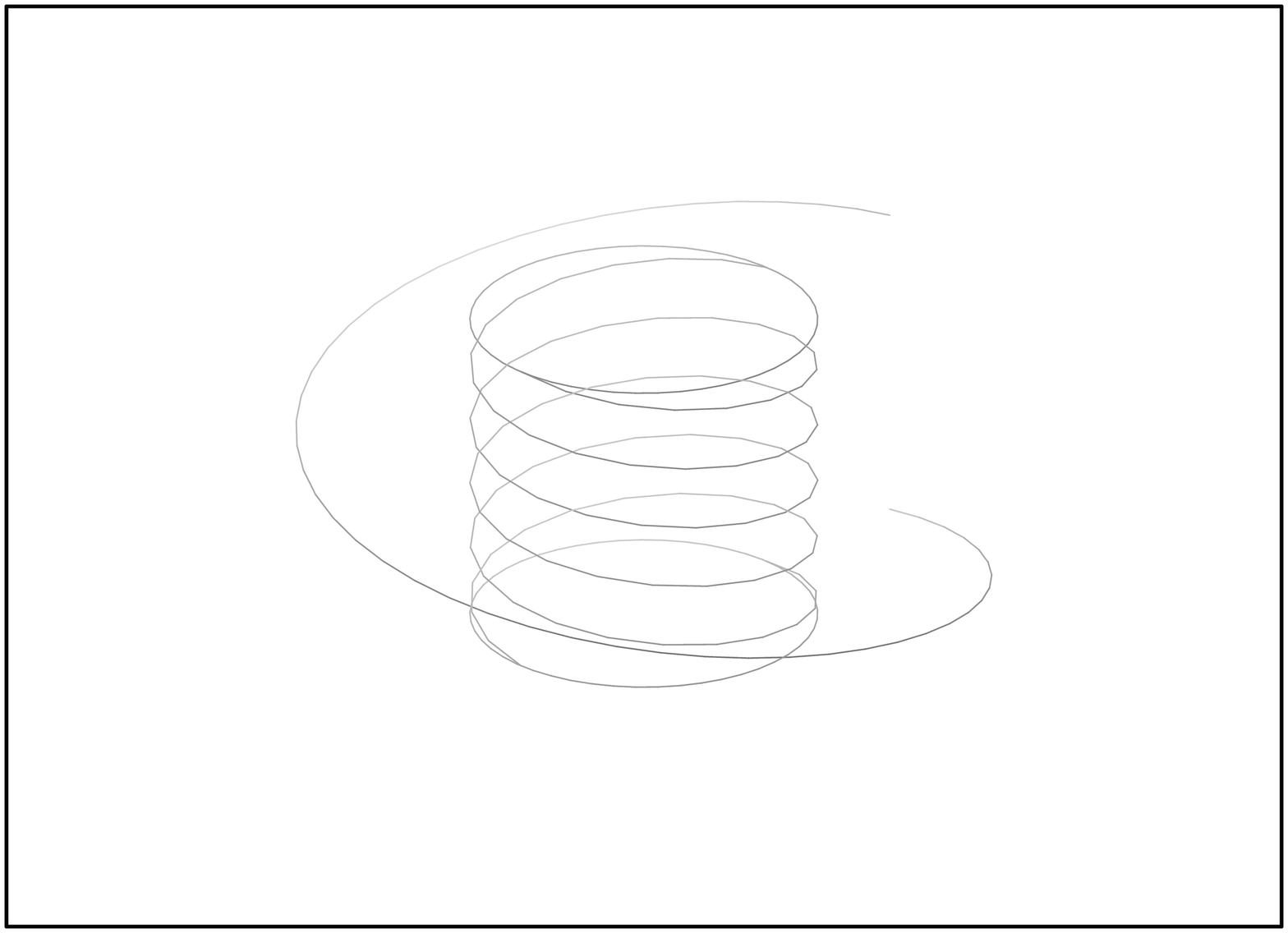}
\caption{The curve of type $g_1$ where $\nu=5$.}\label{braid g_1}
\end{center}
}

\end{proof}

\subsection{Local braid monodromy of singularity type $g_2$}

\bpr \label{Prop g_2}
Let $S$ be the curve $(y^2+x^\nu)(y-ax)$. Then, the local braid monodromy at the origin is generated by the braid $(\sigma _2 \sigma _1 \sigma _2^\nu \sigma _1 \sigma _2$).
\epr

\begin{proof}
The proof is similar to the proof of Proposition \ref{Prop g_1}.
\end{proof}

\subsection{Local braid monodromy of singularity type $g_3$}

\bpr \label{Prop g_3}
Let $S$ be the curve $(y)(y^3-x^2)$. Then, the local braid monodromy at the origin is generated by the braid $(\sigma _1 \sigma _2 \sigma _3 \sigma _2)^2$.
\epr

\begin{proof}
Take $E$ to be the unit disk, $D$ be a closed disk sufficiently large such that $S | _E \subset E \times D$.
$N=(0+0i) \in E$. We take $\Gamma=\{x=e^{2\pi i t} \ | \ t \in [0,1]\}$ to be the element of $\pi _1(E \setminus N,1)$. Thus, when lifting $\Gamma$ into $S$ we get $4$ paths:
$$\delta_k(t)=(e^{2\pi i t},e^{\frac{4}{3}\pi i t+\frac{2}{3}\pi ik})$$
$$\delta_4(t)=(e^{2 \pi i t},0)$$
Where, $t \in[0,1]$ and Where $k=0,1,2$.

Projecting them to $D$ we get two paths:
$$a_k(t)=e^{\frac{4}{3}\pi i t+\frac{2}{3}\pi ik}$$
$$a_4(t)=0$$
\noindent
Where $k=0,1,2$.
Which induce the braid $(\sigma _1 \sigma _2 \sigma _3 \sigma _2)^2$, as can be seen in Figure \ref{braid g_3}.

\myfigure{
\begin{center}
\includegraphics[angle=90,scale=0.20]{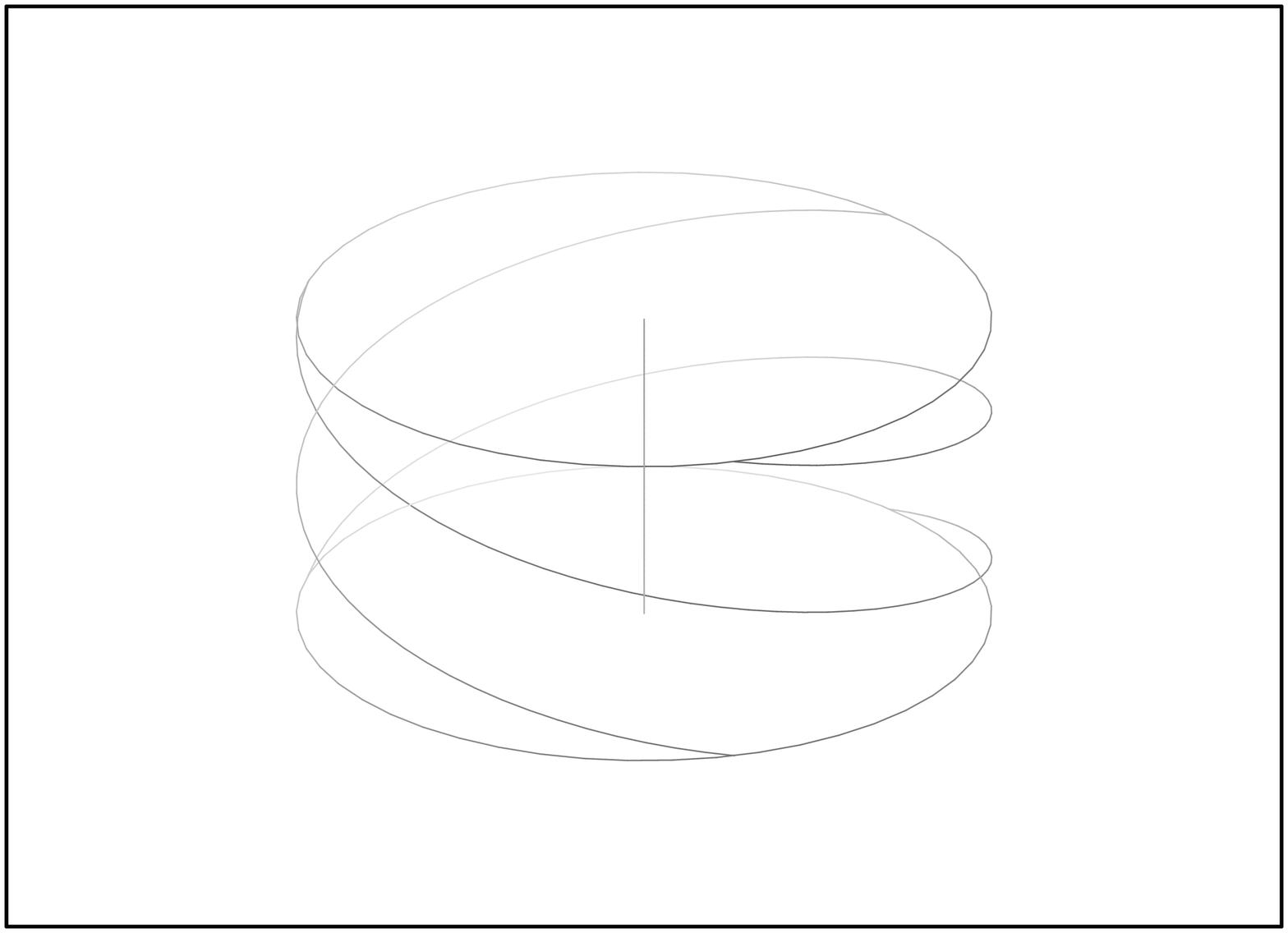}
\caption{The curve of type $g_3$.}\label{braid g_3}
\end{center}
}
\end{proof}

%%%%%%%%%%%%%%%%%%%%%%%%%%%%%%%%%%%%%%%%%%%%%%%%%%%%%%%%%%%%%%%%%%%%%%%%%%%%%%%%%%%%%%%%%%%%%%%%
\section{A model for the fiber} \label{A model for the fiber}

\subsection{Model definition}
In the algorithm we use a model of the fiber above a real point $x_0$.
In this fiber we distinguish $n$ (the degree of the projection to the first coordinate $\pi _1$) real and complex points.

Since the model should be diffeomorphic to $D$ we take the model to be a disk $\widetilde{D} \subset \C$ which is centered at the point $n+1$ and have radius $n+1$.
The distinguished real points will be in the set $K'=\{1, \cdots, n\}$, and the distinguished complex points will be in the set $K''=\{n+1 \pm i, \cdots ,n+1 \pm \lceil \frac{n}{2} \rceil i\}$.

Usually, we shell assume that we have $d$ (an even number) complex points and therefore, we identify in $\widetilde{D}$ the set of points $K_d=\{1, \cdots , (n-d),(n+1)\pm i, \cdots ,(n+1)\pm \frac{d}{2}i\}$.

For example, if $n=9$ and $d=4$ (i.e., we have $4$ complex points), we get $K_4=\{1,2,3,4,5,10+i,10+2i,10-i,10-2i\}$ and the model is $\widetilde{D}$ with the set of points $K_4$ identified, as can be seen in figure \ref{Model example}

\myfigure{
\begin{center}
\includegraphics[scale=0.75]{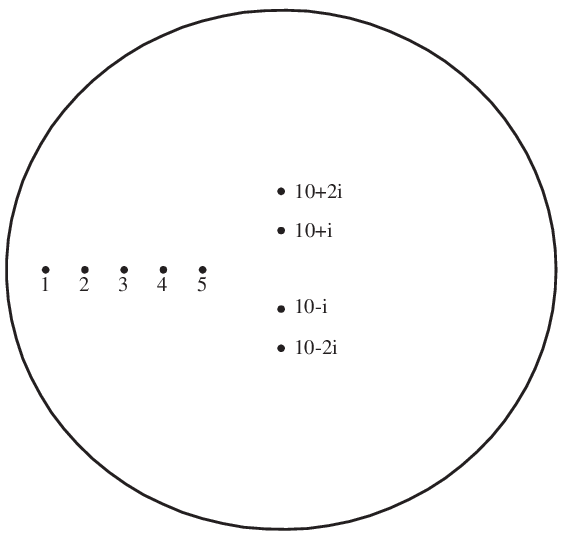}
\caption{$E_4$}\label{Model example}
\end{center}
}

\subsection{Half-twists}
We need to define some braids using our model, hence we need to describe an element of the $\mathcal{MCG}$ of the disk $\widetilde{D}$ punctured at the points $K_d$.

Let $\sigma$ be a simple path in $\widetilde{D}$ connecting two points of $K_d$ which do not intersect any point of $K_d$ besides its two end points. Take a small neighborhood $U$ of $\sigma$ and a small neighborhood $U'$ of $U$, which contain from $K_d$ only the two end points of $\sigma$. We define a diffeomorphism which switchs the two end points of $\sigma$ along the curve counterclockwise, and is the identity outside $U'$. We call the diffeomorphism constructed the \emph{(positive) half-twist} defined by $\sigma$, and we denote it by $H(\sigma)$. Note, different diffeomorphisms defined on the same curve $\sigma$, using the above construction, will induce the same braid, hence half-twists are well defined.

\bde
Let $\sigma$ be a non self intersecting path in $\widetilde{D} \setminus \partial \widetilde{D}$ starting and ending at two different points of $K_d$, which go through all the points of $K_d$. The points of $K_d$ induce a partition of $\sigma$ into $n-1$ consecutive paths $\sigma _1,\cdots ,\sigma _{n-1}$. We call the ordered system of half-twists $H(\sigma _1),\cdots ,H(\sigma _{n-1})$ a \emph{frame} of $B_n$.
\ede

It is well known that if $(H(\sigma _1),...,H(\sigma _{n-1}))$ is a frame of $B_n$, then the braid group is generated by its elements, and the relations are just like the relations in Artin presentation for the braid group \cite{BGTI}.

\subsection{Diffeomorphisms in the model}

We will need to define braids using the model $\widetilde{D}$ of the punctured disk $D \setminus K(x)$. For this matter we will give here a few diffeomorphisms which we will use as building blocks. Note that all diffeomorphisms are the identity on $\partial \widetilde{D}$. Moreover, instead of giving tedious formal definition, we will describe a continuous family of diffeomorphisms (starting with the identity), by the action they induce on a closed or a finite set of points in the interior of $\widetilde{D}$. The braid we define will be the last diffeomorphism in the family.

\bde \label{R}
Let $R^n\left<s,r \right>$ be the diffeomorphism which rotates counterclockwise in an angle of $2\pi n$ a closed disk centered at $s$ and of radius $r$.
\ede

\bde
$L\left<c,k\right>$ is the diffeomorphism which maps the point $c \in \C$ to the point $k\in \C$ and the point $\bar{c}$ (the complex conjugate of $c$) to the point $\bar{k}$ (the complex conjugate of $k$), along a simple path (usually a straight line) which does not intersect $K''$.
\ede

\bde
$M\left<(h,j),(k,l)\right>$ is the diffeomorphism which maps all points in the real segment $[h,j]$ onto the real segment $[k,l]$ (note that $j-h=l-k$).
\ede

\bde
$C$ is the diffeomorphism which maps the points along the straight line between $n+1+i$ to $n+1+\frac{d}{2}i$ onto the straight line between $n+1+2i$ and $n+1+(\frac{d}{2}+1)i$, and the straight line between $n+1-i$ to $n+1-\frac{d}{2}i$ onto the straight line between $n+1-2i$ and $n+1-(\frac{d}{2}+1)i$.
\ede

\bde \label{A}
$A^n\left<c,r_1,r_2\right>$ is the diffeomorphism which rotates counterclockwise by angle of $2 \pi n$ the annulus centered at $c$ between the radiuses $r_1$ and $r_2$.
\ede

We now turn to define several diffeomorphisms, each one of them in connection with a specific type of point defined in Section \ref{Preliminaries}. Note, that since usually $\Delta$ denotes a half-twist in the braid group the powers of the $\Delta$ will be twice the angles of rotations used in Definitions \ref{R} and \ref{A}.
\begin{small}

\bform
\Delta _{a_1}^{\frac{1}{2}}\left<k,l\right>=
L\left<n+1+i,k+\frac{1}{2}+\frac{i}{2}\right> \circ
M\left<(k,n-d),(k+2,n-d+2)\right> \circ 
R^{\frac{1}{4}}\left<k+\frac{1}{2},\frac{1}{2}\right> \circ
C^{-1}
\eform

\bform
\Delta _{a_2}^{\frac{1}{2}}\left<k,l\right>=
R^{\frac{1}{4}}\left<k+\frac{1}{2},\frac{1}{2}\right> \circ
C \circ
L\left<k+\frac{1}{2}+\frac{i}{2},n+1+i\right> \circ 
M\left<(k+2,n-d),(k,n-d-2)\right>
\eform

\bform
\Delta _b^2\left<k,l\right> = R^1\left<k+\frac{1}{2},\frac{1}{2}\right>
\eform

\bform
\Delta _c\left<k,l\right> = R^{\frac{1}{2}}\left<\frac{k+l}{2},\frac{l-k}{2}\right>
\eform

\noindent
$
\Delta _{d_1}^{\frac{\nu}{2}}\left<k\right>=
\left \{
\begin{array}{ll}
L\left<n+1+i,k+\frac{1}{2}+\frac{i}{2}\right> \circ 
M\left<(k,n-d),(k+2,n-d+2)\right> \circ & \text{if } \nu \text{ is odd} \\
R^{\frac{\nu}{4}}\left<k+\frac{1}{2},\frac{1}{2}\right> \circ C^{-1} & \\ \\
R^{\frac{\nu}{4}}\left<k+\frac{1}{2},\frac{1}{2}\right> & \text{if }\nu \text{ is even}
\end{array}
\right.
$

\noindent
$
\Delta _{d_2}^{\frac{\nu}{2}}\left<k,l\right>=
\left \{
\begin{array}{ll}
R^{\frac{\nu}{4}}\left<k+\frac{1}{2},\frac{1}{2}\right> \circ C \circ
L\left<k+\frac{1}{2}+\frac{i}{2},n+1+i\right> \circ & \text{if } \nu \text{ is odd} \\
M\left<(k+2,n-d),(k,n-d-2)\right> & \\ \\
L\left<n+1+i,k+\frac{1}{2}+\frac{i}{2}\right> \circ 
M\left<(k,n-d),(k+2,n-d+2)\right> \circ & \text{if } \nu \text{ is even} \\
R^{\frac{\nu}{4}}\left<k+\frac{1}{2},\frac{1}{2}\right> \circ & \\
M\left<(k+2,n-d+2),(k,n-d)\right> \circ L\left<k+\frac{1}{2}+\frac{i}{2},n+1+i\right> &
\end{array}
\right.
$

\bform
\Delta _{d_3}^{\frac{2}{3}}\left<k\right>=
L\left<n+1+i,k+\frac{1}{4}+\frac{\sqrt{3}}{4}i\right> \circ
M\left<(k,n-d),(k+1,n-d+1)\right> \circ \\ 
R^\frac{1}{3}\left<k+\frac{1}{2},\frac{1}{2}\right> \circ 
M\left<(k+1,n-d+1),(k,n-d)\right> \circ
L\left<k+\frac{1}{4}+\frac{\sqrt{3}}{4}i,n+1+i\right> 
\eform

\bform
\Delta _{e_1}^{\frac{3}{2}}\left<k,l\right>=
L\left<n+1+i,k+\frac{1}{2}+\frac{i}{2}\right> \circ
R^{\frac{3}{4}}\left<k+\frac{1}{2},\frac{1}{2}\right> \circ \\
L\left<k+\frac{1}{2}+\frac{i}{2},n+1+i\right> 
\eform

\bform
\Delta _{e_2}^{\frac{2}{3}}\left<k,l\right>_{\<a<b\>}=
L\left<n+1+i,l-\frac{1}{4}+\frac{\sqrt{3}}{4}i\right> \circ
L\left<n+1+2i,k+\frac{1}{4}+\frac{\sqrt{3}}{4}i\right> \circ
R^{\frac{1}{3}}\left<k+\frac{1}{2},\frac{1}{2}\right> \circ
L\left<k+\frac{1}{4}+\frac{\sqrt{3}}{4}i,n+1+2i\right> \circ
L\left<l-\frac{1}{4}+\frac{\sqrt{3}}{4}i,n+1+i\right> 
\eform

\bform
\Delta _{e_2}^{\frac{2}{3}}\left<k,l\right>_{\<a>b\>}=
L\left<n+1+i,k+\frac{1}{4}+\frac{\sqrt{3}}{4}i\right> \circ
L\left<n+1+2i,l-\frac{1}{4}+\frac{\sqrt{3}}{4}i\right> \circ
R^{\frac{1}{3}}\left<k+\frac{1}{2},\frac{1}{2}\right> \circ
L\left<l-\frac{1}{4}+\frac{\sqrt{3}}{4}i,n+1+2i\right> \circ
L\left<k+\frac{1}{4}+\frac{\sqrt{3}}{4}i,n+1+i\right> 
\eform

For the next diffeomorphism we need to define $\mathcal{D}$ as follows:
\bform
\mathcal{D}=L\left<n+1+i,k+1+i\right> \circ
M\left<(k+1,n-d),(k+3,n-d+2)\right> \circ \\
M\left<(k,k),(k+1,k+1)\right>
\eform

\noindent
$\Delta _{f_1}^{\frac{\nu}{2}}\left<k,l\right>= \left \{
\begin{array}{ll}
\mathcal{D} \circ R^{\frac{\nu}{4}}\left<k+1,1\right> \circ C^{-1} & \text{if } \nu \text{ is odd} \\ \\
R^{\frac{\nu}{4}}\left<k+1,1\right> & \text{if } \nu \text{ is even } 
\end{array}
\right.
$

\noindent
$\Delta _{f_2}^{\frac{\nu}{2}}\left<k,l\right>= \left \{
\begin{array}{ll}
R^{\frac{\nu}{4}}\left<k+1,1\right> \circ
C \circ
L\left<k+1+i,n+1+i\right> \circ & \text{if } \nu \text{ is odd} \\
M\left<(k+1,k+1),(k,k)\right> \circ & \\
M\left<(k+3,n-d),(k+1,n-d-2)\right>  & \\ \\
\mathcal{D} \circ R^{\frac{\nu}{4}}\left<k+1,1\right> \circ \mathcal{D}^{-1} & \text{if } \nu \text{ is even} 
\end{array}
\right.
$

\noindent
$\Delta _{g_1}^{\frac{\nu}{2}}\left<k,l\right>= \left \{
\begin{array}{ll}
M\left<(k+1,n-d),(k+3,n-d+2)\right> \circ & \text{if } \nu \text{ is odd} \\
L\left<n+1+i,k+1+\frac{i}{2}\right> \circ & \\
A^{\frac{1}{2}}\left<k+1,\frac{3}{4},\frac{5}{4}\right> \circ
R^{\frac{\nu}{4}}\left<k+1,\frac{1}{2}\right> \circ & \\
M\left<(k-\frac{1}{2},k+\frac{1}{2}),(k,k+1)\right> \circ C^{-1} & \\ \\
M\left<(k,k),(k+\frac{1}{2},k+\frac{1}{2})\right> \circ
A^{\frac{1}{2}}\left<k+\frac{3}{2},\frac{3}{4},\frac{5}{4}\right> \circ & \text{if } \nu \text{ is even}\\
R^{\frac{\nu}{4}}\left<k+\frac{3}{2},\frac{1}{2}\right> \circ
M\left<(k+1,k+2),(k,k+1)\right> \circ  & \\
M\left<(k+\frac{5}{2},k+\frac{5}{2}),(k+2,k+2)\right> &
\end{array}
\right.
$

\noindent
$\Delta _{g_2}^{\frac{\nu}{2}}\left<k,l\right>= \left \{
\begin{array}{ll}
M\left<(k,k),(k+\frac{1}{2},k+\frac{1}{2})\right> \circ
A^{\frac{1}{2}}\left<k+\frac{3}{2},\frac{3}{4},\frac{5}{4}\right> \circ & \text{if } \nu \text{ is odd} \\
R^{\frac{\nu}{4}}\left<k+\frac{3}{2},\frac{1}{2}\right> \circ C \circ
L\left<k+\frac{3}{2}+\frac{i}{2},n+1+i\right> \circ & \\
M\left<(k+\frac{5}{2},k+\frac{5}{2}),(k,k)\right> \circ  \\
M\left<(k+3,n-d),(k+1,n-d-2)\right> & \\ \\
M\left<(k+1,n-d),(k+3,n-d+2)\right> \circ & \text{if } \nu \text{ is even} \\
L\left<n+1+i,k+1+\frac{i}{2}\right> \circ \\
A^{\frac{1}{2}}\left<k+1,\frac{3}{4},\frac{5}{4}\right> \circ
R^{\frac{\nu}{4}}\left<k+1,\frac{1}{2}\right> \circ & \\
L\left<k+1+\frac{i}{2},n+1+i\right> \circ & \\
M\left<(k+2,n-d+2),(k,n-d)\right> & \\
\end{array}
\right.
$

\bform
\Delta _{g_3}^{\frac{2}{3}}\left<k,l\right>=
M\left<(k+1,k+1),(k+\frac{1}{2},k+\frac{1}{2})\right> \circ
L\left<n+1+i,k-\frac{1}{4}+\frac{\sqrt{3}}{4}i\right> \circ \\
R^{\frac{1}{3}}\left<k,\frac{1}{2}\right> \circ
L\left<k-\frac{1}{4}+\frac{\sqrt{3}}{4}i,n+1+i\right> \circ \\
M\left<(k+\frac{1}{2},k+\frac{1}{2}),(k+1,k+1)\right>
\eform
\end{small}

\bform
\Delta _{com}^1\<k,l\>=
R^{\frac{1}{2}}\<n+1+(k+\frac{1}{2})i,\frac{1}{2}\> \circ
R^{-\frac{1}{2}}\<n+1-(k+\frac{1}{2})i,\frac{1}{2}\>
\eform

\section{Description of the algorithm} \label{Describing the algorithm}

In this section we develop the theory needed for the algorithm. We prove the theorems which form the heart of the algorithm, and at the end of this section we give the algorithm's steps.

\subsection{definitions}

\bde \label{Lefchetz pair}
Let $0 < \alpha <\epsilon$ be two very small real numbers. Let $x_0 \in N$ and $A_0=(x_0,y_0)$ the singular point of the projection map $\pi _1$ above $x_0$. Let $x_0' \in E \cap \RR$ be a point such that $|x_0-x_0'|<\alpha$. Then it is possible to enumerate the points of $K(x_0') \cap \RR$ and  get that $K(x_0') \cap \RR=\{y_1< \cdots <y_n\}$. For sufficiently small $\epsilon$ there exists a unique pair of numbers $k$ and $l$ ($1 \leq k \leq l \leq n$) such that $K(x_0') \cap (y_0-\epsilon,y_0+\epsilon)=\{y_k<\cdots <y_l\}$. We call $(k,l)$ the \emph{Lefschetz pair} associated with the point $x_0$ and denote it by $\mathcal{L}(x_0)$.
\ede

\brem
Definition \ref{Lefchetz pair} does not hold when $A_0$ is of type $d_2,f_2$ and $\nu$ is even. In these cases we enumerate the real points of $K(x_0)=\{y_1 < \cdots <y_n\}$, and take $k=l=j$ where $y_j$ is the second coordinate of $A_0$.
\erem

\brem
If $A_0$ is a point of type
\ben
\item
$a_1,b,c,d_1,d_3,e_1,e_2,f_1,g_1,g_3$ we take $x_0'>x_0$.
\item
$a_2,d_2 (\text{odd } \nu),f_2 (\text{odd } \nu),g_2$ we take $x_0'<x_0$.
\een
\erem

\brem \label{k=l+1}
If $A_0$ is a point of type
\ben
\item
$a_1,a_2,b,d_1,d_2 (\text{odd } \nu),e_1,e_2,g_3$ we have $l=k+1$.
\item
$f_1,f_2 (\text{odd } \nu),g_1,g_2 (\text{odd } \nu)$ we have $l=k+2$.
\item
$c$ we have $l=k+n_l-1$ where $n_l$ is the number of intersecting (locally) non singular branches.
\item
$d_3,d_2 (\text{even } \nu),f_2 (\text{even } \nu),g_2 (\text{even } \nu)$ we have $l=k$.
\een
\erem

\bde \label{complex lefschetz pair}
Let $x_0 \in B$, and let $(x_0,y_0),(x_0,y_0')\in S$ be the pair of points for which $0<\Im(y_0)=\Im(y_0')$. Take $x_0' \in E \cap \RR$ to be a point such that $|x_0'-x_0|<\delta$ for a very small $0<\delta$. Then it is possible to enumerate by the $\Im$ coordinate the points of $K(x_0') \setminus \RR$ and get $K(x_0') \setminus \RR=\{y_1< \cdots <y_d\}$. Since the curve $S$ is defined over $\RR$ we look at the set $A=\{\Im(y_1),\cdots ,\Im(y_d)\}\cap \RR^+=\{r_1, \cdots, r_\frac{d}{2}\}=\{\Im(y_{\frac{d}{2}+1}), \cdots , \Im(y_{d})\}$. For sufficiently small $0<\epsilon$ there exists a unique pair of numbers $k$ and $l$ ($1 \leq k \leq l \leq \frac{d}{2}$) such that $A \cap ((\Im(y_0)-\epsilon,\Im(y_0)+\epsilon)=\{r_k,r_l\}$. We call the pair $(k,l)$ the \emph{complex Lefschetz pair} associated with the point $x_0$ and denote it by $\mathcal{L}(x_0)$.
\ede

\brem \label{complex symmetrical}
Note that if $x_0 \in B$ and we have $\mathcal{L}(x_0)=(k,l)$ we always have $l=k+1$. Moreover, if we take $A'=\{\Im(y_1),\cdots ,\Im(y_d)\}\cap \RR^-=\{r_1, \cdots, r_\frac{d}{2}\}$ in a decreasing order, we get an analogous definition to Definition \ref{complex lefschetz pair}, and the two complex Lefschetz pairs achieved from the two analogous definitions are the same.
\erem

\bde \label{orientation}
Let $x_0 \in B$ with complex Lefschetz pair $\mathcal{L}(x_0)=(k,l)$, let $(x_0,y_0),(x_0,y_0')\in S$ be the pairs of points for which $0<\Im(y_0)=\Im(y_0')$ and let $A$ be as in Definition \ref{complex lefschetz pair}. Suppose that $\Re(y_0) < \Re(y_0')$. Then, there exists a small enough $0 < \delta$ such that for every $x \in (x_0-\delta,x_0+\delta)$ if $(x,y)$ is the lift of $x$ to $S$ which is connected to the point $(x_0,y_0)$ and $(x,y')$ is the lift of $x$ to $S$ which is connected to the point $(x_0,y_0')$ we have $\Re(y) < \Re(y')$, and $\delta$ is small enough to determine the complex Lefschetz pair.

Now, take $x_0' \in E \cap \RR$ such that $0<x_0-x_0'<\delta$ and $x_0'' \in E \cap \RR$ such that $0<x_0''-x_0<\delta$. Denote by $\{z_1', \cdots , z_{\frac{d}{2}}'\}=\{z \in K(x_0') \setminus \RR \ | \ \Im(z) \in \RR^+\}$ and by $\{z_1'', \cdots ,z_\frac{d}{2}''\}=\{z \in K(x_0'') \setminus \RR \ | \ \Im(z) \in \RR^+\}$.

We assign to $x_0$ an orientation $\mathcal{O}(x_0) \in \{1,-1,0\}$ as follows:

\ben
\item
If $\Re(z_k')>\Re(z_l')$ and $\Re(z_k'')<\Re(z_l'')$ then, $\mathcal{O}(x_0)=1$.
\item
If $\Re(z_k')<\Re(z_l')$ and $\Re(z_k'')>\Re(z_l'')$ then, $\mathcal{O}(x_0)=-1$.
\item
Otherwise $\mathcal{O}(x_0)=0$.
\een
\ede

We are now going to use the model described in Section \ref{A model for the fiber} in order to describe the algorithm.

\ble \label{model lemma}
Let $\mathcal{R}$ be a connected component in $E \setminus (N \cup B)$, and suppose that we have $n-d$ real points in $K(x)$ for $x \in \mathcal{R}$. Then there exist a continuous family of diffeomorphisms $\{\beta _x \ | \ x \in \mathcal{R}\}$ such that:
\ben
\item
$\beta _x:D \to \widetilde{D}$
\item
$\beta _x(K(x))=K_d$
\item
$\beta _x | _{\RR} = \widetilde{D}\cap \RR$
\item
For all $x,x' \in \mathcal{R}$ we have $\beta _x(y)=\beta _{x'}(y) \in \partial \widetilde{D}$ for all $y \in \partial D$.
\een
\ele

\begin{proof}
The proof follows immediately when noticing that for any two points $x,x'$ there are exactly $n-d$ real points, and that since $\mathcal{R}$ does not contain any point in $B$ the order in terms of $\Im$ coordinate of the non real points in $K(x)$ is fixed. Therefore, the points in $K(x)$ can be identified with the points in $E$, the real points of $K(x)$ with the points in $K'$, and the non real points of $K(x)$ with the points in $K''$.

Now use a diffeomorphism which will transform the real points of $K(x)$ to $K'$ maintaining the order described in Definition \ref{Lefchetz pair}, the non real points of $K(x)$ to the non real points in $\widetilde{D}$ ordered by their $\Im$ coordinate and fixing the image on $\partial D$. Since $S$ is a curve, above $\mathcal{R}$ its elements are continuous so we can make sure that the family of diffeomorphisms is continuous as needed.
\end{proof}

\brem
Although $\beta$ defined in Lemma \ref{model lemma} is not unique, every two families of diffeomorphisms $\beta$ and $\beta '$ are homotopic. This fact will ensure us to have the correct result in the algorithm as we will describe later.
\erem

\bde
Let $s(t)$ be a curve in $E \setminus N$ parameterized by $t \in [0,1]$. $s(t)$ has $n$ lifts into $S$. Projecting these paths to $D$ we get $n$ paths in $D$ starting at $K(s(0))$ and ending at $K(s(1))$. These paths induce a homotopy class of diffeomorphisms called the \emph{Lefschetz diffeomorphism} induced by $s(t)$ and denoted by $\Psi _s$.
\ede

\bpr \label{Model proposition}
Let $x_0 \in N$ and denote $A_0=(x_0,y_0)$ the singular point of $\pi _1$ above $x_0$. Let $s(t)$ be a parameterized semicircle of radius $0< \alpha$ (very small) below the real line centered at $x_0$, where $s(0)=-\alpha$ and $s(1)=\alpha$. Let $\Psi _s$ be the Lefschetz diffeomorphism induced by $s$, and $\mathcal{L}=(k,l)$ be the Lefschetz pair associated with $x_0$. Then,
\begin{footnotesize}
$$\beta _{x_0-\alpha}^{-1}\circ \Psi _s \circ \beta _{x_0+\alpha}=\left \{ 
\begin{array}{ll}
\Delta _{a_1}^{\frac{1}{2}}\left<k,l\right>           & \text{ if } A_0 \text{ is of type } a_1 \\
\Delta _{a_2}^{\frac{1}{2}}\left<k,l\right>           & \text{ if } A_0 \text{ is of type } a_2 \\
\Delta _b^2\left<k,l\right>                           & \text{ if } A_0 \text{ is of type } b \\
\Delta _{c}\left<k,l\right>                           & \text{ if } A_0 \text{ is of type } c \\
\Delta _{d_1}^{\frac{\nu}{2}}\left<k,l\right>         & \text{ if } A_0 \text{ is of type } d_1 \\
\Delta _{d_2}^{\frac{\nu}{2}}\left <k \right >        & \text{ if } A_0 \text{ is of type } d_2 \\
\Delta _{d_3}^{\frac{2}{3}}\left<k\right>             & \text{ if } A_0 \text{ is of type } d_3 \\
\Delta _{e_1}^{\frac{3}{2}}\left<k,l\right>         & \text{ if } A_0 \text{ is of type } e_1 \\
\Delta _{e_2}^{\frac{1}{3}}\left<k,l\right>_{\<a<b\>} & \text{ if } A_0 \text{ is of type } e_2 \text { and } a<b \\
\Delta _{e_2}^{\frac{1}{3}}\left<k,l\right>_{\<a>b\>} & \text{ if } A_0 \text{ is of type } e_2 \text { and } a>b \\
\Delta _{f_1}^{\frac{\nu}{2}}\left<k,l\right>         & \text{ if } A_0 \text{ is of type } f_1 \\
\Delta _{f_2}^{\frac{\nu}{2}}\left<k,l\right>         & \text{ if } A_0 \text{ is of type } f_2 \\
\Delta _{g_1}^{\frac{\nu}{2}}\left<k,l\right>         & \text{ if } A_0 \text{ is of type } g_1 \\
\Delta _{g_2}^{\frac{\nu}{2}}\left<k,l\right>         & \text{ if } A_0 \text{ is of type } g_2 \\
\Delta _{g_3}^{\frac{2}{3}}\left<k,l\right>           & \text{ if } A_0 \text{ is of type } g_3 
\end{array}
\right.
$$
\end{footnotesize}
\epr

\begin{proof}

In order to simplify notations we assume that $D=\widetilde{D}$, the equations of $S$ are all simple, and that the model suitable above the point $\alpha$ is $E_d$ as described in Section \ref{A model for the fiber}.

\noindent
The proof for points of types $a_1,a_2,b,c$ can be found in Proposition $1.2$ of \cite{BGTII}.

\noindent
\underline{$A_0$ is of type $d_1$}: \\
By Remark \ref{k=l+1} we know that $k+1=l$. Moreover, the component of $S$ passing through $A_0$ has equation $(y-(k+\frac{1}{2}))^2-x^\nu=0$, and the other components of $S$ have equations $y=c$ where $c \in K_d \setminus \{k,k+1\}=\{1, \cdots , k-1,k+2,\cdots, n-d,n+1\pm i,\cdots,n+1\pm \frac{d}{2}i\}$. 
Now when $\nu$ is odd, we have that $E_d$ is the model suitable above the point $\alpha$, and that $E_{d+2}$ is the suitalbe model above the point $-\alpha$. In addition, $K(-\alpha)=\{1, \cdots , k-1,k+2,\cdots, n-d,n+1\pm i,\cdots,n+1\pm \frac{d}{2}i\} \cup \{k+\frac{1}{2} \pm i \sqrt{\alpha ^\nu}\}$ and $K({\alpha})=\{1, \cdots , k-1,k+2,\cdots, n-d,n+1\pm i,\cdots,n+1\pm \frac{d}{2}i\} \cup \{k+\frac{1}{2} \pm  \sqrt{\alpha ^\nu}\}$.\\
Therefore, $\beta _{-\alpha}^{-1}:(D,K_{d+2}) \to (D,K(-\alpha))$ is induced by the motion: $n+1 \pm i \mapsto k+\frac{1}{2}\pm i\sqrt{\alpha^\nu}$, and so $\beta _{\alpha}$ is induced by the motion: $k+\frac{1}{2} \pm \sqrt{\alpha^\nu} \mapsto k+\frac{1}{2} \pm \frac{1}{2}$ and $n+1+ij \mapsto n+1+i(j-1)$ where $j=2,\cdots,\frac{d+2}{2}$. \\
Above $s$, $S$ is the union of $n$ paths. Projecting them to $D$ we get the following motions:

$\begin{array}{ll}
y=c & \text{ where } t:0 \mapsto 1 \text{ and } c \in K_d \setminus \{k,k+1\} \\
y=k+\frac{1}{2}+\sqrt{\alpha ^\nu(e^{i\pi +i\pi t})^\nu} & \text{ where } t:0 \mapsto 1 \\
y=k+\frac{1}{2}-\sqrt{\alpha ^\nu(e^{i\pi +i\pi t})^\nu} & \text{ where } t:0 \mapsto 1 
\end{array}$

Therefore, $\Psi _s:(D,K(-\alpha) \to (D,K(\alpha))$ is induced from the motion of the points $k+\frac{1}{2}\pm i \sqrt{\alpha^\nu}$ along a circle of radius $\sqrt{\alpha^\nu}$ centered at $k+\frac{1}{2}$ counterclockwise by angle of $\frac{\nu \pi}{2}$ to the points $k+\frac{1}{2}\pm \sqrt{\alpha^\nu}$.

When $\nu$ is even, the model suitable above the point $-\alpha$ is $E_d$ and so, $K(-\alpha)=K(\alpha)=\{1, \cdots , k-1,k+2,\cdots, n-d,n+1\pm i,\cdots,n+1\pm \frac{d}{2}i\} \cup \{k+\frac{1}{2} \pm \sqrt{\alpha^\nu}\}$. \\
Therefore, $\beta _{-\alpha}^{-1}:(D,K_d) \to (D,K(-\alpha))$ is induced by the motion: $k \mapsto k+\frac{1}{2}-\sqrt{\alpha^\nu}$ and $k+1 \mapsto k+\frac{1}{2}+\sqrt{\alpha^\nu}$. Moreover, $\beta _\alpha=\beta _{-\alpha}$.\\
Above $s$, $S$ is the union of $n$ paths. Projecting them to $D$ we get the following motions:

$\begin{array}{ll}
y=c & \text{ where } t:0 \mapsto 1 \text{ and } c \in K_d \setminus \{k,k+1\} \\
y=k+\frac{1}{2}+\sqrt{\alpha ^\nu(e^{i\pi +i\pi t})^\nu} & \text{ where } t:0 \mapsto 1 \\
y=k+\frac{1}{2}-\sqrt{\alpha ^\nu(e^{i\pi +i\pi t})^\nu} & \text{ where } t:0 \mapsto 1 
\end{array}$

Therefore, $\Psi _s:(D,K(-\alpha) \to (D,K(\alpha))$ is induced from the motion of the points $k+\frac{1}{2}\pm\sqrt{\alpha^\nu}$ along a circle of radius $\sqrt{\alpha^\nu}$ centered at $k+\frac{1}{2}$ counterclockwise by angle of $\frac{\nu \pi}{2}$ to the points $k+\frac{1}{2}\pm \sqrt{\alpha^\nu}$.

Hence, in any case of $\nu$ (odd or even), the composition of these motions is homotopic to the composition of the motions that induce $\Delta _{d_1}^{\frac{\nu}{2}}\left<k,l\right>$. Therefore, the induced diffeomorphisms are homotopic, and we have $\beta _{x_0-\alpha}^{-1}\circ \Psi _s \circ \beta _{x_0+\alpha}=\Delta _{d_1}^{\frac{\nu}{2}}\left<k,l\right>$.

\noindent
\underline{$A_0$ is of type $d_2$}: \\
The proof is similar to the case where $A_0$ is of type $d_1$.

\noindent
\underline{$A_0$ is of type $d_3$}: \\
Since $E_d$ is the model suitable above the point $\alpha$, $E_d$ is also the model suitable above the point $-\alpha$. The component of $S$ passing through $A_0$ has equation $(y-k)^3-x^2=0$, and other components of $S$ have equations $y=c$ where $c\in K_d \setminus \{k,n+1\pm i\}=\{1, \cdots ,k-1,k+1,\cdots,n-d,n+1\pm 2i,\cdots ,n+1 \pm \frac{d}{2}i\}$. 

Now, $K(\alpha)=K(-\alpha)=\{1,\cdots,k-1,k+1,\cdots ,n-d,n+1\pm 2i,\cdots ,n+1\pm i\frac{d}{2}\}\cup\{k+\sqrt[3]{\alpha ^2},k-\frac{1}{2}\sqrt[3]{\alpha ^2}\pm i\frac{\sqrt{3}}{2}\sqrt[3]{\alpha ^2}\}$. \\
$\beta _{-\alpha}^{-1}:(D,K_d) \to (D,K(-\alpha))$ is induced from the motion $n+1\pm i \mapsto k-\frac{\sqrt[3]{\alpha ^2}}{2}\pm i\frac{\sqrt{3}}{2}\sqrt[3]{\alpha ^2}$, and $k \mapsto k+\sqrt[3]{\alpha ^2}$. Moreover, $\beta _\alpha=\beta _{-\alpha}$.

Above $s$, $S$ is the union of $n$ paths. Projecting them to $D$ we get the following motions:

$\begin{array}{ll}
y=c & \text{where } t:0 \mapsto 1 \\
    & c\in K_d \setminus \{k,n+1\pm i\}\\
y=k+(\alpha e^{i\pi+i\pi t})^\frac{2}{3} & \text{ where } t:0 \mapsto 1 \\
y=k-\frac{1}{2}(\alpha e^{i\pi+i\pi t})^\frac{2}{3}+i\frac{\sqrt{3}}{2}(\alpha e^{i\pi+i\pi t})^\frac{2}{3} & \text{ where } t: 0 \mapsto 1 \\
y=k-\frac{1}{2}(\alpha e^{i\pi+i\pi t})^\frac{2}{3}-i\frac{\sqrt{3}}{2}(\alpha e^{i\pi+i\pi t})^\frac{2}{3} & \text{ where } t: 0 \mapsto 1 
\end{array}$

Therefore, $\Psi _s:(D,K(-\alpha) \to (D,K(\alpha))$ is induced from the motion of the points $k+\sqrt[3]{\alpha ^2}$ and $k-\frac{1}{2}\sqrt[3]{\alpha ^2}\pm i\frac{\sqrt{3}}{2}\sqrt[3]{\alpha ^2}$ along a circle of radius $\sqrt[3]{\alpha^2}$ centered at the point $k$ counterclockwise by angle of $\frac{2\pi}{3}$ to themselves.

Hence, the composition of these three motions is homotopic to the composition of the motions that induce $\Delta _{d_3}^{\frac{2}{3}}\left<k\right>$. Therefore, the induced diffeomorphisms are homotopic so we have $\beta _{x_0-\alpha}^{-1}\circ \Psi _s \circ \beta _{x_0+\alpha}=\Delta _{d_3}^{\frac{2}{3}}\left<k\right>$.

\noindent
\underline{$A_0$ is of type $e_1$}: \\
Since $E_d$ is the model suitable above the point $\alpha$ it is also the model suitable above the point $-\alpha$.
By Remark \ref{k=l+1} we know that $k+1=l$. Moreover, the component of $S$ passing through $A_0$ has equation $((y-(k+\frac{1}{2}))^4-x^6=0$, and other components of $S$ have equations $y=c$ where $c\in K_d \setminus \{k,k+1,n+1\pm i\}$.

Now, $K(\alpha)=K(-\alpha)=\{1,\cdots,k-1,k+2,\cdots ,n-d,n+1\pm 2i,\cdots ,n+1\pm i\frac{d}{2}\}\cup\{k+\frac{1}{2}\pm\sqrt{\alpha^3},k+\frac{1}{2}\pm i\sqrt{\alpha^3}\}$.

$\beta _{-\alpha}^{-1}:(D,K_d) \to (D,K(-\alpha))$ is induced from the motion $k \mapsto k+\frac{1}{2}-\sqrt{\alpha^3}$, $k+1 \mapsto k+\frac{1}{2}+\sqrt{\alpha^3}$, $n+1\pm i \mapsto k+\frac{1}{2}\pm i\sqrt{\alpha^3}$ in shortest lines. Moreover, $\beta _\alpha=\beta _{-\alpha}$.

Above $s$, $S$ is the union of $n$ paths. Projecting them to $D$ we get the following motions:

$\begin{array}{ll}
y=c \text{ where } t:0 \mapsto 1 &  c\in K_d \setminus \{k,k+1\,n+1\pm i\} \\
y=k+\frac{1}{2}+\sqrt{(\alpha e^{i\pi+i\pi t})^3} & \text{ where } t: 0 \mapsto 1 \\
y=k+\frac{1}{2}-\sqrt{(\alpha e^{i\pi+i\pi t})^3} & \text{ where } t: 0 \mapsto 1 \\
y=k+\frac{1}{2}+i\sqrt{(\alpha e^{i\pi+i\pi t})^3} & \text{ where } t: 0 \mapsto 1 \\
y=k+\frac{1}{2}-i\sqrt{(\alpha e^{i\pi+i\pi t})^3} & \text{ where } t: 0 \mapsto 1 \\
\end{array}$

Therefore, $\Psi _s:(D,K(-\alpha) \to (D,K(\alpha))$ is induced from the motion of the points $k+\frac{1}{2}\pm\sqrt{\alpha^3}$ and $k+\frac{1}{2}\pm i\sqrt{\alpha^3}$ along a circle of radius $\sqrt{\alpha^3}$ centered at $k+\frac{1}{2}$ counterclockwise by angle of $\frac{3\pi}{2}$ to themselves.

Hence, the composition of these three motions is homotopic to the composition of the motions that induce $\Delta _{e_1}^{\frac{3}{2}}\left<k,l\right>$. Therefore, the induced diffeomorphisms are homotopic so we have $\beta _{x_0-\alpha}^{-1}\circ \Psi _s \circ \beta _{x_0+\alpha}=\Delta _{e_1}^{\frac{3}{2}}\left<k,l\right>$.

\noindent
\underline{$A_0$ is of type $e_2$}: \\
We deal with the two cases for points $A_0$ of type $e_2$ simultaneously. By Remark \ref{k=l+1} we know that $k+1=l$, and we also know that the model suitable above the points $\alpha$ and $-\alpha$ is $E_d$. Moreover, the component of $S$ passing through $A_0$ has equation $((y-(k+\frac{1}{2}))^3+ax^2)((y-(k+\frac{1}{2}))^3-bx^2)=0$, and other components of $S$ have equations $y=c$ where $c\in K_d \setminus \{k,k+1,n+1\pm i,n+1 \pm 2i\}$.

Now, $K(\alpha)=K(-\alpha)=\{1,\cdots,k-1,k+2,\cdots ,n-d,n+1\pm 3i,\cdots ,n+1\pm i\frac{d}{2}\}\cup\{
k+\frac{1}{2}+\sqrt[3]{b\alpha ^2},
k+\frac{1}{2}-\frac{\sqrt[3]{b\alpha ^2}}{2}\pm i\frac {\sqrt{3}}{2}\sqrt[3]{b\alpha ^2},
k+\frac{1}{2}-\sqrt[3]{a\alpha ^2},
k+\frac{1}{2}+\frac{\sqrt[3]{a\alpha ^2}}{2}\pm i\frac {\sqrt{3}}{2}\sqrt[3]{a\alpha ^2}\}$

In order to know $\beta _{-\alpha}^{-1}:(D,K_d) \to (D,K(-\alpha))$ we need to separate between two cases: 

\noindent
\underline{$a>b$} \\
$\beta _{-\alpha}^{-1}:(D,K_d) \to (D,K(-\alpha))$ is induced from the motion 
$k         \mapsto k+\frac{1}{2}-\sqrt[3]{a\alpha^2},
 k+1       \mapsto k+\frac{1}{2}+\sqrt[3]{b\alpha^2},
 n+1\pm i  \mapsto k+\frac{1}{2}-\frac{1}{2}\sqrt[3]{b\alpha^2}\pm i\frac{\sqrt{3}}{2}\sqrt[3]{b\alpha^2},
 n+1\pm 2i \mapsto k+\frac{1}{2}+\frac{1}{2}\sqrt[3]{a\alpha^2}\pm i\frac{\sqrt{3}}{2}\sqrt[3]{a\alpha^2}$ in shortest lines. 
 
\noindent
\underline{$a<b$} \\
$\beta _{-\alpha}^{-1}:(D,K_d) \to (D,K(-\alpha))$ is induced from the motion 
$k         \mapsto k+\frac{1}{2}-\sqrt[3]{a\alpha^2},
 k+1       \mapsto k+\frac{1}{2}+\sqrt[3]{b\alpha^2},
 n+1\pm i  \mapsto k+\frac{1}{2}+\frac{1}{2}\sqrt[3]{a\alpha^2}\pm i\frac{\sqrt{3}}{2}\sqrt[3]{a\alpha^2},
 n+1\pm 2i \mapsto k+\frac{1}{2}-\frac{1}{2}\sqrt[3]{b\alpha^2}\pm i\frac{\sqrt{3}}{2}\sqrt[3]{b\alpha^2}$ in shortest lines. 

Moreover, $\beta _\alpha=\beta_{-\alpha}$.

Above $s$, $S$ is the union of $n$ paths. Projecting them to $D$ we get the following motions:

$\begin{array}{ll}
y=c & \text{where } t:0 \mapsto 1 \\
    &  c\in K_d \setminus \{k,k+1, \\
    &  n+1\pm i,n+1\pm 2i\} \\
k+\frac{1}{2}-\sqrt[3]{a\alpha^2e^{2 \pi i t}} & \text{ where } t:0 \mapsto 1 \\
k+\frac{1}{2}+\sqrt[3]{b\alpha^2e^{2 \pi i t}} & \text{ where } t:0 \mapsto 1 \\
k+\frac{1}{2}+\frac{1}{2}\sqrt[3]{a\alpha^2e^{2 \pi i t}}+i\frac{\sqrt{3}}{2}\sqrt[3]{a\alpha^2e^{2 \pi i t}} & \text{ where } t:0 \mapsto 1 \\
k+\frac{1}{2}+\frac{1}{2}\sqrt[3]{a\alpha^2e^{2 \pi i t}}-i\frac{\sqrt{3}}{2}\sqrt[3]{a\alpha^2e^{2 \pi i t}} & \text{ where } t:0 \mapsto 1 \\
k+\frac{1}{2}-\frac{1}{2}\sqrt[3]{b\alpha^2e^{2 \pi i t}}+i\frac{\sqrt{3}}{2}\sqrt[3]{b\alpha^2e^{2 \pi i t}} & \text{ where } t:0 \mapsto 1 \\
k+\frac{1}{2}-\frac{1}{2}\sqrt[3]{b\alpha^2e^{2 \pi i t}}-i\frac{\sqrt{3}}{2}\sqrt[3]{b\alpha^2e^{2 \pi i t}} & \text{ where } t:0 \mapsto 1 \\
\end{array}$

Therefore, $\Psi _s:(D,K(-\alpha) \to (D,K(\alpha))$ is induced from the motion of the points 
$
 k+\frac{1}{2}-\sqrt[3]{a\alpha^2},
 k+\frac{1}{2}+\sqrt[3]{b\alpha^2},
 k+\frac{1}{2}+\frac{1}{2}\sqrt[3]{a\alpha^2}\pm i\frac{\sqrt{3}}{2}\sqrt[3]{a\alpha^2},
 k+\frac{1}{2}-\frac{1}{2}\sqrt[3]{b\alpha^2}\pm i\frac{\sqrt{3}}{2}\sqrt[3]{b\alpha^2}$
centered at $k+\frac{1}{2}$ counterclockwise by angle of $\frac{2\pi}{3}$ to themselves.

The composition of these three motions is homotopic to the composition of the motions that induce $\Delta _{e_2}^{\frac{1}{3}}\left<k,l\right>_{\<a<b\>}$ or $\Delta _{e_2}^{\frac{1}{3}}\left<k,l\right>_{\<a>b\>}$ (depending whether $a<b$ or $b<a$). Therefore, the induced diffeomorphisms are homotopic and we have $\beta _{x_0-\alpha}^{-1}\circ \Psi _s \circ \beta _{x_0+\alpha}=\Delta _{e_2}^{\frac{1}{3}}\left<k,l\right>$ (the $\<a<b\>$ or $\<b<a\>)$ have been omitted to denote simultaneously both cases). 

\noindent
\underline{$A_0$ is of type $f_1$}: \\
By Remark \ref{k=l+1} we know that $k+2=l$. Moreover, the component of $S$ passing through $A_0$ has equation $(y-(k+1))((y-(k+1))^2-x^\nu)=0$, and the other component of $S$ have equations $y=c$ where $c \in K_d \setminus \{k,k+1,k+2\}=\{1, \cdots , k-1,k+3,\cdots, n-d,n+1\pm i,\cdots,n+1\pm \frac{d}{2}i\}$. 

Now, when $\nu$ is odd, we know that if the model suitable above the point $\alpha$ is $E_d$ then $E_{d+2}$ is the suitable model above the point $-\alpha$. In addition, $K(-\alpha)=\{1, \cdots , k-1,k+3,\cdots, n-d,n+1\pm i,\cdots,n+1\pm \frac{d}{2}i\} \cup \{k+1, k+1 \pm i \sqrt{\alpha ^\nu}\}$ and $K({\alpha})=\{1, \cdots , k-1,k+3,\cdots, n-d,n+1\pm i,\cdots,n+1\pm \frac{d}{2}i\} \cup \{k+1, k+1 \pm  \sqrt{\alpha ^\nu}\}$.\\
Therefore, $\beta _{-\alpha}^{-1}:(D,K_{d+2}) \to (D,K(-\alpha))$ is induced by the motion: $j \mapsto j+2$ where $j=k+1,\cdots,n-d$, $k \mapsto k+1$ and $n+1 \pm i \mapsto k+1\pm i\sqrt{\alpha^\nu}$. Moreover, $\beta _{\alpha}$ is induced by the motion: $k+1 \pm \sqrt{\alpha^\nu} \mapsto k+1 \pm 1$ and $n+1+ij \mapsto n+1+i(j-1)$ where $j=2,\cdots,\frac{d+2}{2}$. \\
Above $s$, $S$ is the union of $n$ paths. Projecting them to $D$ we get the following motions:

$\begin{array}{ll}
y=c & \text{ where } t:0 \mapsto 1 \text{ and } c \in K_d \setminus \{k,k+1,k+2\} \\
y=k+1 & \text{ where } t:0 \mapsto 1 \\
y=k+1+\sqrt{(\alpha e^{i\pi +i\pi t})^\nu} & \text{ where } t:0 \mapsto 1 \\
y=k+1-\sqrt{(\alpha e^{i\pi +i\pi t})^\nu} & \text{ where } t:0 \mapsto 1 
\end{array}$

Therefore, $\Psi _s:(D,K(-\alpha) \to (D,K(\alpha))$ is induced from the motion of the points $k+1\pm i \sqrt{\alpha^\nu}$ along a circle of radius $\sqrt{\alpha^\nu}$ centered at $k+1$ counterclockwise by angle of $\frac{\nu \pi}{2}$ to the points $k+1\pm \sqrt{\alpha^\nu}$.

When $\nu$ is even, we have that the model suitable above the points $\alpha$ and $-\alpha$ is $E_d$. $K(-\alpha)=K(\alpha)=\{1, \cdots , k-1,k+3,\cdots, n-d,n+1\pm i,\cdots,n+1\pm \frac{d}{2}i\} \cup \{ k+1, k+1 \pm \sqrt{\alpha^\nu}\}$. \\
Therefore, $\beta _{-\alpha}^{-1}:(D,K_d) \to (D,K(-\alpha))$ is induced by the motion: $k \mapsto k+1-\sqrt{\alpha^\nu}$ and $k+2 \mapsto k+1+\sqrt{\alpha^\nu}$. Moreover, $\beta _\alpha=\beta _{-\alpha}$.\\
Above $s$, $S$ is the union of $n$ paths. Projecting them to $D$ we get the following motions:

$\begin{array}{ll}
y=c & \text{ where } t:0 \mapsto 1 \text{ and } c \in K_d \setminus \{k,k+1,k+2\} \\
y=k+1 & \text{ where } t:0 \mapsto 1 \\
y=k+1+\sqrt{(\alpha e^{i\pi +i\pi t})^\nu} & \text{ where } t:0 \mapsto 1 \\
y=k+1-\sqrt{(\alpha e^{i\pi +i\pi t})^\nu} & \text{ where } t:0 \mapsto 1 
\end{array}$

Therefore, $\Psi _s:(D,K(-\alpha) \to (D,K(\alpha))$ is induced from the motion of the points $k+1\pm\sqrt{\alpha^\nu}$ along a circle of radius $\sqrt{\alpha^\nu}$ centered at $k+1$ counterclockwise by angle of $\frac{\nu \pi}{2}$ to the points $k+1,k+\frac{1}{2}\pm \sqrt{\alpha^\nu}$.

Hence, in any case of $\nu$ (odd or even), the composition of these motions is homotopic to the composition of the motions that induce $\Delta _{f_1}^{\frac{\nu}{2}}\left<k,l\right>$. Therefore, the induced diffeomorphisms are homotopic. We have $\beta _{x_0-\alpha}^{-1}\circ \Psi _s \circ \beta _{x_0+\alpha}=\Delta _{f_1}^{\frac{\nu}{2}}\left<k,l\right>$.

\noindent
\underline{$A_0$ is of type $f_2$}: \\
The proof is similar to the case where $A_0$ is of type $f_1$.

\noindent
\underline{$A_0$ is of type $g_1$}: \\
By Remark \ref{k=l+1} we know that $k+2=l$. Moreover, the component of $S$ passing through $A_0$ has equation $(y-(k+1)-ax)((y-(k+1))^2-x^\nu)=0$, where $0 < a$ and the other components of $S$ have equations $y=c$ where $c \in K_d \setminus \{k,k+1,k+2\}=\{1, \cdots , k-1,k+3,\cdots, n-d,n+1\pm i,\cdots,n+1\pm \frac{d}{2}i\}$. 

Now when $\nu$ is odd, we know that since the model suitable above the point $\alpha$ is $E_d$, the model suitable above the point $-\alpha$ is $E_{d+2}$. In addition, $K(-\alpha)=\{1, \cdots , k-1,k+3,\cdots, n-d,n+1\pm i,\cdots,n+1\pm \frac{d}{2}i\} \cup \{k+1-a\alpha, k+1 \pm i \sqrt{\alpha ^\nu}\}$ and $K(\alpha)=\{1, \cdots , k-1,k+3,\cdots, n-d,n+1\pm i,\cdots,n+1\pm \frac{d}{2}i\} \cup \{k+1+a\alpha, k+1 \pm  \sqrt{\alpha ^\nu}\}$.\\
Therefore, $\beta _{-\alpha}^{-1}:(D,K_{d+2}) \to (D,K(-\alpha))$ is induced by the motion: $j \mapsto j+2$ where $j=k+1,\cdots,n-d$, $k \mapsto k+1-a\alpha$ and $n+1 \pm i \mapsto k+1\pm i\sqrt{\alpha^\nu}$. Moreover, $\beta _{\alpha}$ is induced by the motion: $k+1 \pm \sqrt{\alpha^\nu} \mapsto k+\frac{1}{2} \pm \frac{1}{2},k+1+a\alpha \mapsto k+2$ and $n+1\pm ij \mapsto n+1\pm i(j-1)$ where $j=2,\cdots,\frac{d}{2}$. \\
Above $s$, $S$ is the union of $n$ paths. Projecting them to $D$ we get the following motions:

$\begin{array}{ll}
y=c & \text{ where } t:0 \mapsto 1 \text{ and } c \in K_d \setminus \{k,k+1,k+2\} \\
y=k+1+a\alpha e^{i\pi+i\pi t} & \text{ where } t:0 \mapsto 1 \\
y=k+1+\sqrt{(\alpha e^{i\pi +i\pi t})^\nu} & \text{ where } t:0 \mapsto 1 \\
y=k+1-\sqrt{(\alpha e^{i\pi +i\pi t})^\nu} & \text{ where } t:0 \mapsto 1 
\end{array}$

Therefore, $\Psi _s:(D,K(-\alpha) \to (D,K(\alpha))$ is induced from the motion of the points $k+1\pm i \sqrt{\alpha^\nu}$ along a circle of radius $\sqrt{\alpha^\nu}$ centered at $k+1$ counterclockwise by angle of $\frac{\nu \pi}{2}$ to the points $k+1\pm \sqrt{\alpha^\nu}$, and the motion of the point $k+1-a\alpha$ to the point $k+1+a\alpha$ along a circle of radius $a\alpha$ counterclockwise by angle of $\pi$.

When $\nu$ is even, we have that $E_d$ is the model suitable for both points $\alpha$ and $-\alpha$. $K(-\alpha)=K(\alpha)=\{1, \cdots , k-1,k+3,\cdots, n-d,n+1\pm i,\cdots,n+1\pm \frac{d}{2}i\} \cup \{ k+1-a\alpha, k+1 \pm \sqrt{\alpha^\nu}\}$. \\
Therefore, $\beta _{-\alpha}^{-1}:(D,K_d) \to (D,K(-\alpha))$ is induced by the motion: $k \mapsto k+1-a\alpha$, $k+1 \mapsto k+1-\sqrt{\alpha^\nu}$ and $k+2 \mapsto k+1+\sqrt{\alpha^\nu}$. Moreover, $\beta _\alpha=\beta _{-\alpha}$.\\
Above $s$, $S$ is the union of $n$ paths. Projecting them to $D$ we get the following motions:

$\begin{array}{ll}
y=c & \text{ where } t:0 \mapsto 1 \text{ and } c \in K_d \setminus \{k,k+1,k+2\} \\
y=k+1+a\alpha e^{i\pi+i\pi t} & \text{ where } t:0 \mapsto 1 \\
y=k+1+\sqrt{(\alpha e^{i\pi +i\pi t})^\nu} & \text{ where } t:0 \mapsto 1 \\
y=k+1-\sqrt{(\alpha e^{i\pi +i\pi t})^\nu} & \text{ where } t:0 \mapsto 1 
\end{array}$

Therefore, $\Psi _s:(D,K(-\alpha) \to (D,K(\alpha))$ is induced from the motion of the points $k+1\pm\sqrt{\alpha^\nu}$ along a circle of radius $\sqrt{\alpha^\nu}$ centered at $k+1$ counterclockwise by angle of $\frac{\nu \pi}{2}$ to themselves, and the motion of the point $k+1-a\alpha$ to the point $k+1+a\alpha$  along a circle of radius $a\alpha$ centered at $k+1$ counterclockwise by angle of $\pi$.

Hence, in any case of $\nu$ (odd or even), the composition of these motions is homotopic to the composition of the motions that induce $\Delta _{g_1}^{\frac{\nu}{2}}\left<k,l\right>$. Therefore, the induced diffeomorphisms are homotopic. We have $\beta _{x_0-\alpha}^{-1}\circ \Psi _s \circ \beta _{x_0+\alpha}=\Delta _{g_1}^{\frac{\nu}{2}}\left<k,l\right>$.

\noindent
\underline{$A_0$ is of type $g_2$}: \\
The proof is similar to the case where $A_0$ is of type $g_1$.

\noindent
\underline{$A_0$ is of type $g_3$}: \\
By Remark \ref{k=l+1} we know that $k+1=l$. We also know that the model suitable above both points $\alpha$ and $-\alpha$ is $E_d$. Moreover, the component of $S$ passing through $A_0$ has equation $(y-k)((y-k)^3-x^2)=0$, and other components of $S$ have equations $y=c$ where $c\in K_d \setminus \{k,k+1,k+2,k+3\}=\{1, \cdots ,k-1,k+4,\cdots,n-d,n+1\pm i,\cdots ,n+1 \pm \frac{d}{2}i\}$. 

Now, $K(\alpha)=K(-\alpha)=\{1,\cdots,k-1,k+4,\cdots ,n-d,n+1\pm i,\cdots ,n+1\pm i\frac{d}{2}\}\cup\{k,k+\sqrt[3]{\alpha ^2},k-\frac{1}{2}\sqrt[3]{\alpha ^2}\pm i\frac{\sqrt{^3}}{2}\sqrt[3]{\alpha ^2}\}$.

$\beta _{-\alpha}^{-1}:(D,K_d) \to (D,K(-\alpha))$ is induced from the motion $n+1\pm i \mapsto k-\frac{1}{2}\sqrt[3]{\alpha ^2}\pm i\frac{\sqrt{^3}}{2}\sqrt[3]{\alpha ^2}$, $j \mapsto j+2$ where $j=k+2,\cdots, n-d$ and $k+1 \mapsto k+\sqrt[3]{\alpha ^2}$ in shortest lines. Moreover, $\beta _\alpha=\beta _{-\alpha}$.

Above $s$, $S$ is the union of $n$ paths. Projecting them to $D$ we get the following motions:

$\begin{array}{ll}
y=c \text{ where } t:0 \mapsto 1 &  c\in K_d \setminus \{k,k+1,k+2,k+3\} \\
y=k+\sqrt[3]{(\alpha e^{i \pi  +i \pi t})^2} & \text{ where } t:0 \mapsto 1 \\
y=(k-\frac{1}{2}\sqrt[3]{(\alpha e^{i \pi  +i \pi t})^2}+i\frac{\sqrt{3}}{2}\sqrt[3]{(\alpha e^{i \pi  +i \pi t})^2} & \text{ where } t: 0 \mapsto 1 \\
y=(k-\frac{1}{2}\sqrt[3]{(\alpha e^{i \pi  +i \pi t})^2}-i\frac{\sqrt{3}}{2}\sqrt[3]{(\alpha e^{i \pi  +i \pi t})^2} & \text{ where } t: 0 \mapsto 1 \\
y=k & \text{ where } t:0 \mapsto 1 
\end{array}$

Therefore, $\Psi _s:(D,K(-\alpha) \to (D,K(\alpha))$ is induced from the motion of the points $k,k+\sqrt[3]{\alpha ^2}$ and $(k-\frac{1}{2}\sqrt[3]{\alpha ^2}\pm i\frac{\sqrt{^3}}{2}\sqrt[3]{\alpha ^2})$ along a circle of radius $\sqrt[3]{\alpha^2}$ centered at $k$ counterclockwise by angle of $\frac{2\pi}{3}$ to themselves.

The composition of these three motions is homotopic to the composition of the motions that induce $\Delta _{g_3}^{\frac{1}{3}}\left<k,l\right>$. Therefore, the induced diffeomorphisms are homotopic so we have $\beta _{x_0-\alpha}^{-1}\circ \Psi _s \circ \beta _{x_0+\alpha}=\Delta _{g_3}^{\frac{1}{3}}\left<k,l\right>$.

\end{proof}

\bpr \label{Lefschetz diffeomorphism for B}
Let $x_0 \in B$. Let $u$ be a parameterized (from $x_0-\alpha$ to $x_0+\alpha$) the straight line segment $[x_0-\alpha,x_0+\alpha]$ ($\alpha$ very small), and let $\Psi _u$ be the Lefschetz diffeomorphism induced by $u$, and $\mathcal{L}=(k,l)$ be the complex Lefschetz pair associated with $x_0$. Then,
$$\Delta _{com}^{\mathcal{O}(x_0)}\<k,l\>=\beta _{x_0-\alpha}^{-1} \Psi _u \beta _{x_0+\alpha}$$
\epr

\begin{proof}
First it follows from Remark \ref{complex symmetrical} that the diffeomorphism must be symmetrical along the real axis; Therefore, we concentrate on the half of the complex plane which consists of the points with positive imaginary coordinate, keeping in mind that what ever happens there also happens in the other half-plane in a mirror like image. Recall from Definition \ref{complex lefschetz pair} that if $(x_0,y_0),(x_0,y_0')\in S$ is one of the pairs of points for which $0<\Im(y_0)=\Im(y_0')$ and $A$ as in Definition \ref{complex lefschetz pair} we take small enough $0<\alpha$ such that the complex Lefschetz pair is well defined and that for every $x \in (x_0-\alpha,x_0+\alpha)$ if $(x,y)$ is the lift of $x$ to $S$ which is connected to the point $(x_0,y_0)$ and $(x,y')$ is the lift of $x$ to $S$ which is connected to the point $(x_0,y_0')$ we have $\Re(y) < \Re(y')$ when $\Re(y_0) < \Re(y_0')$ or $\Re(y) > \Re(y')$ when $\Re(y_0) > \Re(y_0')$ (It can't be that $\Re(y_0) = \Re(y_0')$ or else the point $x_0$ will be in $N$ instead of $B$). Without lose of generality we assume that $\Re(y_0) < \Re(y_0')$.

Now as in Definition \ref{orientation}, let $x_0' \in E \cap \RR$ such that $0<x_0-x_0'<\alpha$ and $x_0'' \in E \cap \RR$ such that $0<x_0''-x_0<\alpha$. Denote by $\{z_1', \cdots , z_{\frac{d}{2}}'\}=\{z \in K(x') \setminus \RR \ | \ \Im(z) \in \RR^+\}$ and by $\{z_1'', \cdots ,z_\frac{d}{2}''\}=\{z \in K(x'') \setminus \RR \ | \ \Im(z) \in \RR^+\}$. 

$\beta _{x_0-\alpha}^{-1}:(D,K_d) \to (D,K(\alpha))$ is induced by the motion of the point $n+1+ki \mapsto z_k'$, $n+1+li \mapsto z_l'$ and the points $n+1-ki,n+1-li$ goes to the complex conjugates of $z_k'$ and $z_l'$ along straight lines (we assume that all other components of $S$ have equations $y=c$ where $c \in K_d \setminus \{n+1\pm ik,n+1 \pm il\}$. 

$\beta _{x_0+\alpha}:(D,K(\alpha)) \to (D,K_d)$ is induced by the motion of the points $z_k'',z_l''$ and their complex conjugates to the points $n+1\pm ki,n+1\pm li$ along straight lines.

We look at the following cases:

\noindent
\underline{$\Re(z_k')>\Re(z_l')$ and $\Re(z_k'')<\Re(z_l'')$} \\
In this case the point $(x_0',z_k')$ is connected to the point $(x_0,y_0')$ and the point $(x_0',z_l')$ is connected to the point $(x_0,y_0)$. Because $k<l$ we know by Lemma \ref{model lemma} that $\Im(z_k')<\Im(z_l')$. Now, since $\Re(z_k'')<\Re(z_l'')$ the point $(x_0'',z_k'')$ is connected to the point $(x_0,y_0)$ and the point $(x_0'',z_l'')$ is connected to the point $(x_0,y_0')$. However, $\Im(z_k'') < \Im(z_l'')$ hence we have a motion in which the point on the left goes down and the point on the right goes up crossing each other $\Im$ coordinate once (see Figure \ref{O=1/2} (left)). Therefore, the composition $\beta _{x_0-\alpha}^{-1} \Psi _u \beta _{x_0+\alpha}$ is homotopic to the diffeomorphism defining $\Delta _{com}^1\<k,l\>$ (see Figure \ref{O=1/2} (right)).

\myfigure{
\begin{center}
\includegraphics[scale=0.70]{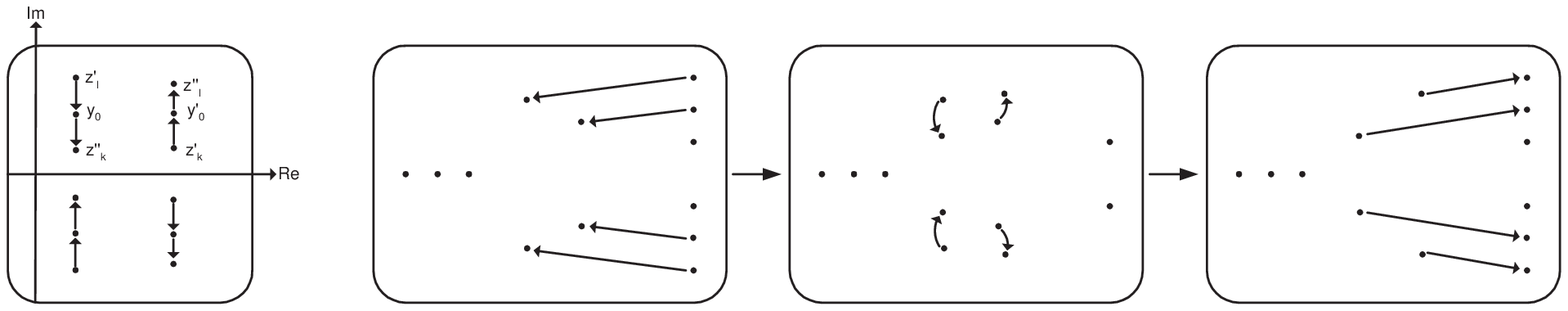}
\caption{Positive motion in the fiber (left), $\beta _{x_0-\alpha}^{-1} \Psi _u \beta _{x_0+\alpha}$ (right)}\label{O=1/2}
\end{center}
}

\noindent
\underline{$\Re(z_k')<\Re(z_l')$ and $\Re(z_k'')>\Re(z_l'')$}
In this case the point $(x_0',z_k')$ is connected to the point $(x_0,y_0)$ and the point $(x_0',z_l')$ is connected to the point $(x_0,y_0')$. Because $k<l$ we know by Lemma \ref{model lemma} that $\Im(z_k')<\Im(z_l')$. Now, since $\Re(z_k'')>\Re(z_l'')$ the point $(x_0'',z_k'')$ is connected to the point $(x_0,y_0')$ and the point $(x_0'',z_l'')$ is connected to the point $(x_0,y_0)$. However, $\Im(z_k'') < \Im(z_l'')$ hence we have a motion in which the point on the left goes up and the point on the right goes down crossing each other $\Im$ coordinate once (see Figure \ref{O=-1/2} (left)). Therefore, the composition $\beta _{x_0-\alpha}^{-1} \Psi _u \beta _{x_0+\alpha}$ is homotopic to the diffeomorphism defining $\Delta _{com}^{-1}\<k,l\>$ (see Figure \ref{O=-1/2} (right)).

\myfigure{
\begin{center}
\includegraphics[scale=0.70]{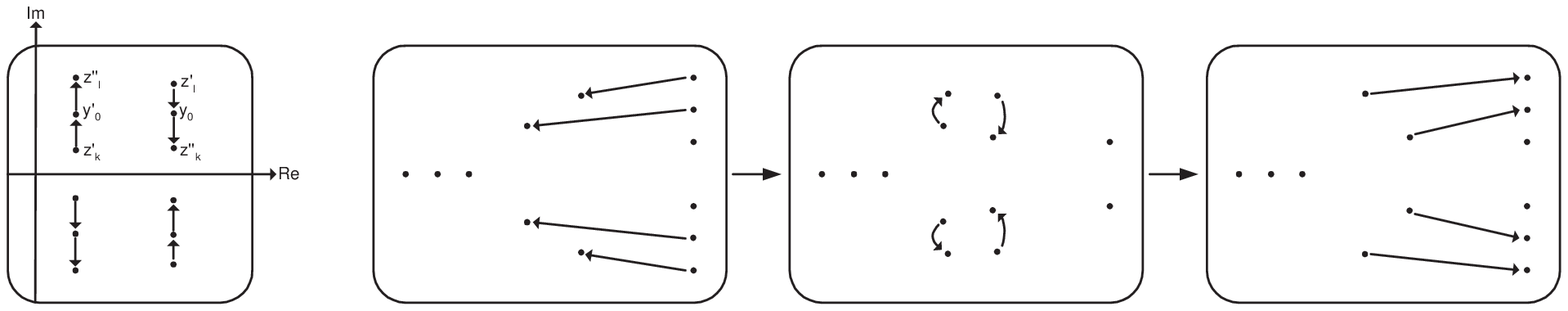}
\caption{Negative motion in the fiber (left), $\beta _{x_0-\alpha}^{-1} \Psi _u \beta _{x_0+\alpha}$ (right)}\label{O=-1/2}
\end{center}
}

\noindent
\underline{Otherwise}
In these cases we have either $\Re(z_k')<\Re(z_l')$ and $\Re(z_k'')<\Re(z_l'')$ or $\Re(z_k')>\Re(z_l')$ and $\Re(z_k'')>\Re(z_l'')$. Hence the point on the left and on the right move until their $\Im$ coordinates coincide and then retract back. In this case there is no crossing therefore the diffeomorphism induced by the motion which induce $\beta _{x_0-\alpha}^{-1} \Psi _u \beta _{x_0+\alpha}$ is trivial as can be seen in Figure \ref{O=0} (one of two possible cases).

\myfigure{
\begin{center}
\includegraphics[scale=0.70]{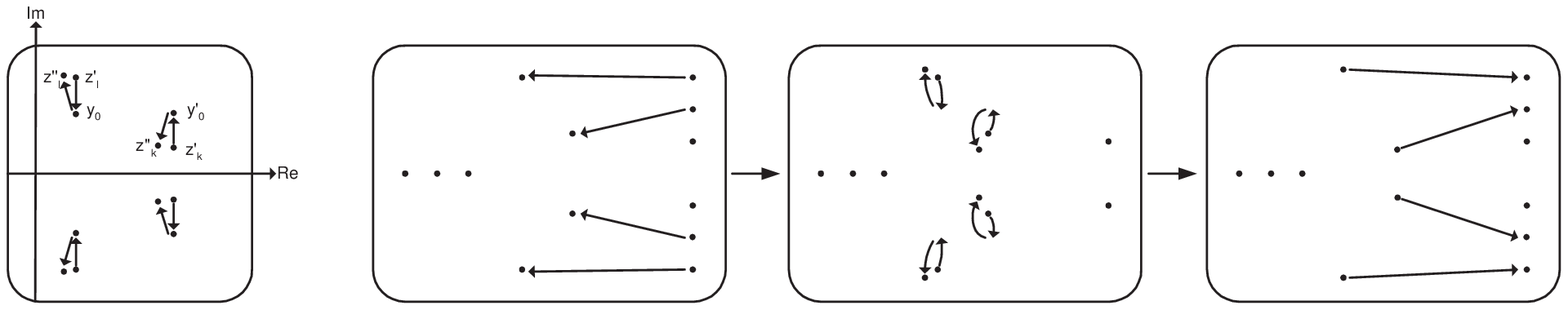}
\caption{Trivial motion in the fiber (left), $\beta _{x_0-\alpha}^{-1} \Psi _u \beta _{x_0+\alpha}$ (right)}\label{O=0}
\end{center}
}

\end{proof}

We need to be able to describe the braids of the local braid monodromy at each point of the set $N$ in a way that will make it easy to conjugate them. This can be achieved by using a product of half-twists that define the resulted braid. This product is the parallel of the skeleton definition given in \cite{BGTI} and \cite{BGTII}.

\bde
Let $\widetilde{K} \subset \widehat{K} \subset \widetilde{D}$. A \emph{skeleton} in $(\widetilde{D},\widehat{K},\widetilde{K})$ is a consecutive sequence of paths each starting and ending in one of the points of $\widetilde{K}$. Each path starts where it's previous ends. A skeleton is denoted by $\xi$.
\ede

\bde
We call a product of half-twists that describes the local braid monodromy at the point $x_0 \in N$ the \emph{generalized skeleton} of $x_0$, and denote it by $\xi_{x_0}$. 
\ede

Recall that half-twists are defined by paths, so in order to describe the proper generalized skeleton for each type of singularity, it is enough to describe the paths defining its half-twists, and give their product form. 

\bde \label{skeletons}
Let $x_0 \in N$, and $A_0=(x_0,y_0)$ be the singular point of $\pi _1$ above $x_0$. Let $\mathcal{L}_{x_0}=(k,l)$ be the Lefschetz pair associated with the point $x_0$. We define the generalized skeletons:
\ben
\item 
If $A_0$ is of singularity type $a_1$, the generalized skeleton assigned to it is $h_1$, where $h_1$ is described in Figure \ref{Generalized skeleton} (a).
\item
If $A_0$ is of singularity type $a_2$, the generalized skeleton assigned to it is $h_1$, where $h_1$ is described in Figure \ref{Generalized skeleton} (b).
\item
If $A_0$ is of singularity type $b$, the generalized skeleton assigned to it is $h_1^2$, where $h_1$ is described in Figure \ref{Generalized skeleton} (a).
\item
If $A_0$ is of singularity type $c$, the generalized skeleton assigned to it is $(h_1\cdots h_{l-k})(h_1\cdots h_{l-k-1})\cdots(h_1 h_2)(h_1)$, where $h_i$ are described in Figure \ref{Generalized skeleton} (c).
\item
If $A_0$ is of singularity type $d_1$, the generalized skeleton assigned to it is $h_1^\nu$, where $h_1$ is described in Figure \ref{Generalized skeleton} (a).
\item
If $A_0$ is of singularity type $d_2$, the generalized skeleton assigned to it is $h_1^\nu$, where $h_1$ is described in Figure \ref{Generalized skeleton} (b).
\item
If $A_0$ is of singularity type $d_3$, the generalized skeleton assigned to it is $(h_1h_2)^2$, where $h_i$ are described in Figure \ref{Generalized skeleton} (d).
\item
If $A_0$ is of singularity type $e_1$, the generalized skeleton assigned to it is $(h_2h_1h_3)^6$, where $h_1,h_2,h_3$ are described in Figure \ref{Generalized skeleton} (e).
\item
If $A_0$ is of singularity type $e_2$, the generalized skeleton assigned to it is $h_2h_4h_1h_3h_5$, where $h_i$ are described in Figure \ref{Generalized skeleton} (f) when $b<a$, and in Figure \ref{Generalized skeleton} (g) when $a<b$.
\item
If $A_0$ is of singularity type $f_1$, the generalized skeleton assigned to it is $(h_1h_2h_1)^\nu$, where $h_1,h_2$ are described in Figure \ref{Generalized skeleton} (c).
\item
If $A_0$ is of singularity type $f_2$, the generalized skeleton assigned to it is $(h_2h_1h_2)^\nu$, where $h_1,h_2$ are described in Figure \ref{Generalized skeleton} (h).
\item
If $A_0$ is of singularity type $g_1$, the generalized skeleton assigned to it is $h_2h_1h_2^\nu h_1h_2$, where $h_1,h_2$ are described in Figure \ref{Generalized skeleton} (c).
\item
If $A_0$ is of singularity type $g_2$, the generalized skeleton assigned to it is $h_2h_1h_2^\nu h_1h_2$, where $h_1$ is described in Figure \ref{Generalized skeleton} (h).
\item
If $A_0$ is of singularity type $g_3$, the generalized skeleton assigned to it is $(h_1h_2h_3h_2)^\nu$, where $h_1,h_2,h_3$ are described in Figure \ref{Generalized skeleton} (i).
\een
\ede

\myfigure{
\begin{center} 
\includegraphics[scale=0.65]{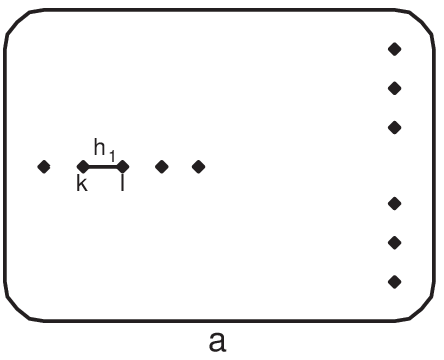}
\includegraphics[scale=0.65]{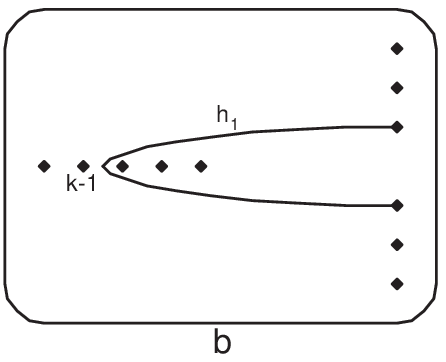}
\includegraphics[scale=0.65]{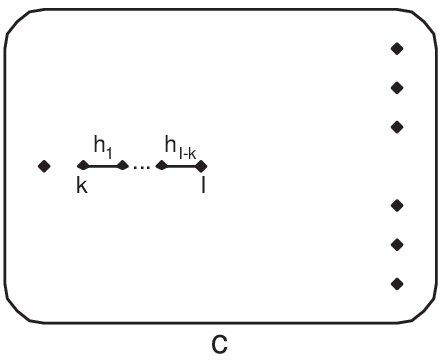}\\
\includegraphics[scale=0.65]{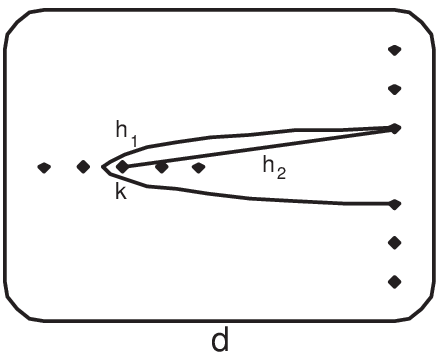}
\includegraphics[scale=0.65]{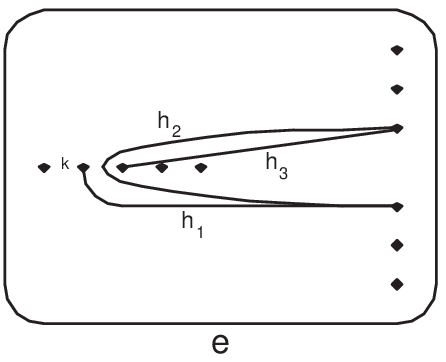}
\includegraphics[scale=0.65]{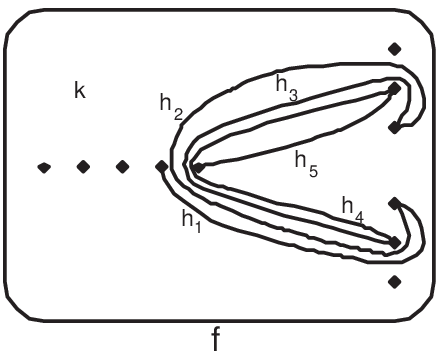}\\
\includegraphics[scale=0.65]{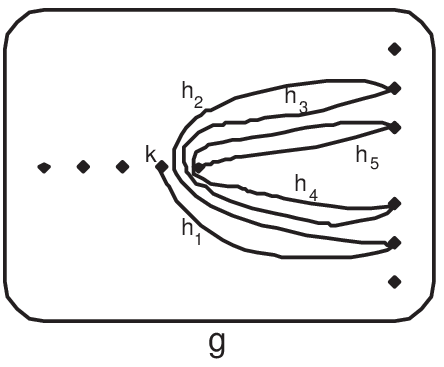}
\includegraphics[scale=0.65]{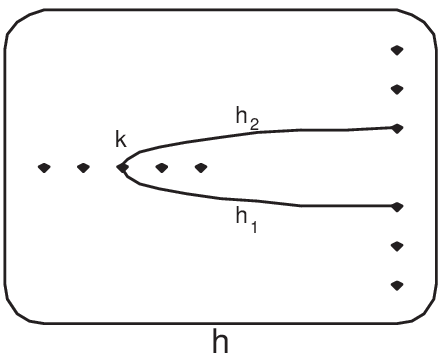}
\includegraphics[scale=0.65]{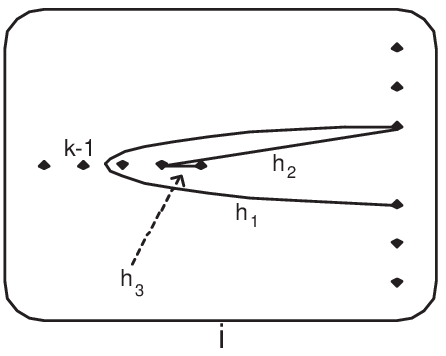}
\caption{Half-twists used in products defining generalized skeletons}
\label{Generalized skeleton}
\end{center}
}

\bth
Let $x_0 \in N$, and $A_0=(x_0,y_0)$ be the singular point of $\pi _1$ above $x_0$. Let $\mathcal{L}_{x_0}=(k,l)$ be the Lefschetz pair associated with the point $x_0$. Then, the generalized skeletons defined in Definition \ref{skeletons}, describe the local braid monodromy at the point $x_0$.
\eth

\begin{proof}
The proof follows immediately from Proposition \ref{Prop d_2} to Proposition \ref{Prop g_3} and the proofs in \cite{BGTII}.
\end{proof}

%%%%%%%%%%%%%%%%%%%%%%%%%%%%%%%%%%%%%%%%%%%%%%%%%%%%%%%%%%%%%%%%%%%%%%%%%%%%%%%%%%%%%%%%%%%%%
\subsection{Description of a generating set for $\pi _1(E \setminus N)$}

$ $\\

We need to compute the braid monodromy homomorphism. Recall that $\varphi$ is a homomorphism from $\pi _1(E \setminus N)$, which is the fundamental group of a punctured disk, to the braid group $B_n$. Therefore, in order to describe the mapping, we choose a generating set for $\pi _1(E \setminus N)$, and compute the image of $\varphi$ on these elements. The image of $\varphi$ on this generating set not only describe completely the braid monodromy homomorphism but also serve as an important step in the computation of the fundamental group of the complement of the curve.

We fix as a base point for the fundamental group of the origin of $\varphi$ the point $x_0$ which is the larger of the two points from $\partial E \cap \RR$. Moreover we enumerate the points of $N$ such that $\{x_1 > \cdots > x_p\}$.  Then, we take an ordered system of paths below the real axis. Each of the paths circles exactly one point in $N$, as can be seen in Figure \ref{generating set}. We denote the elements of this generating set by $\Gamma _1, \cdots \Gamma _p$.

\myfigure{
\begin{center}
\includegraphics[scale=0.8]{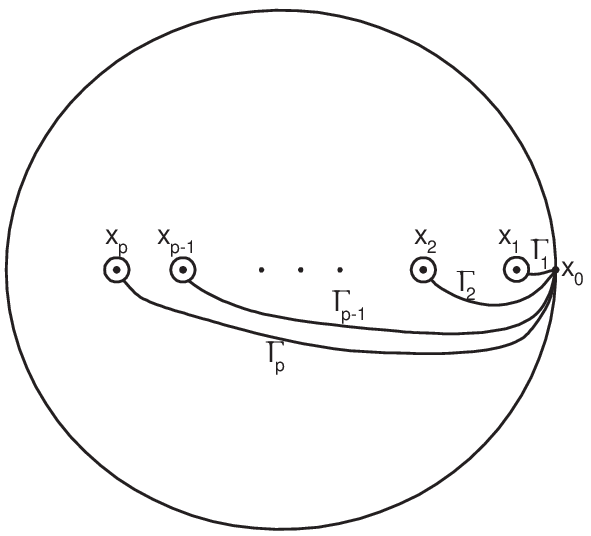}
\caption{The generating set for $\pi _1(E \setminus N)$} \label{generating set}
\end{center}
}

Now we turn to describe these paths.

To each point $x_j \in N$ we attach a real value $0 < \alpha_j$ which is not only sufficiently small in order to determine the Lefschetz pair $\mathcal{L}(x_j)=(k_j,l_j)$ for the point $x_j$ but also small enough such that $(x_j-\alpha _j,x_j+\alpha _j)\cap (N \cup B)=\{x_j\}$. Denote by $t_0$ the real segment between $x_1+\alpha _1$ and $x_0$, and for any $1\leq j < p$, denote $t_
j(t)$ $t \in [0,1]$ the real segment $[x_{j+1}+\alpha_{j+1},x_j-\alpha _j]$.

Denote also by $s_j(t)$ the parametrization of a semi circle below the real axis centered at $x_j$ with radius $\alpha _j$, and $r_j(t)$ the full circle oriented counterclockwise centered at $x_j$ with radius $\alpha _j$ for $t \in [0,1]$, such that $s_j(0)=x_j-\alpha _j$ and $s_j(1)=x_j+\alpha _j$, and $r_j(0)=r_j(1)=x_j+\alpha _j$.
We also denote the type of singularity of the singular point of $\pi _1$ above $x_j$ by $p_j$ (i.e., $p_j \in \{a_1,a_2,b,c,d_1,d_2,d_3,e_1,e_2,f_1,f_2,g_1,g_2,g_3\}$.

The closed disk $E$ with the notations above can be seen in Figure \ref{pi1}.

\myfigure{
\begin{center}
\includegraphics[scale=0.5]{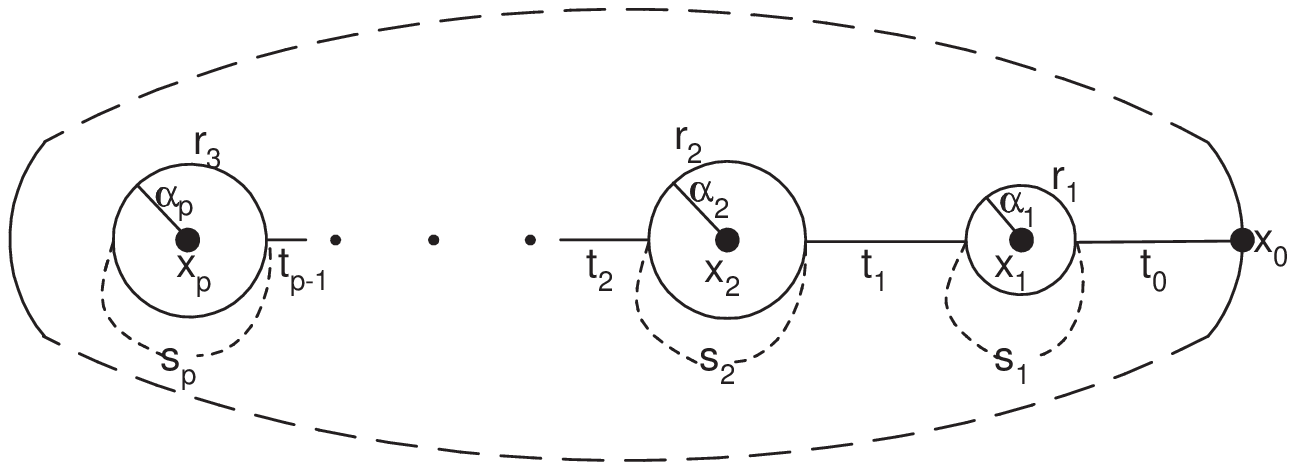}
\caption{$\pi _1(E \setminus N)$ with notations} \label{pi1}
\end{center}
}
$\Gamma _j$ is now homotopic to the following path:
$$\Gamma _j=t_0^{-1}\left (\prod _{k=1}^{j-1}(s_k^{-1}t_k^{-1})\right )r_j\left (\prod _{k=j-1}^{1}(t_ks_k)\right )t_0.$$
We take a closer look at the segment $t_k$. Denote the ordered set of points $t_k \cap B=\{b_{k,1} > \cdots > b_{k,n_k}\}$ the points of the set $B$ in the segment $t_k$. To each of these points, using definition \ref{complex lefschetz pair}, we attached a positive small real number $0 < \delta _{k,l}$ where $l=1,\cdots ,n_k$. The set of points $\{x_k-\alpha _k,b_{k,1}+\delta _{k,1},b_{k,1}-\delta _{k,1}, \cdots , b_{k,n_k}+\delta _{k,n_k},b_{k,n_k}-\delta _{k,n_k},x_{k+1}+\alpha _{k+1}\}$ induce a partition on $t_k$. Here too, we take $\delta _{k,l}$ (and $\alpha$) small enough such that $(b_{k,l}-\delta_{k,l},b_{k,l}+\delta _{k,l}) \cap (N \cup B)=\{b_{k,l}\}$, and $(b_{k,l}-\delta_{k,l},b_{k,l}+\delta_{k,l}) \cap (x_j-\alpha_j,x_j+\alpha_j)=\emptyset$ For all $k,l,j$.

Denote:
\ben
\item
$t_{k,0}=[b_{k,1}+\delta _{k,1},x_k-\alpha _k] \ \forall k \neq 0$.
\item
$t_{k,l}=[b_{k,l+1}+\delta _{k,l+1},b_{k,l}-\delta _{k,l}] \ \forall l=1,\cdots ,n_k-1$.
\item
$t_{k,n_k}=[x_{k+1}+\alpha _{k+1},b_{k,n_k}-\delta _{k,n_k}]$.
\item
$u_{k,l}=[b_{k,l}-\delta _{k,l},b_{k,l}+\delta _{k,l}] \ \forall l=1,\cdots ,n_k$.
\item
$t_{0,0}=[b_{0,1}+\delta_{0,1},x_0]$.
\een

This partition implies that:
$$t_k=\left (\prod _{l=n_k}^{1}t_{k,l}u_{k,l}\right )t_{k,0}.$$
Hence,
\begin{small}
$$\Gamma _j=$$
$$\left (t_{0,0}^{-1}(\prod _{l=1}^{n_0}u_{0,l}^{-1}t_{0,l}^{-1})\right)
\prod _{k=1}^{j-1}(s_k^{-1} t_{k,0}^{-1}(\prod _{l=1}^{n_k}u_{k,l}^{-1}t_{k,l}^{-1}))
r_j
\prod _{k=j-1}^{1}((\prod _{l=n_k}^1t_{k,l}u_{k,l})t_{k,0} s_k)
\left((\prod _{l=n_0}^1t_{0,l}u_{0,l})t_{0,0}\right).$$
\end{small}

%%%%%%%%%%%%%%%%%%%%%%%%%%%%%%%%%%%%%%%%%%%%%%%%%%%%%%%%%%%%%%%%%%%%%%%%%%%%%%%%%%%%%%%%%%%%%%%%%%%%%%%%%%%%%%%%%%%%%%%%%%
\subsection{The computation of $\varphi(\Gamma _j)$}

Using the definitions above it follows that if
\mytable{
\begin{tabular}{rcl}
$\mathcal{Q}$ & = & $\left (\prod _{k=j-1}^1\Psi _{t_k} \Psi _{s_k}\right )\Psi_{t_0}.$ \\
              & = & $\prod _{k=j-1}^1 \left ((\prod _{l=n_k}^1 \Psi_{t_{k,l}}\Psi_{u_{k,l}})\Psi_{t_{k,0}}\Psi_{s_k}\right )
                    \left (\prod_{l=n_0}^1 \Psi_{t_{0,l}} \Psi _{u_{0,l}}\right ) \Psi_{t_{0,0}}$. \\
\end{tabular}
}
Then, the Lefschetz diffeomorphism of $\Gamma _j$ is the product:
$$\Psi_{\Gamma _j}=\mathcal{Q}^{-1}\Psi _{r_j}\mathcal{Q}$$
which implies that
$$\varphi(\Gamma _j)=\mathcal{Q}^{-1}\Delta _{p_j}\mathcal{Q}.$$

Where $\Delta _{p_j}$ is the local braid monodromy defined in Section \ref{A model for the fiber}, which complies with the topological type of the singular point $x_j$ denoted $p_j$, and with the Lefschetz pair associated with it $\mathcal{L}(x_j)$.

In the model defined above this is realized as follows:
\mytable{
\begin{tabular}{rcl}
$\Psi _{\Gamma _j}$ & = & $\beta _{x_0}^{-1}\mathcal{Q}^{-1}\Delta _{p_j}\mathcal{Q}\beta _{x_0}$ \\
                    & = & $(\beta _{x_0}^{-1}\mathcal{Q}^{-1}\beta _{x_j+\alpha _j})(\beta _{x_j+\alpha _j}^{-1}\Delta _{p_j}\beta _{x_j+\alpha _j})(\beta _{x_j+\alpha _j}^{-1}\mathcal{Q}\beta _{x_0}).$ \\
\end{tabular}
}

We turn now to $\beta _{x_j+\alpha _j}^{-1}\mathcal{Q}\beta _{x_0}$.

\mytable{
\begin{tabular}{rcl}
$\beta _{x_j+\alpha _j}^{-1}\mathcal{Q}\beta _{x_0}$ & = & $\beta _{x_j+\alpha_j}^{-1}\left (\prod _{k=j-1}^1\Psi _{t_k}\Psi_{s_k}\right ) \Psi_{t_0} \beta _{x_0}$ \\
                    & = & $\left ( \prod_{k=j-1}^1(\beta _{x_{k+1}+\alpha_{k+1}}^{-1} \Psi_{t_k} \beta _{x_k-\alpha _k})
                           (\beta _{x_k-\alpha _k}^{-1} \Psi _{s_k} \beta _{x_k+\alpha _k}) \right )
                           (\beta _{x_1+\alpha _1}^{-1} \Psi_{t_0} \beta _{x_0})$
\end{tabular}
}

By Proposition \ref{Model proposition} we know that $\beta _{{x_k}-\alpha _k}^{-1} \Psi _{s_k} \beta _{{x_k}+\alpha _k}=\Delta _{p_k}\<k,l\>$ where $p_k$ is the singularity type as defined above, and the power is omitted. Therefore, we need to study $\beta _{x_{k+1}+\alpha _{k+1}}^{-1}\Psi _{t_k}\beta _{x_k+\alpha _k}$. However, $t_k$ is actually a product so if one identifies the point $b_{k,n_k+1}+\delta_{k,n_k+1}$ as the point $x_{k+1}+\alpha _{k+1}$ and takes $\alpha _0=0$ we have:

\mytable{
\begin{tabular}{rcl}
$\beta _{x_{k+1}+\alpha _{k+1}}^{-1} \Psi _{t_k} \beta _{x_k-\alpha _k}$ & = & 
            $\beta _{x_{k+1}+\alpha _{k+1}}^{-1}\left ( \prod _{l=n_k}^1 \Psi_{t_{k,l}} \Psi_{u_{k,l}} \right ) \Psi _{t_{k_0}} \beta _{x_k-\alpha _k}$\\ 
 & = & $\left ( \prod _{l=n_k}^1(\beta _{b_{k,l+1}+\delta _{k,l+1}}^{-1} \Psi _{t_{k,l}} \beta _{b_{k,l}-\delta_{k,l}})
                               (\beta _{b_{k,l}-\delta_{k,l}}^{-1} \Psi _{u_{k,l}} \beta _{b_{k,l}+\delta _{k,l}}) \right )$ \\
 &   &                        $(\beta _{b_{k,1}+\delta _{kl,1}} \Psi _{t_{k,0}} \beta _{x_k-\alpha_k})$.
\end{tabular}
}

However, by Lemma \ref{model lemma} we have that the following four diffeomorphisms all induce the trivial braid:
\ben
\item
$\beta _{b_{0,1}+\delta_{0,1}}^{-1} \Psi_{t_{0,0}} \beta _{x_0}$
\item
$\beta _{b_{k,l+1}+\delta_{k,l+1}}^{-1} \Psi_{t_{k,l}} \beta _{b_{k,l}-\delta_{k,l}}$ for $k=0,\cdots ,p-1$ and $l=1,\cdots ,n_j-1$.
\item
$\beta _{b_{k,1}+\delta_{k,1}}^{-1} \Psi_{t_{k,0}} \beta _{x_k-\alpha _k}$ for $k=1,\cdots ,p-2$.
\item
$\beta _{x_{k+1}+\alpha _{k+1}}^{-1} \Psi_{t_{k,n_k}} \beta _{b_{k,n_k}-\delta_{k,n_k}}$ for $k=1,\cdots ,p-1$.
\een
Therefore, by Propositions \ref{Model proposition} and \ref{Lefschetz diffeomorphism for B} as a braid we get the following diffeomorphism of the model:

$$\beta _{x_j+\alpha _j}^{-1} \mathcal{Q} \beta _{x_0}$$
$$=\prod _{k=j-1}^1 \left ( \left (\prod_{l=n_k}^1 (\beta _{b_{k,l}-\delta_{k,l}}^{-1} \Psi _{u_{k,l}}\beta _{b_{k,l}+\delta_{k,l}})\right )(\beta _{x_k-\alpha _k}^{-1} \Psi _{s_k} \beta _{x_k+\alpha _k})\right)(\prod_{l=n_0}^1(\beta _{b_{0,l}-\delta_{0,l}}^{-1} \Psi _{u_{0,l}} \beta _{b_{0,l}+\delta _{0,l}}) )$$
$$=\prod _{k=j-1}^1 \left ( \left (\prod_{l=n_k}^1 \Delta _{com}^{\mathcal{O}(b_{k,l})} \<k_{k,l},l_{k,l}\> \right ) \Delta _{p_k}^{e_k}\<k_{x_k},l_{x_k}\> \right)(\prod_{l=n_0}^1\Delta _{com}^{\mathcal{O}(b_{0,l})}\<k_{0,l},l_{0,l}\> )$$
where $\Delta _{p_k}^{e_k}\<k_{x_k},l_{x_k}\>$ is the diffeomorphism in Proposition \ref{Model proposition}.

Recall that by Definition \ref{skeletons} each $\Delta _{p_j}$ is a generalized skeleton which is a product of half-twists. Therefore, it is possible to write $\Delta _{x_j}=\prod _m (H_m)$ where $H_m$ are the half-twists from Definition \ref{skeletons}.

Moreover, it was proved in \cite{BGTI} that on a half-twist conjugation is done by activating on the path defining the half-twist the diffeomorphism conjugating it. Therefore, we can calculate the braid monodromy $\varphi(\Gamma _j)$ using the following:

\bth \label{The algorithm theorem}
Let $S,N,x_j,p_j,\mathcal{L}(x_j)=\<k_{x_j},l_{x_j}\>,\mathcal{L}(b_{k,l})=\<k_{k,l},l_{k,l}\>,e_k,\alpha _j$ and $\delta _{k,l}$ be as above. Then,
\begin{footnotesize}
$$\varphi(\Gamma _j)=\prod_m\left(\prod _{k=j-1}^1 \left ( \left (\prod_{l=n_k}^1 \Delta _{com}^{\mathcal{O}(b_{k,l})} \<k_{k,l},l_{k,l}\> \right ) \Delta _{p_k}^{e_k}\<k_{x_k},l_{x_k}\> \right)(\prod_{l=n_0}^1\Delta _{com}^{\mathcal{O}(b_{0,l})}\<k_{0,l},l_{0,l}\> )\right)(H_m)$$
\end{footnotesize}
(as a braid).
\eth
\hfill
$\qed$

\brem
It is always possible to choose the point $x_0$ in such a way that $B \cap [x_1,x_0]=\emptyset$. Hence, we may assume that in $t_0$ there are no points that belong to the set $B$.
\erem

\brem
The exact location of the points in the set $N \cup B$ is not important for the algorithm. What is important, however, is the order of the points with respect to each other.
\erem

\subsection{The algorithm}

Now when the theory was established we can formulate the algorithm for computing the braid monodromy of a curve $S$. It's main idea is to find all the critical points which effect the braid monodromy results, the singular points of the curve and the points where complex switching occur. Then, using the computed information take each element of the generalized skeleton and compute for the result of the action induced by the diffeomorphisms conjugations as formulated in Theorem \ref{The algorithm theorem}. The product of the results as a product of braids is the braid monodromy for this generator of the fundamental group. The ordered set $\varphi(\Gamma _j)$ ($j=1,\cdots ,n$) define, as a product, the \emph{braid monodromy factorization} which describes the braid monodromy homomorphism.

\balg {Algorithm for computing the braid monodromy of a curve $S$} \label{The algorithm}
\ben
\item
Find in $S$ all the singular points with respect to the generic projection $\pi _1$.
\item
Compute for each point found in step $1$ its projection to the first coordinate. The set of all projections is the set $N$.
\item
Compute for each point $x_j \in N$ its Lefschetz pair $\mathcal{L}(x_j)=(k_{x_j},l_{x_j})$.
\item
Find for each $x_j \in N$ the diffeomorphism $\Delta _{p_j}^{e_j}$ induced by $x_j$ on the model by Proposition \ref{Model proposition}.
\item
Find in $S$ all the complex switching points with respect to $\pi _1$.
\item
Compute for each point found in step $4$ its projection to the first coordinate. The set of all projections is the set $B$.
\item
Compute for each point $b_{k,l} \in B$ its complex Lefschetz pair $\mathcal{L}(b_{k,l})=(k_{k,l},l_{k,l})$, and find the orientation $\mathcal{O}(b_{k,l})$ induced by the point using Definition \ref{orientation}.
\item
For each point $x_j$ find using Definition \ref{skeletons} the generalized skeleton $\prod _m H_m$ defining its local braid monodromy.
\item
For each point $x_j$ compute the product which is the braid monodromy result on $\Gamma _j$ using the generalized skeleton from Step 7 and Theorem \ref{The algorithm theorem}:
\begin{footnotesize}
$$\varphi(\Gamma _j)=\prod_m\left(\prod _{k=j-1}^1 \left ( \left (\prod_{l=n_k}^1 \Delta _{com}^{\mathcal{O}(b_{k,l})} \<k_{k,l},l_{k,l}\> \right ) \Delta _{p_k}^{e_k}\<k_{x_k},l_{x_k}\> \right)(\prod_{l=n_0}^1\Delta _{com}^{\mathcal{O}(b_{0,l})}\<k_{0,l},l_{0,l}\> )\right)(H_m)$$
\end{footnotesize}
\een
\ealg

In the next Section we give examples for the computation of braid monodromy using the algorithm.

%%%%%%%%%%%%%%%%%%%%%%%%%%%%%%%%%%%%%%%%%%%%%%%%%%%%%%%%%%%%%%%%%%%%%%%%%%%%%%%%%%%%%%%%%%%%%%%%
%%%%%%%%%&&&%%%%%%%%%%%%%%%%%%%%%%%
\section{Two examples for the algorithm use}
\subsection{The computation of the braid monodromy of a cusp}
In order to demonstrate the algorithm we compute the braid monodromy for the image of the curve $S$ defined by the equation $y^2-x^3=0$ under a linear transformation. In order to properly modify the curve using a linear transformation we take it's homoganization to $\C\PP^2$ which is defined by the homogenous polynomial: $y^2w-x^3$.
Then, we use the matrix 
$\left (
\begin{array}{ccc}
  5 & 3 & 1 \\
  3 & 2 & 1 \\
  2 & 1 & 1\\
\end{array}
\right )
$
defining a linear transformation, to look at the image of the curve under this transformation, which is the curve defined by ${(3x+2y+w)}^2(2x+y+w)-(5x+3y+w)^3=0$. In order to compute it's braid monodromy we choose generically a line at infinity and look at the curve's affine part which is ${(3x+2y+1)}^2(2x+y+1)-(5x+3y+1)^3=0$.

Following steps 1 and 2 of Algorithm \ref{The algorithm}, we find the points of the set $N$. Figure \ref{curve1} shows the real part of the curve with it's singular points identified. This also serve us in Step 2 so we know the order of the projections to the first coordinate of the singular points of $S$.
\myfigure{
\begin{center}
\includegraphics[scale=0.35,angle=270]{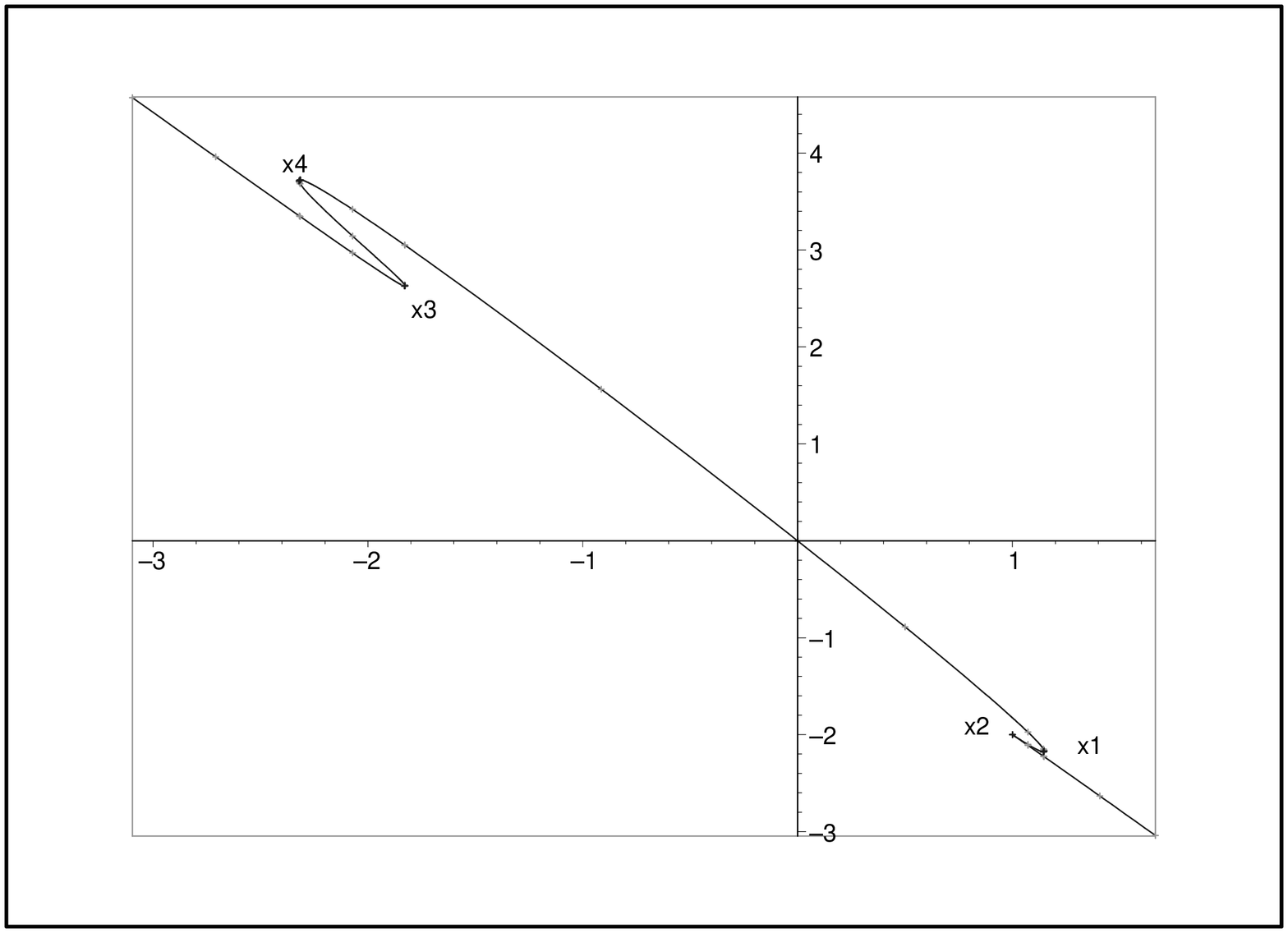}
\caption{The curve defined by ${(3x+2y+1)}^2(2x+y+1)-(5x+3y+1)^3=0$.}\label{curve1}
\end{center}
}

By Step 3, we compute the Lefschetz pair for each one of the singular points.
$\mathcal{L}(x_1)=(2,3)$, $\mathcal{L}(x_2)=(1,2)$, $\mathcal{L}(x_3)=(1,2)$, $\mathcal{L}(x_4)=(2,3)$.

In Step 4, we find the diffeomorphisms induced by each point depending on it's Lefschetz pair and singularity type: $\Delta _{a_2}^{\frac{1}{2}}\<2,3\>, \Delta _{d_1}^{\frac{3}{2}}\<1,2\>, \Delta _{a_2}^{\frac{1}{2}}\<1,2\>, \Delta _{a_1}^{\frac{1}{2}}\<2,3\>$ for the points $x_1,x_2,x_3,x_4$ respectively.

Note that since our curve is of degree $3$, we could have at most 2 conjugate complex points and therefore, $B=\emptyset$. Hence, Step 5,6 and 7 of the algorithm are trivial.

In Step 8, we find the generalized skeleton for each singular point and so we may summarize what we have until now in Table \ref{tab}.

\mytable{
\begin{center}
\begin{tabular}{||c|c|c|c|c||}
\hline
\hline
& & & & \\
   Point & Type & $\mathcal{L}$ &$\prod _m H_m$ & ${\Delta _{xj} }$  \\
& & & & \\
\hline
$x_1$ & $a_2$ & $\mathcal{L}(x_1)=\<2,3\>$ & {\includegraphics[scale=0.2]{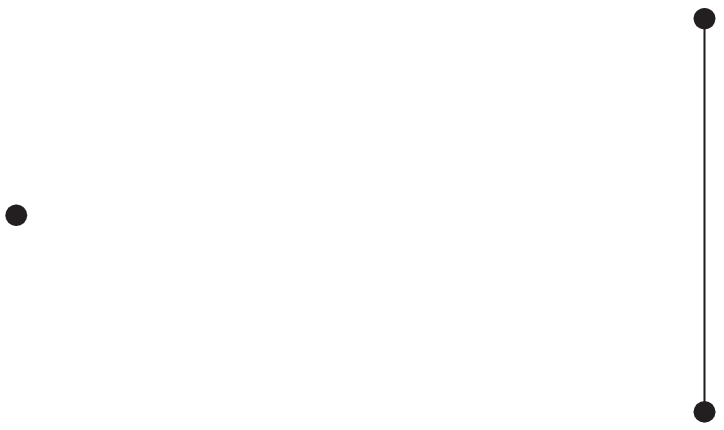}}  & ${\Delta _{a1}^{1/2}  \< 2,3 \> }$  \\
& & & & \\
\hline
$x_2$ & $d_1$ & $\mathcal{L}(x_2)=\<1,2\>$ & {\includegraphics[scale=0.2]{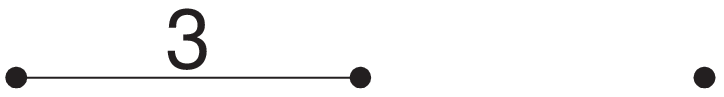}}  & ${\Delta _{d1}^{3/2}  \< 1,2 \> }$  \\
& & & & \\
\hline
$x_3$ & $a_2$ & $\mathcal{L}(x_3)=\<1,2\>$ & {\includegraphics[scale=0.2]{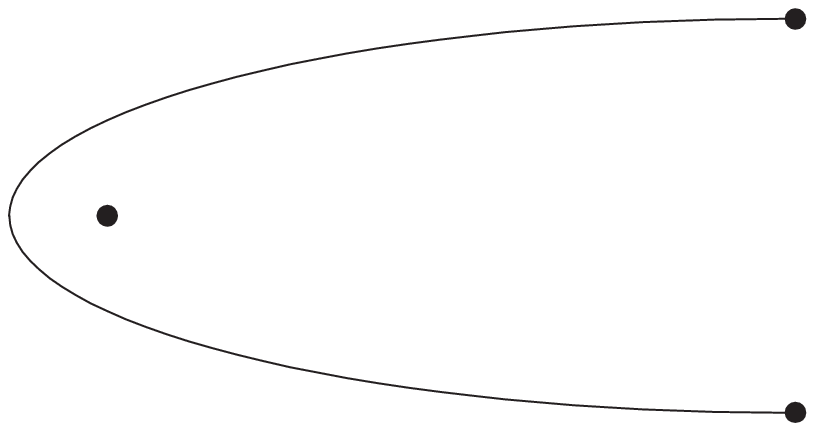}} & ${\Delta _{a1}^{1/2}  \< 1,2 \> }$  \\
& & & & \\
\hline
$x_4$ & $a_1$ & $\mathcal{L}(x_4)=\<2,3\>$ & {\includegraphics[scale=0.2]{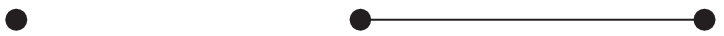}}  & ${\Delta _{a1}^{1/2}  \< 2,3 \> }$  \\
& & & & \\
\hline
\hline
\end{tabular}
\caption{Results of steps 1--8 of Algorithm \ref{The algorithm} on the curve ${(3x+2y+1)}^2(2x+y+1)-(5x+3y+1)^3=0$}\label{tab}
\end{center}}

Finely, we compute the braid monodromy using Step 9:\\
Note, that all generalized skeletons consists of only one half-twist (besides the generalized skeleton of the point $x_2$ which is a half-twist taken to the power $3$). Now, for each $j$ we need to conjugate the generalized skeleton of $x_j$ by all $\Delta _{p_k}^{e_k}\<\mathcal{L}(x_k)\>$, where $k=j-1,\cdots , 1$. This is best summarized Table \ref{BM cusp} (See Appendix B). Each part of table is used for the computation of the braid monodromy of a different $\Gamma _j$. First line in the is the generalized skeleton half-twists determined by the singular point which it's product represents the local braid monodromy at the point $x_j$. The left side of each line (besides the first) is the diffeomorphism which acts on the generalized skeleton elements, where the right side of each line is the result of this action. Finally, taking the product of the generalized skeleton elements at the last line is the result of the braid monodromy of $\Gamma _j$.

\myfigure{
\begin{center}
\includegraphics[scale=0.6]{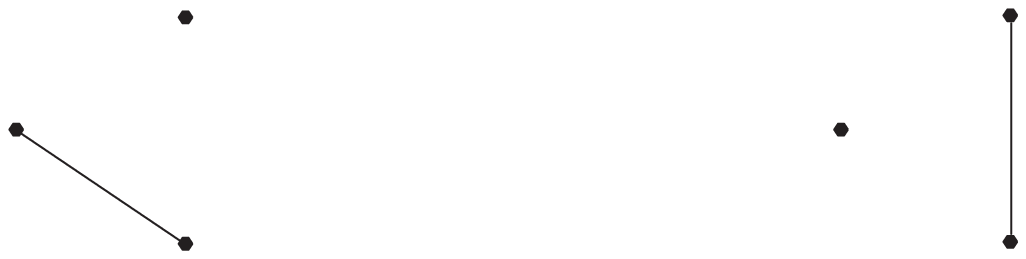} 
\caption{$\sigma _1$ (left) and $\sigma _2$ (right), the generators of $B_3$ in the model above $x_0$.}\label{generators1}
\end{center}
}

If the braid group of the model of the fiber above $x_0$ is generated by the two half-twists $\sigma _1$ and $\sigma _2$ as shown in Figure \ref{generators1}, we conclude that:\\
$\varphi(\Gamma _1)=\sigma _2$ \\
$\varphi(\Gamma _2)=\sigma _1^3$ \\
$\varphi(\Gamma _3)=\sigma _1$ \\
$\varphi(\Gamma _4)=\sigma _1^{-2}\sigma _2\sigma_1^2$ \\

\subsection{Example for the use of the algorithm when $B \neq \emptyset$}

For the second example consider the curve $S$ defined by the polynomial $(x^2-(y+\frac{3}{4})^2-1)((x-\frac{1}{2})^2-y^2-1)$.
We follow the steps of Algorithm \ref{The algorithm} in order to compute it's braid monodromy. \\
In Step 1, we compute all the singular points with respect to the projection $\pi _1$, and in Step 2, we project these points to get the set $N$. 
In our case the set $N$ consists of $6$ points as can be seen in Figure \ref{2ndExamplefig}.
\myfigure{
\begin{center}
\includegraphics[scale=0.35,angle=270]{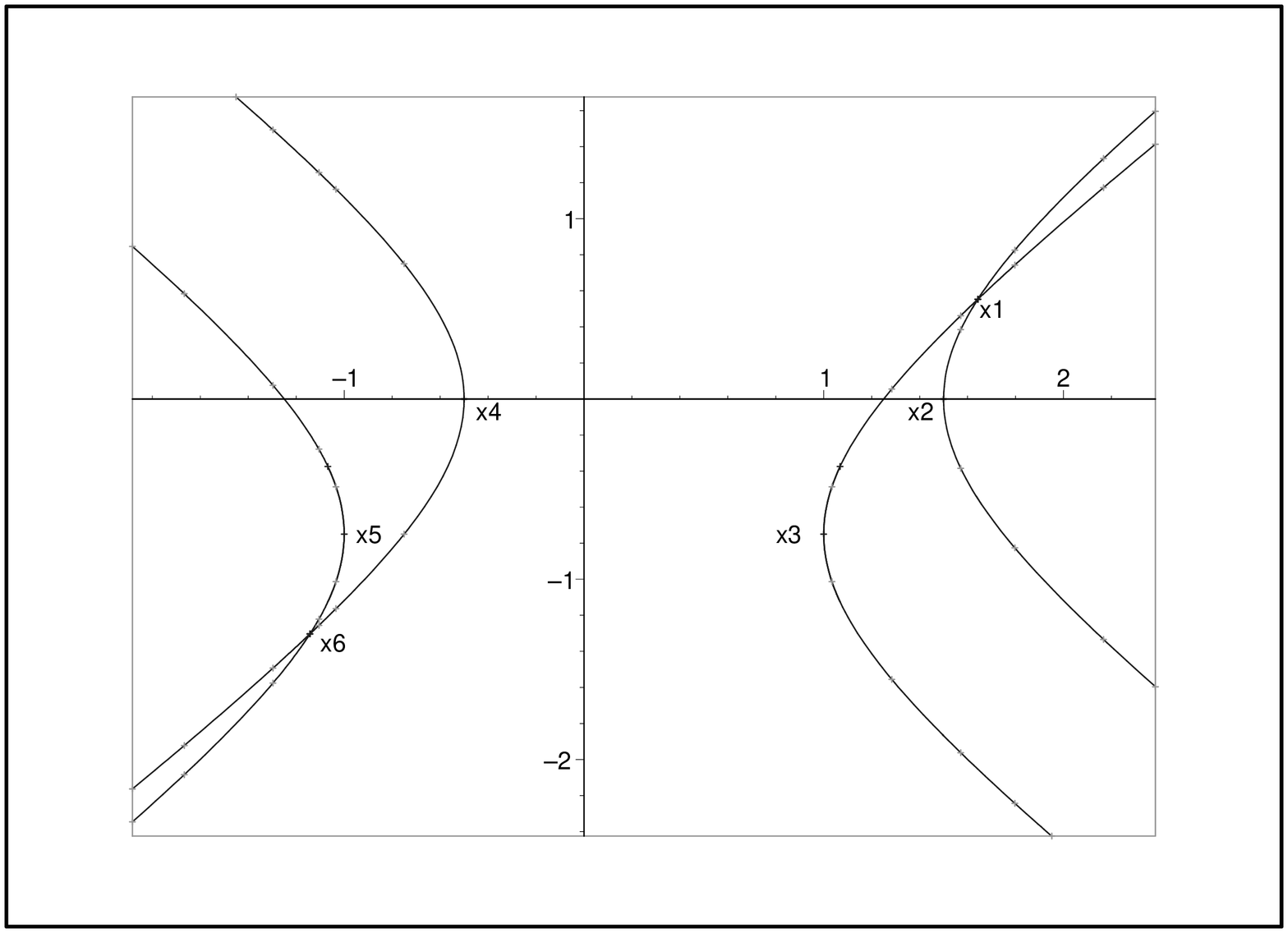} 
\caption{The real part of the curve defined by $(x^2-(y+\frac{3}{4})^2-1)((x-\frac{1}{2})^2-y^2-1)$.}\label{2ndExamplefig}
\end{center}
}
Next step in the algorithm (3) is to compute the Lefschetz pair associated with each point found in Steps 1 and 2. We get $\mathcal{L}(x_1)=(3,4)$, $\mathcal{L}(x_2)=(2,3)$, $\mathcal{L}(x_3)=(1,2)$, $\mathcal{L}(x_4)=(1,2)$, $\mathcal{L}(x_5)=(2,3)$, $\mathcal{L}(x_6)=(1,2)$.\\

In Step 4, we find the diffeomorphisms induced by each singular point: $\Delta _c\<2,3\>$, $\Delta _{a_1}^{\frac{1}{2}}\<2,3\>$, $\Delta _{a_1}^{\frac{1}{2}}\<1,2\>$, $\Delta _{a_2}^{\frac{1}{2}}\<1,2\>$, $\Delta _{a_2}^{\frac{1}{2}}\<2,3\>$, $\Delta _c\<1,2\>$ for the points $x_1,\cdots ,x_6$ respectively.

By Step 5, we find all the complex switching points with respect to $\pi _1$. This can be done for example by solving the equations of the curve to $y$, isolating the $\Im$ coordinate, and finding real $x$ for which we have equalities. It is also evident from Figure \ref{Points in B} which shows only the $\Im$ coordinates of the equation above the real axis, that we have one point in $B$ which is positioned between $x_3$ and $x_4$. Hence, we have completed Step 5 and 6 and we have the set $B=\{b_{3,1}=b\}$. Moreover, since there are at most $4$ complex points at any fiber above the real axis, we know that the complex Lefschetz pair associated with the point $b$ is $\mathcal{L}(b)=\<1,2\>$. As for the orientation induced by the point $b$ we follow Definition \ref{orientation}. To the right of $b$ the points $n+1\pm i$ in the model belong to the imaginary points with smaller real coordinate, and the points $n+1\pm 2i$ belong to the imaginary points with larger real coordinate. To the left of $b$ the points $n+1\pm i$ in the model belong to the imaginary points with larger real coordinate, and the points $n+1\pm 2i$ in the model belong to the imaginary points with the smaller real coordinate. Hence, the orientation induced by $b$ is $\mathcal{O}(b)=1$, and so we have completed Step 7 of the algorithm.

\myfigure{
\begin{center}
\includegraphics[scale=0.35,angle=270]{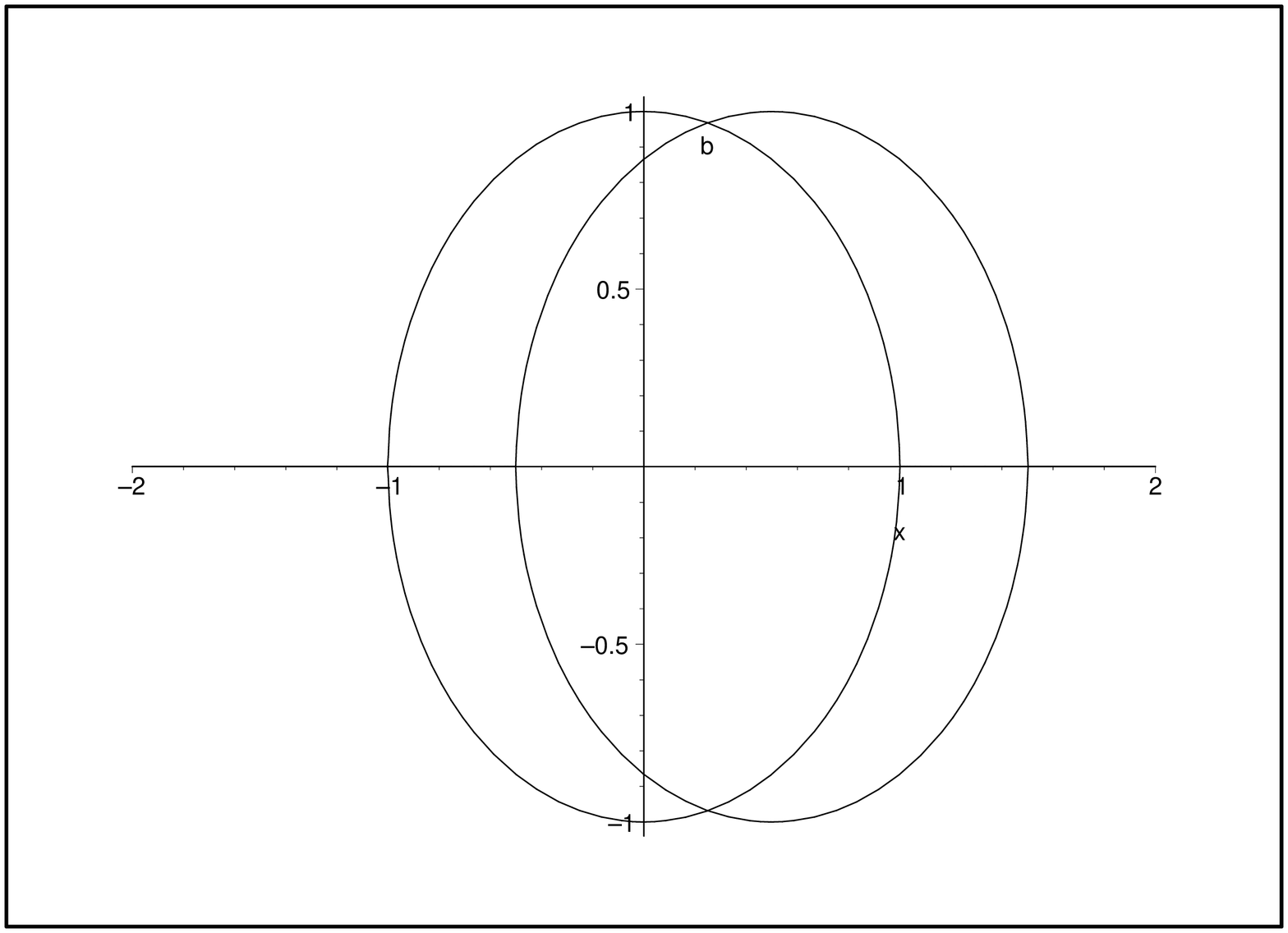} \label{Points in B} 
\caption{The $\Im$ coordinates of the curve defined by $(x^2-(y+\frac{3}{4})^2-1)((x-\frac{1}{2})^2-y^2-1)$ above the real axis.}
\end{center}
}

In Step 8, we find the generalized skeletons induced by each of the singular points found in Step 1. 
Results of the steps of the algorithm up to now are summarized in Table \ref{2ndtab}.

\mytable{
\begin{center}
\begin{tabular}{||c|c|c|c|c||}
\hline
\hline
& & & & \\
   Point & type & $\mathcal{L}$ & $\prod _m H_m$ & ${\Delta _{xj} }$ \\
& & & & \\
\hline
& & & & \\
$x_1$ & $c$ & $\mathcal{L}(x_1)=\<3,4\>$ & \includegraphics[,scale=0.2]{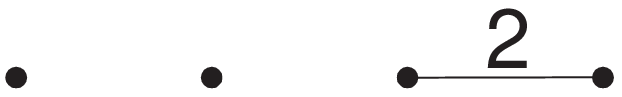}  & ${\Delta _{c} \< 3,4 \> }$ \\ 
& & & & \\
\hline
& & & & \\
$x_2$ & $a_1$ & $\mathcal{L}(x_2)=\<2,3\>$ & \includegraphics[,scale=0.2]{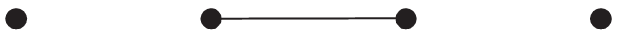} & ${\Delta _{a_1}^{\frac{1}{2}}\< 2,3 \> }$ \\ 
\hline
& & & & \\
$x_3$ & $a_1$ & $\mathcal{L}(x_3)=\<1,2\>$ & \includegraphics[,scale=0.2]{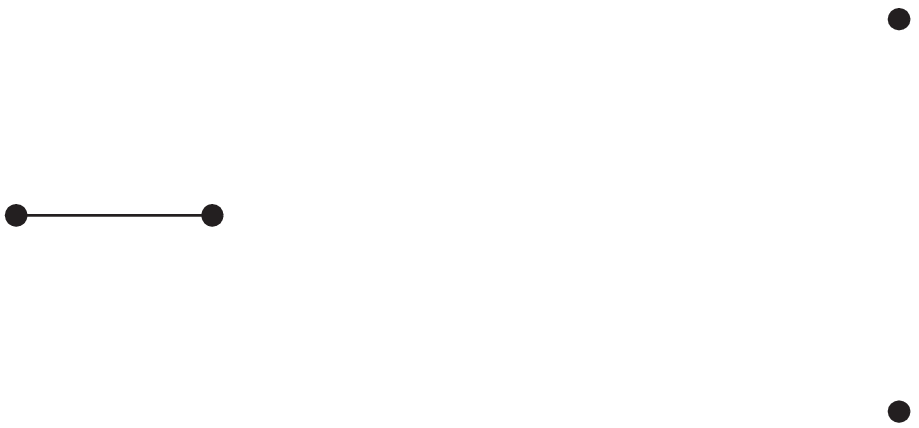} & ${\Delta _{a_1}^{\frac{1}{2}}\< 1,2 \> }$ \\ 
\hline
& & & & \\
$b$ & $com$ & $\mathcal{L}(b)=\<1,2\>$ & No generalized skeleton  & ${\Delta _{com}^1 \< 1,2 \> }$ \\
& & & & \\
\hline
& & & & \\
$x_4$ & $a_2$ & $\mathcal{L}(x_4)=\<1,2\>$ & \includegraphics[,scale=0.2]{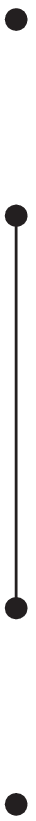} & ${\Delta _{a_2}^{\frac{1}{2}}\< 1,2 \> }$ \\ 
& & & & \\
\hline
& & & & \\
$x_5$ & $a_2$ & $\mathcal{L}(x_5)=\<2,3\>$ & \includegraphics[,scale=0.2]{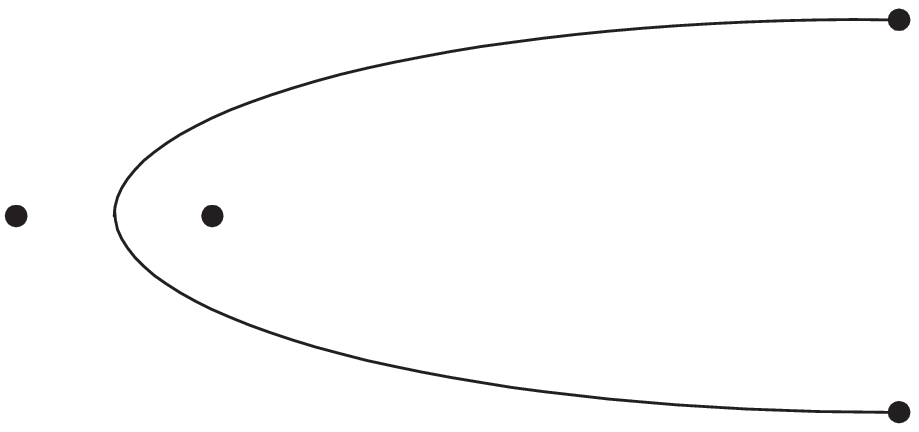} & ${\Delta _{a_2}^{\frac{1}{2}}\< 2,3 \> }$ \\ 
& & & & \\
\hline
& & & & \\
$x_6$ & $c$ & $\mathcal{L}(x_6)=\<1,2\>$ & \includegraphics[,scale=0.2]{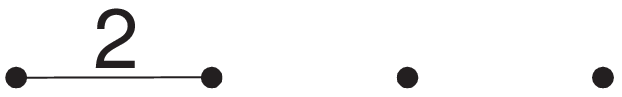} & ${\Delta _{c}\< 1,2 \> }$ \\ 
& & & & \\
\hline
\hline
\end{tabular}
\caption{Results of steps 1--8 of Algorithm \ref{The algorithm} on the curve $(x^2-(y+\frac{3}{4})^2-1)((x-\frac{1}{2})^2-y^2-1)$.}\label{2ndtab}
\end{center}}

\noindent
Finely, we compute the braid monodromy using Step 9:\\
For each $x_j$ we need to conjugate the elements constructing the generalized skeleton of $x_j$ by all $\Delta _{p_k}^{e_k}\<\mathcal{L}(x_k)\>$, where $k=j-1,\cdots ,1$, and $p_k$ is the singularity type of the point $x_k$ or the diffeomorphism induced by the point $b$. This is best summarized Table \ref{BM 2nd} (See Appendix B). In order to demonstrate the use of the algorithm we bring here only full computations for the points $x_5$ and $x_6$, for the rest of the points we give only the final result.

\myfigure{
\begin{center}
\includegraphics[scale=0.5]{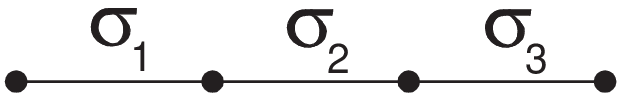} 
\caption{The generators of $B_4$ in the model above $x_0$}\label{generators2}
\end{center}
}

If the braid group of the model of the fiber above $x_0$ is generated by the three half-twists $\sigma _1,\sigma _2$ and $\sigma _3$ as shown in Figure \ref{generators2}, we conclude that:\\
$\varphi(\Gamma _1)=\sigma _3^2$ \\
$\varphi(\Gamma _2)=\sigma _3^{-1}\sigma_2\sigma_3$ \\
$\varphi(\Gamma _3)=\sigma _2^{-1}\sigma_1\sigma_2$ \\
$\varphi(\Gamma _4)=\sigma _2^{-1}\sigma_3\sigma_2$ \\
$\varphi(\Gamma _5)=\sigma_3^2\sigma_2^{-1}\sigma_1\sigma_2\sigma _3^{-2}$ \\
$\varphi(\Gamma _6)=\sigma _3^2$ \\

%%%%%%%%%&&&%%%%%%%%%%%%%%%%%%%%%%%

\clearpage

\section{Appendix A}

\myfigure{
\begin{center}
\includegraphics[scale=0.70]{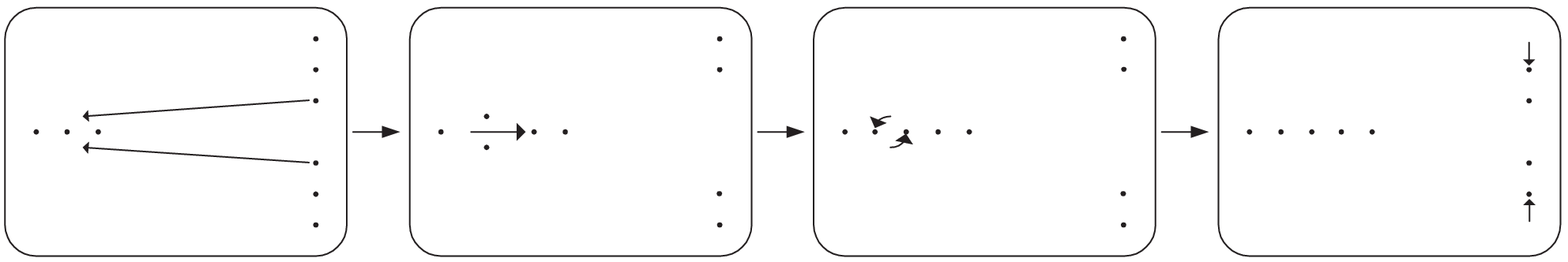}
\caption{$\Delta _{a_1}^{\frac{1}{2}}\left<k,l\right>$}\label{Delta a1}
\end{center}
}

\myfigure{
\begin{center}
\includegraphics[scale=0.70]{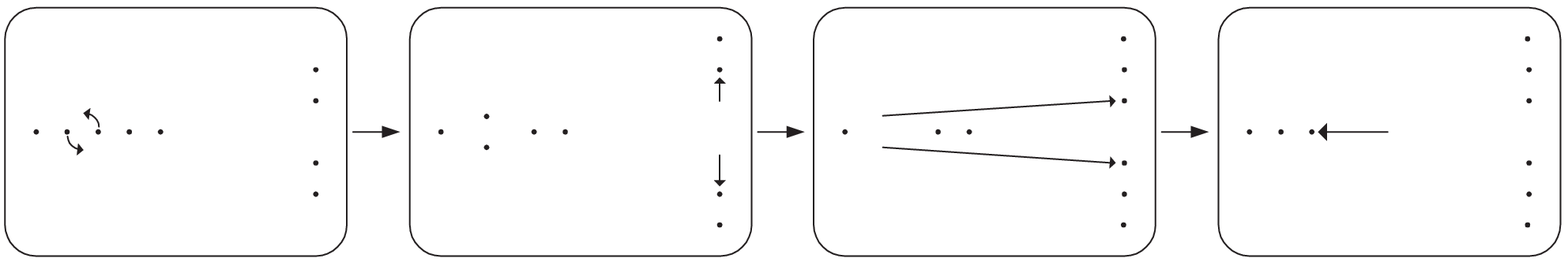}
\caption{$\Delta _{a_2}^{\frac{1}{2}}\left<k,l\right>$}\label{Delta a2}
\end{center}
}

\myfigure{
\begin{center}
\includegraphics[scale=0.70]{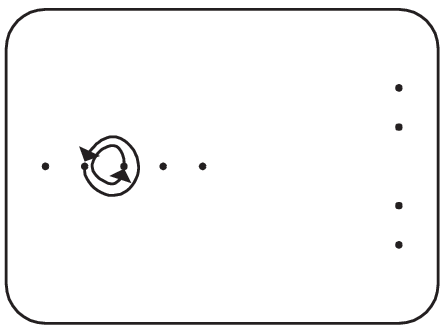}
\caption{$\Delta _b^2\left<k,l\right>$}\label{Delta b}
\end{center}
}

\myfigure{
\begin{center}
\includegraphics[scale=0.70]{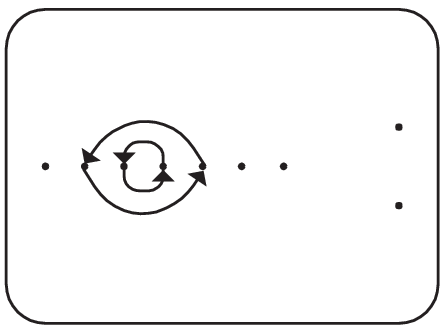}
\caption{$\Delta _c\left<k,l\right>$}\label{Delta c}
\end{center}
}

\clearpage
\myfigure{
\begin{center}
\includegraphics[scale=0.70]{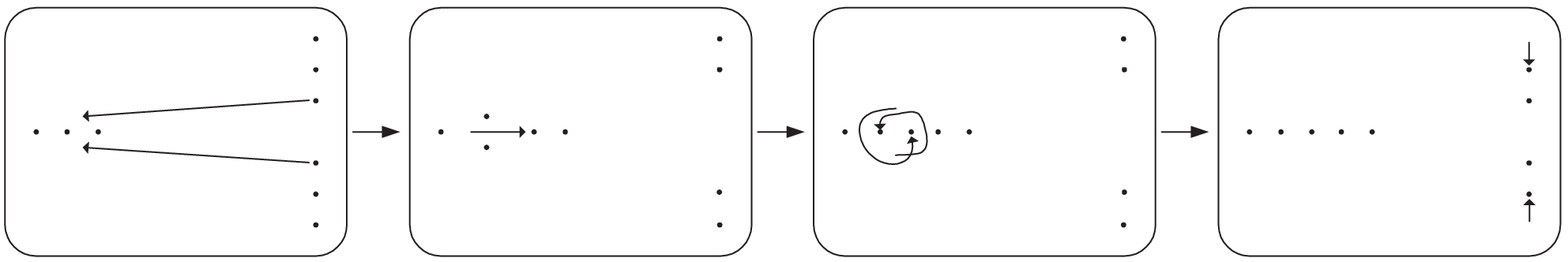}
\caption{$\Delta _{d_1}^{\frac{\nu}{2}}\left<k,l\right>$ with $\nu$ odd}\label{Delta d1 odd}
\end{center}
}

\myfigure{
\begin{center}
\includegraphics[scale=0.70]{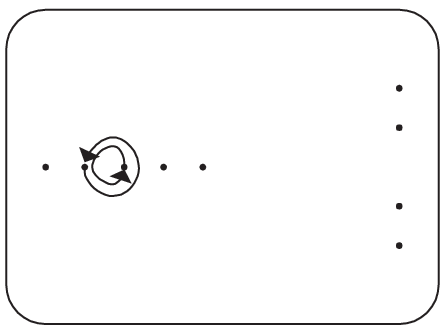}
\caption{$\Delta _{d_1}^{\frac{\nu}{2}}\left<k,l\right>$ with $\nu$ even}\label{Delta d1 even}
\end{center}
}

\myfigure{
\begin{center}
\includegraphics[scale=0.70]{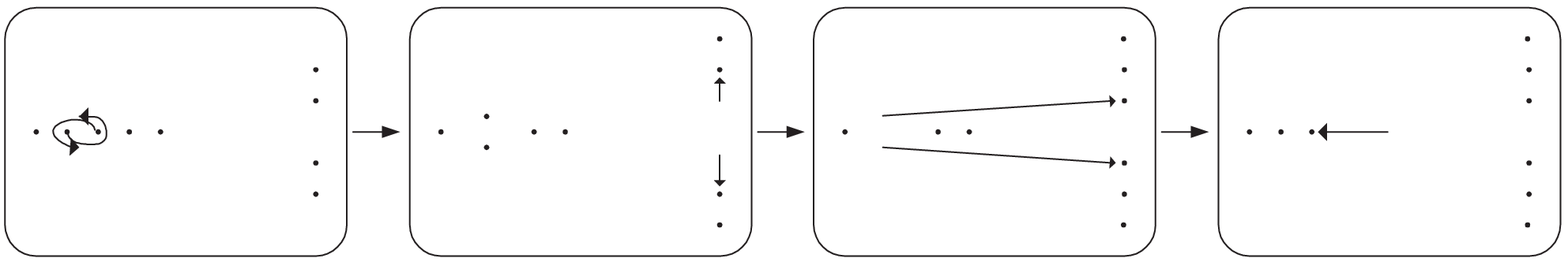}
\caption{$\Delta _{d_2}^{\frac{\nu}{2}}$ with $\nu$ odd}\label{Delta d2 odd}
\end{center}
}

\myfigure{
\begin{center}
\includegraphics[scale=0.70]{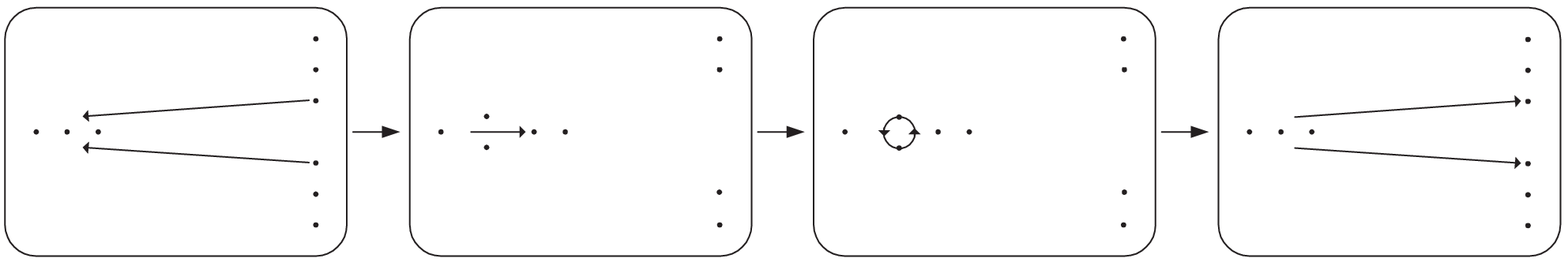}
\caption{$\Delta _{d_2}^{\frac{\nu}{2}}$ with $\nu$ even}\label{Delta d2 even}
\end{center}
}

\myfigure{
\begin{center}
\includegraphics[scale=0.70]{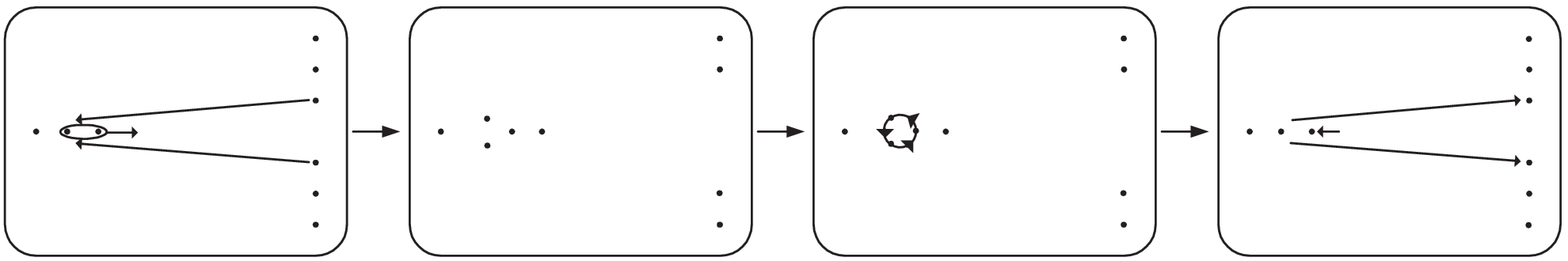}
\caption{$\Delta _{d_3}^{\frac{2}{3}}\left<k,l\right>$}\label{Delta d3}
\end{center}
}

\myfigure{
\begin{center}
\includegraphics[scale=0.70]{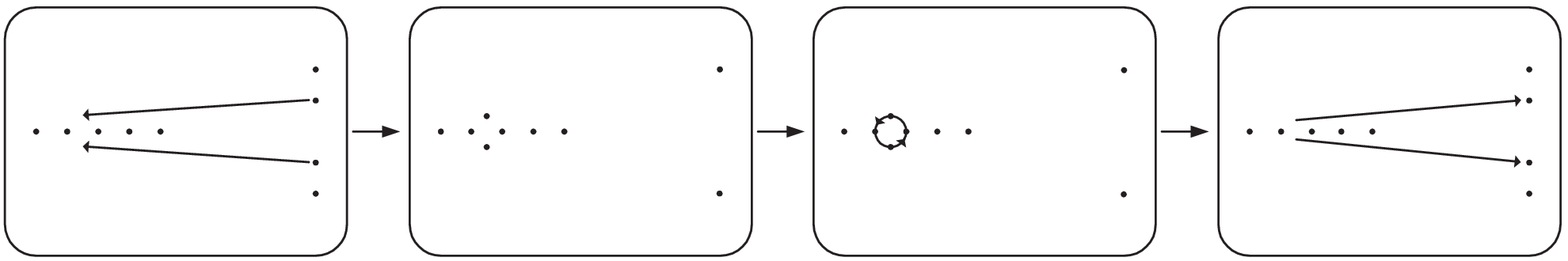}
\caption{$\Delta _{e_1}^{\frac{3}{2}}\left<k,l\right>$}\label{Delta e1}
\end{center}
}

\myfigure{
\begin{center}
\includegraphics[scale=0.70]{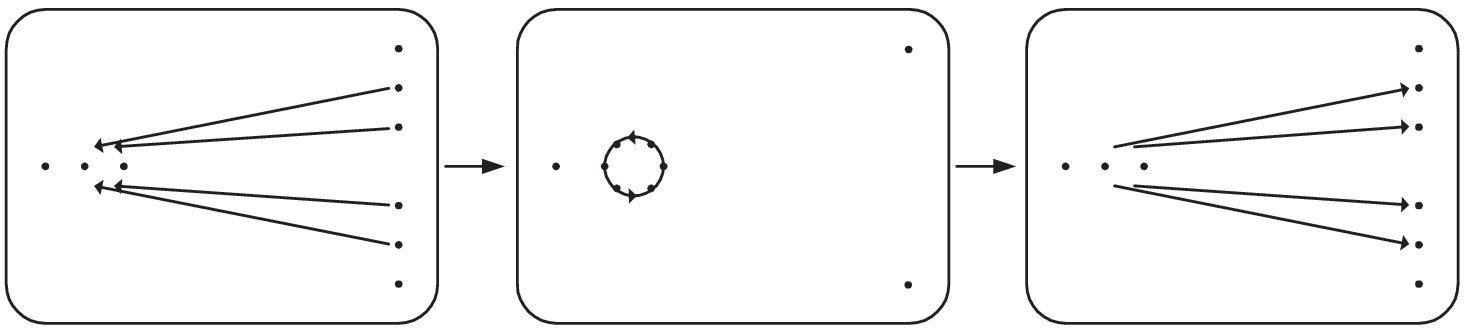}
\caption{$\Delta _{e_2}^{\frac{2}{3}}\left<k,l\right>_{\<a<b\>}$}\label{Delta e2 a<b}
\end{center}
}

\myfigure{
\begin{center}
\includegraphics[scale=0.70]{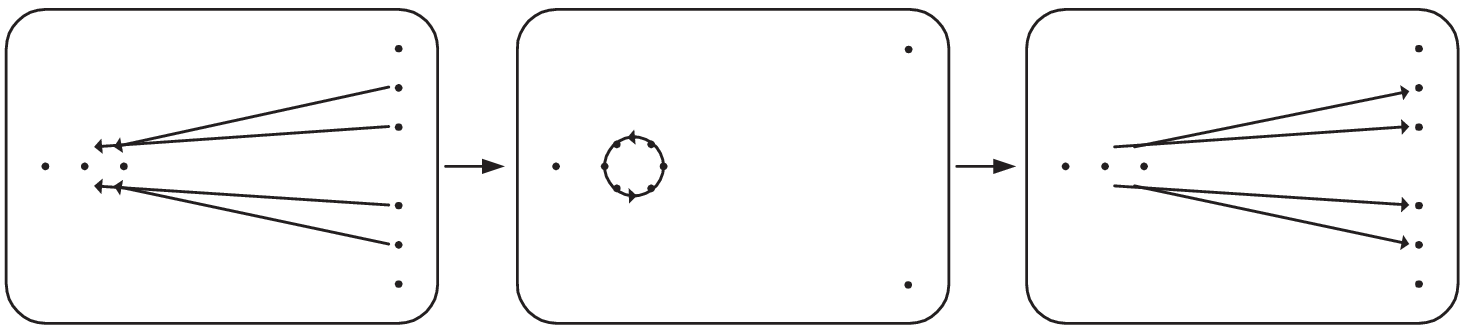}
\caption{$\Delta _{e_2}^{\frac{2}{3}}\left<k,l\right>_{\<a>b\>}$}\label{Delta e2 a>b}
\end{center}
}

\myfigure{
\begin{center}
\includegraphics[scale=0.70]{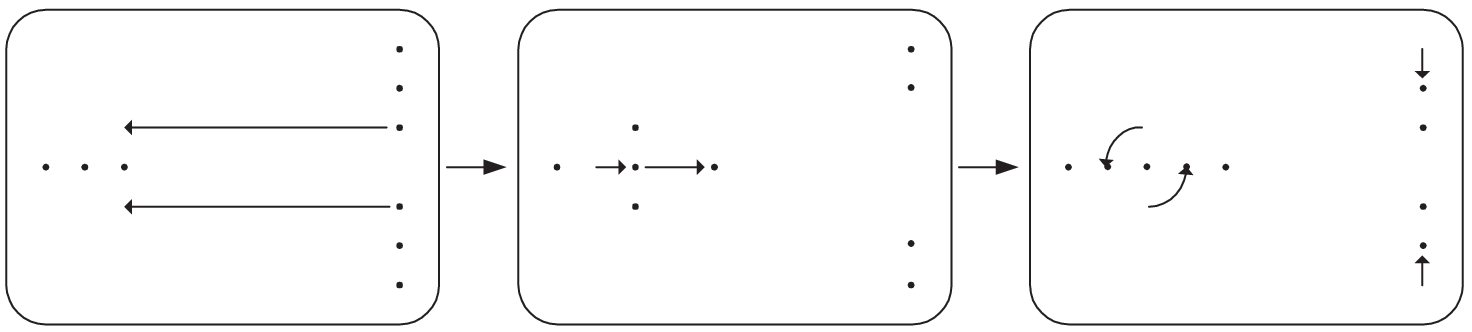}
\caption{$\Delta _{f_1}^{\frac{\nu}{2}}\left<k,l\right>$ with $\nu$ odd}\label{Delta f1 odd}
\end{center}
}

\myfigure{
\begin{center}
\includegraphics[scale=0.70]{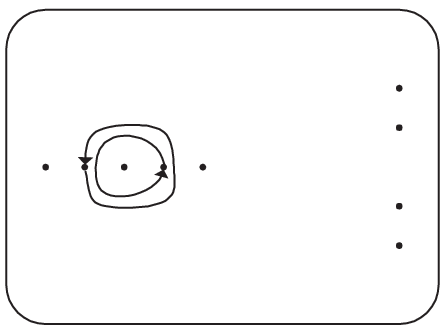}
\caption{$\Delta _{f_1}^{\frac{\nu}{2}}\left<k,l\right>$ with $\nu$ even}\label{Delta f1 even}
\end{center}
}

\myfigure{
\begin{center}
\includegraphics[scale=0.70]{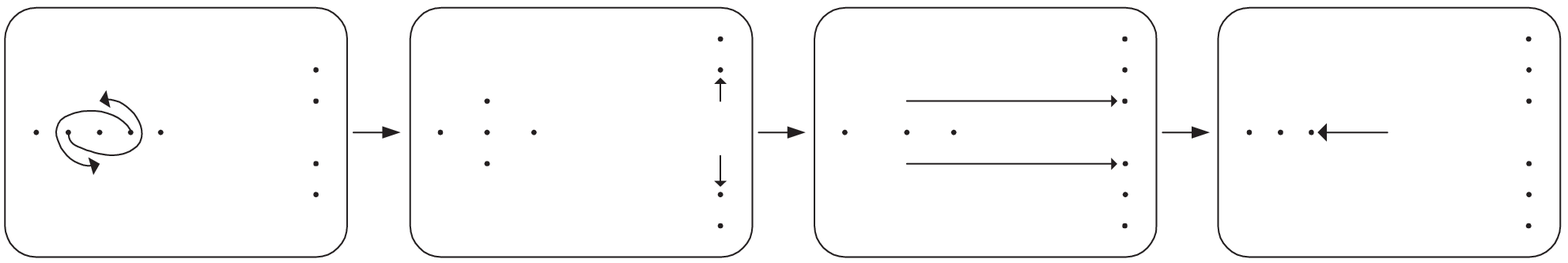}
\caption{$\Delta _{f_2}^{\frac{\nu}{2}}\left<k,l\right>$ with $\nu$ odd}\label{Delta f2 odd}
\end{center}
}

\myfigure{
\begin{center}
\includegraphics[scale=0.70]{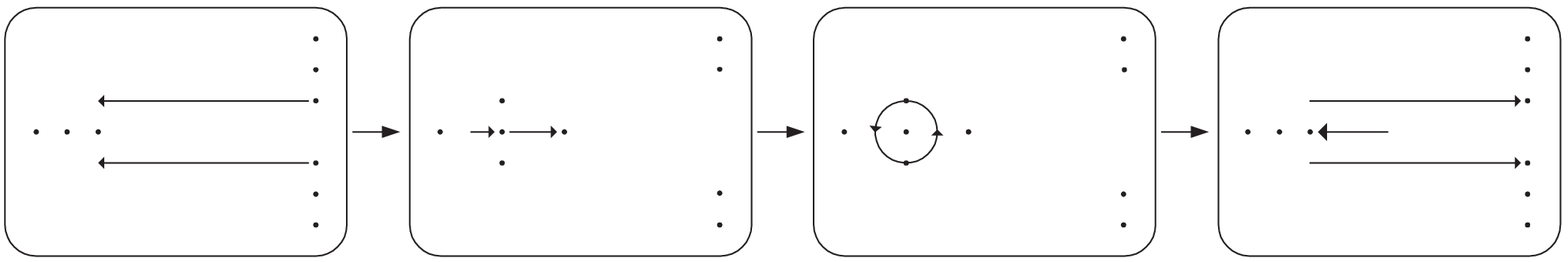}
\caption{$\Delta _{f_2}^{\frac{\nu}{2}}\left<k,l\right>$ with $\nu$ even}\label{Delta f2 even}
\end{center}
}

\myfigure{
\begin{center}
\includegraphics[scale=0.70]{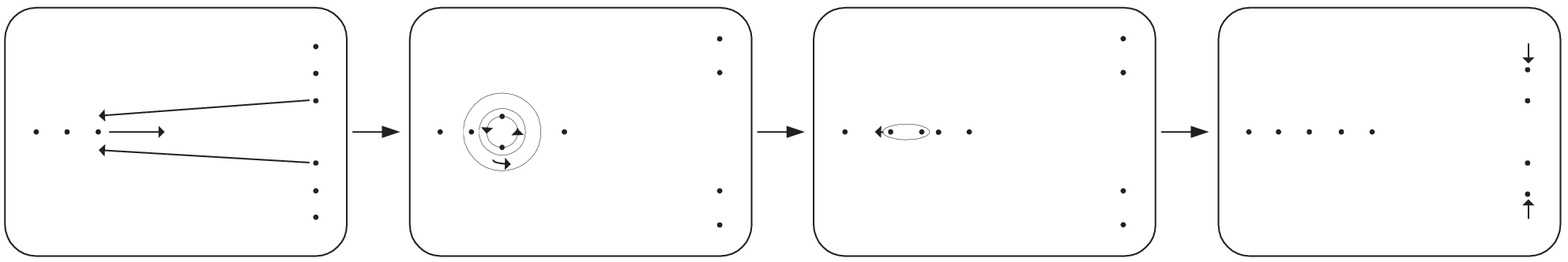}
\caption{$\Delta _{g_1}^{\frac{\nu}{2}}\left<k,l\right>$ with $\nu$ odd}\label{Delta g1 odd}
\end{center}
}

\myfigure{
\begin{center}
\includegraphics[scale=0.70]{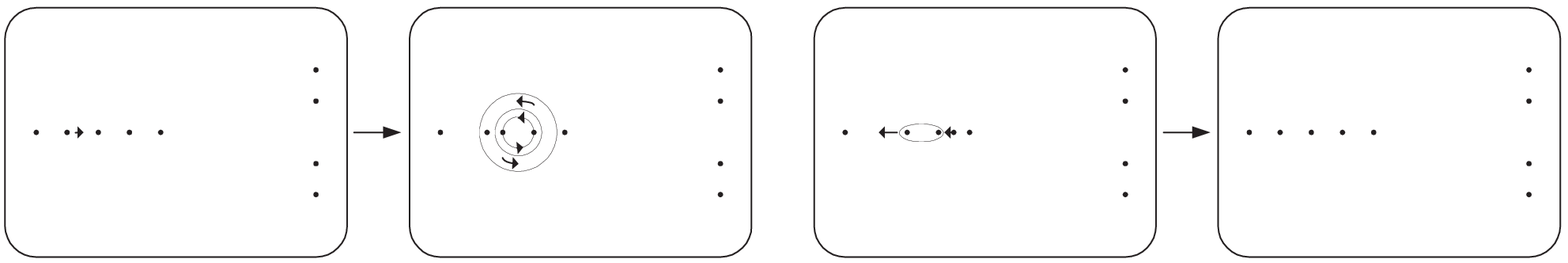}
\caption{$\Delta _{g_1}^{\frac{\nu}{2}}\left<k,l\right>$ with $\nu$ even}\label{Delta g1 even}
\end{center}
}

\myfigure{
\begin{center}
\includegraphics[scale=0.70]{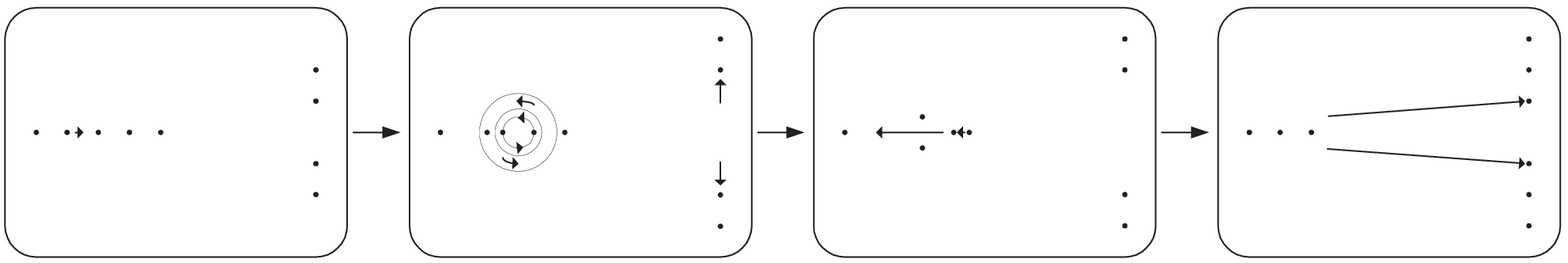}
\caption{$\Delta _{g_2}^{\frac{\nu}{2}}\left<k,l\right>$ with $\nu$ odd}\label{Delta g2 odd}
\end{center}
}

\myfigure{
\begin{center}
\includegraphics[scale=0.70]{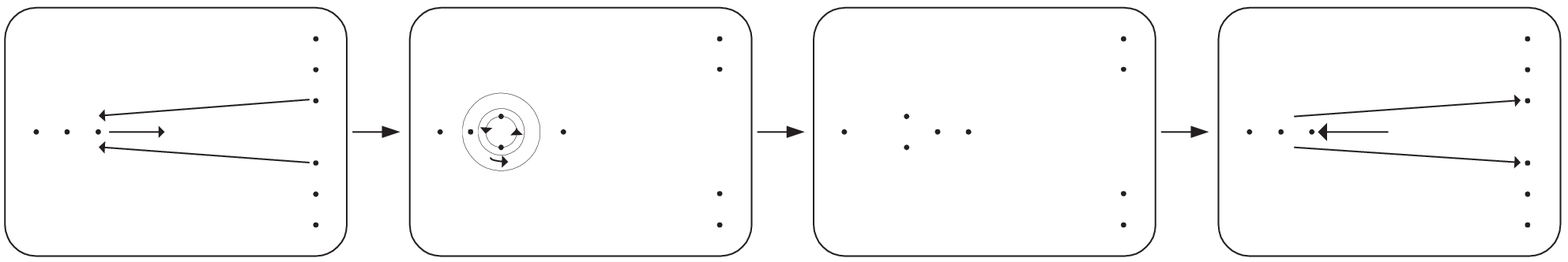}
\caption{$\Delta _{g_2}^{\frac{\nu}{2}}\left<k,l\right>$ with $\nu$ even}\label{Delta g2 even}
\end{center}
}

\myfigure{
\begin{center}
\includegraphics[scale=0.70]{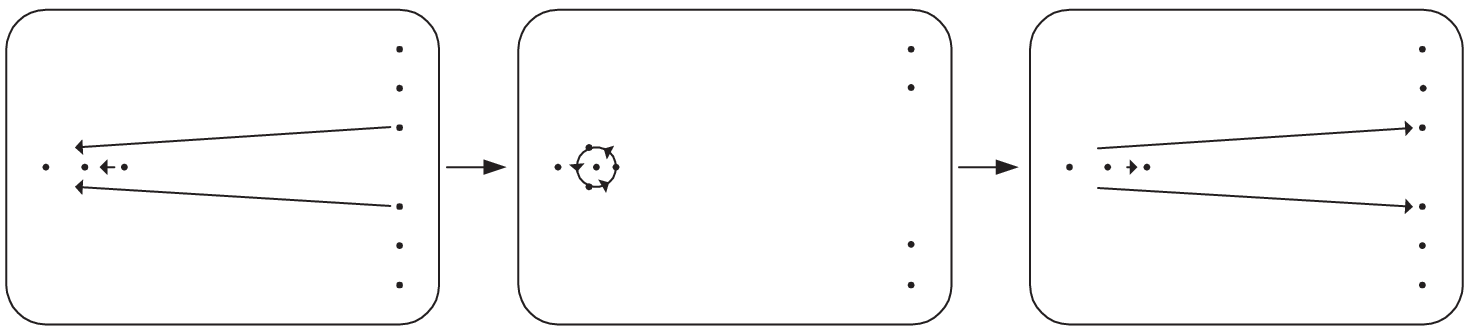}
\caption{$\Delta _{g_3}^{\frac{2}{3}}\left<k,l\right>$}\label{Delta g3}
\end{center}
}

\myfigure{
\begin{center}
\includegraphics[scale=0.70]{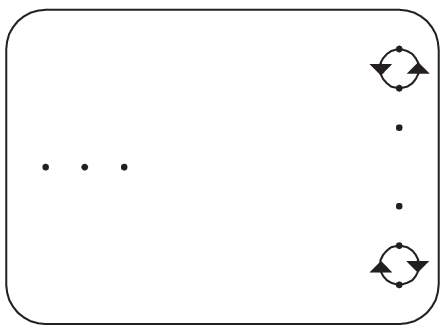}
\caption{$\Delta _{com}^{1}\left<k,l\right>$}\label{Delta com}
\end{center}
}

\clearpage

\section{Appendix B}

\mytable{
\begin{center}
\begin{tabular}{||c|c||c|c||}
\hline
\hline
\multicolumn{2}{||c||}{$\varphi(\Gamma _1)$}& \multicolumn{2}{|c||}{$\varphi(\Gamma _2)$}\\
\hline
Conj. diffeo. & Generalized skeleton & Conj. diffeo. & Generalized skeleton \\
\hline
 & & & \\
 & \includegraphics[scale=0.2]{1stGSX1.eps} & & \includegraphics[scale=0.2]{1stGSX2.eps} \\
\hline
 & & & \\
 & & ${\Delta _{a1}^{1/2}  < 2,3 > }$ & \includegraphics[scale=0.2]{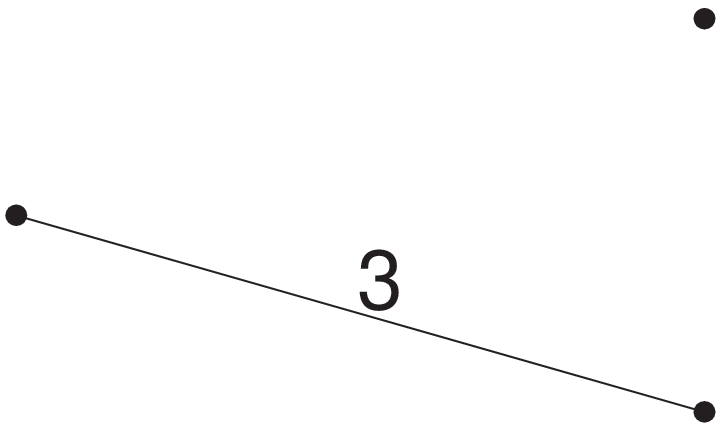} \\
\hline
\hline
\multicolumn{2}{||c||}{$\varphi(\Gamma _3)$}& \multicolumn{2}{|c||}{$\varphi(\Gamma _4)$}\\
\hline
Conj. diffeo. & Generalized skeleton & Conj. diffeo. & Generalized skeleton \\
\hline
 & & & \\
 & \includegraphics[scale=0.2]{1stGSX3.eps} & &\includegraphics[scale=0.2]{1stGSX4.eps}\\
\hline
 & & & \\
 ${\Delta _{d1}^{3/2}  < 1,2 > }$ & \includegraphics[scale=0.2]{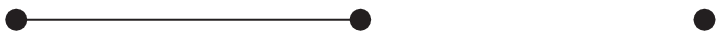} & ${\Delta _{a1}^{1/2}  < 1,2 > }$ & \includegraphics[scale=0.2]{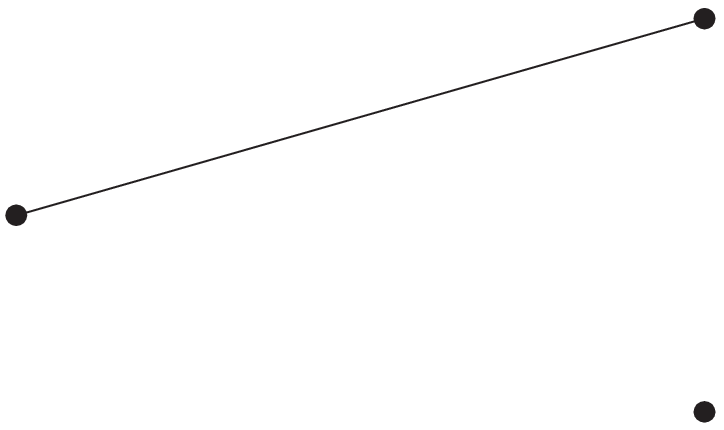} \\
\hline
 & & & \\
 ${\Delta _{a1}^{1/2}  < 2,3 > }$ & \includegraphics[scale=0.2]{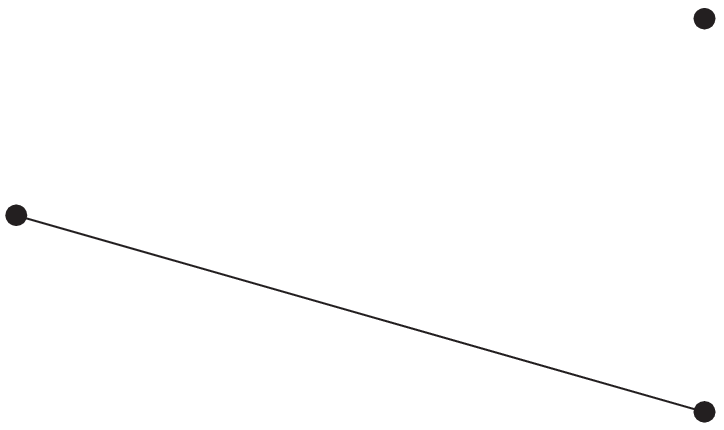} & ${\Delta _{d1}^{3/2}  < 1,2 > }$ & \includegraphics[scale=0.2]{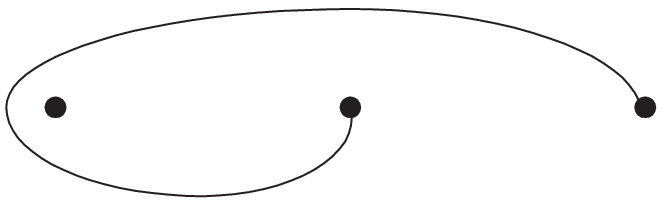} \\
\hline
 & & & \\
 & & ${\Delta _{a1}^{1/2}  < 2,3 > }$ & \includegraphics[scale=0.2]{1stGSX3.eps} \\
\hline
\hline
\end{tabular}
\caption{Braid monodromy for the cusp defined by the polynomial ${(3x+2y+1)}^2(2x+y+1)-(5x+3y+1)^3$ }\label{BM cusp}
\end{center}}

\mytable{
\begin{center}
\begin{tabular}{||c|c||c|c||}
\hline
\hline
\multicolumn{2}{||c||}{$\varphi(\Gamma _1)$}& \multicolumn{2}{|c||}{$\varphi(\Gamma _2)$}\\
\hline
Conj. diffeo. & Generalized skeleton & Conj. diffeo. & Generalized skeleton \\
\hline
 & & & \\
Final result & \includegraphics[scale=0.2]{2ndGSX1.eps} & Final result & \includegraphics[scale=0.2]{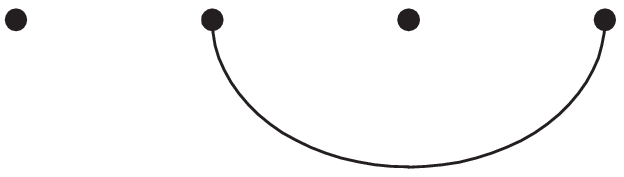} \\
\hline
\hline
\multicolumn{2}{||c||}{$\varphi(\Gamma _3)$}& \multicolumn{2}{|c||}{$\varphi(\Gamma _4)$}\\
\hline
Conj. diffeo. & Generalized skeleton & Conj. diffeo. & Generalized skeleton \\
\hline
 & & & \\
Final result & \includegraphics[scale=0.2]{2ndGSX3.eps} & Final result & \includegraphics[scale=0.2]{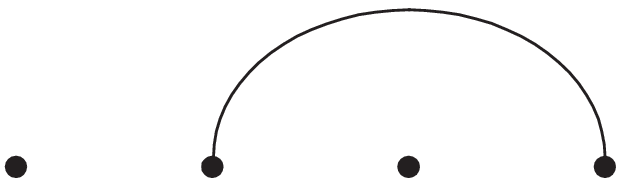} \\
\hline
\hline
\multicolumn{2}{||c||}{$\varphi(\Gamma _5)$}& \multicolumn{2}{|c||}{$\varphi(\Gamma _6)$}\\
\hline
Conj. diffeo. & Generalized skeleton & Conj. diffeo. & Generalized skeleton \\
\hline
 & & & \\
 & \includegraphics[scale=0.2]{2ndGSX5.eps} & &\includegraphics[scale=0.2]{2ndGSX6.eps}\\
\hline
 & & & \\
 ${\Delta _{a_2}^{\frac{1}{2}}  \< 1,2 \> }$ & \includegraphics[scale=0.2]{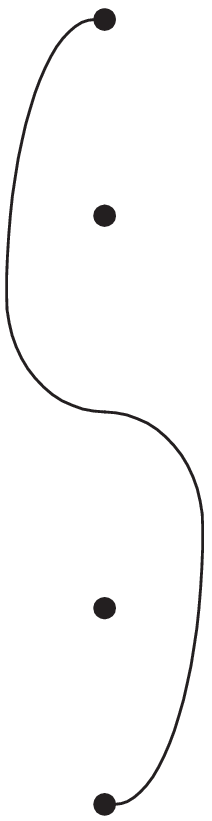} & ${\Delta _{a_2}^{\frac{1}{2}}  < 2,3 > }$ & \includegraphics[scale=0.2]{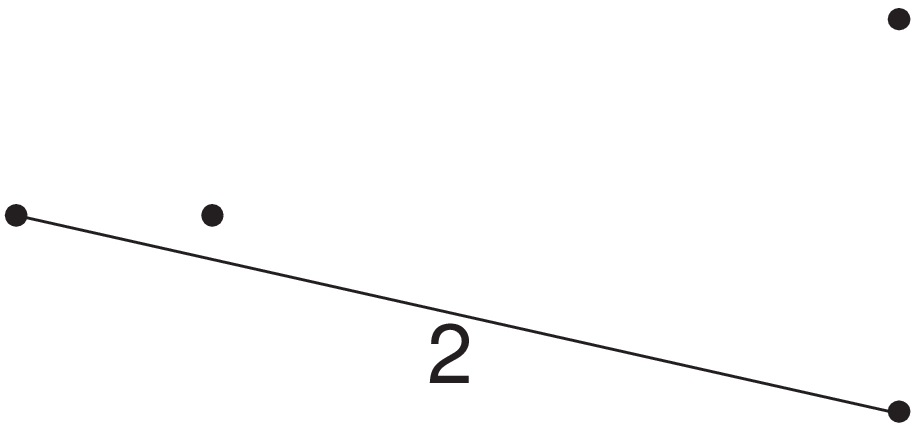} \\
\hline
 & & & \\
 ${\Delta _{com}  \< 1,2 \> }$ & \includegraphics[scale=0.2]{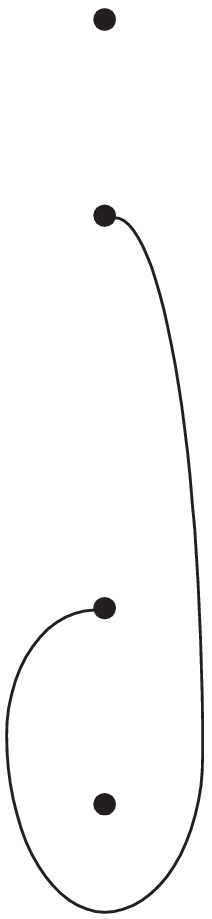} & ${\Delta _{a_2}^{\frac{1}{2}}  \< 1,2 \> }$ & \includegraphics[scale=0.2]{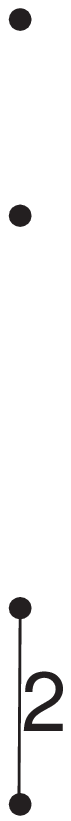} \\
\hline
 & & & \\
 ${\Delta _{a_1}^{\frac{1}{2}}  \< 1,2 \> }$ & \includegraphics[scale=0.2]{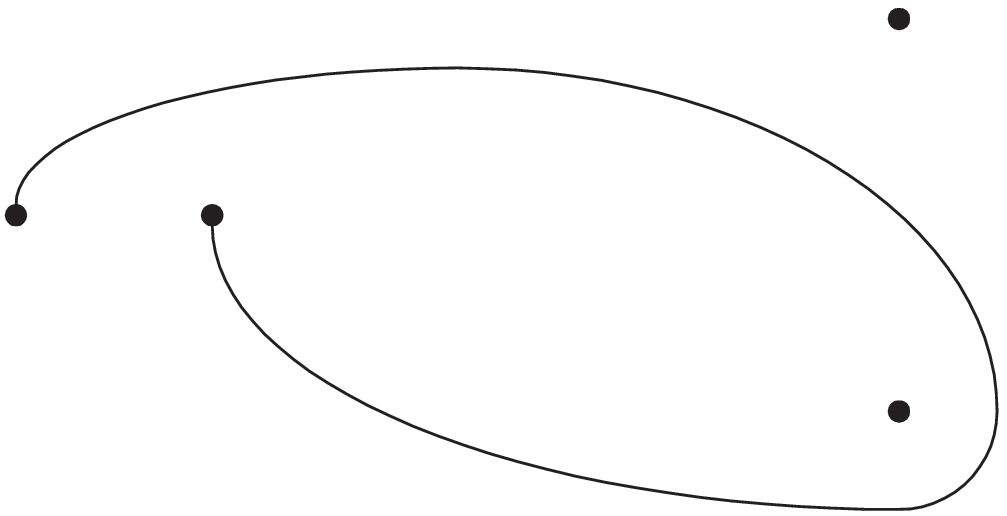} & ${\Delta _{com}  \< 1,2 \> }$ & \includegraphics[scale=0.2]{FinalX6x4.eps} \\
\hline
 & & & \\
 ${\Delta _{a_1}^{\frac{1}{2}}  \< 2,3 \> }$ & \includegraphics[scale=0.2]{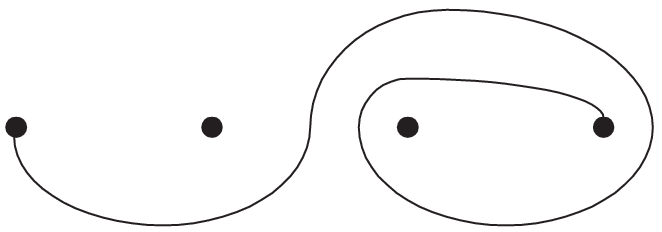} & ${\Delta _{a_1}^{\frac{1}{2}}  \< 1,2 \> }$ & \includegraphics[scale=0.2]{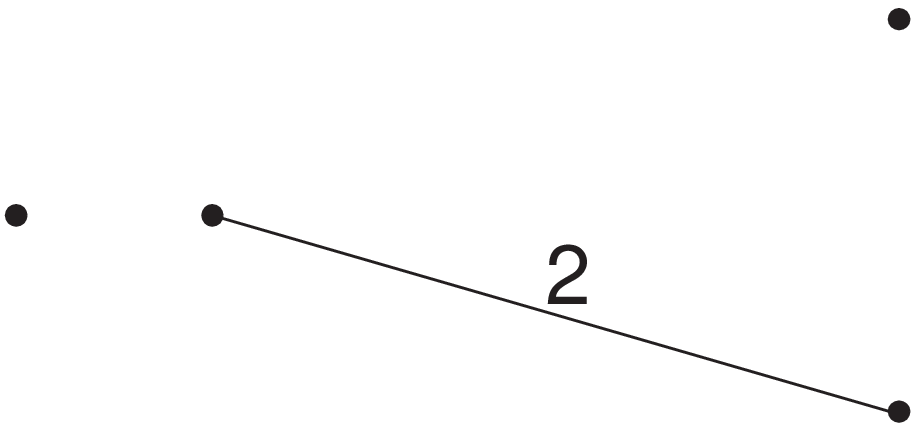} \\
\hline
 & & & \\
 ${\Delta _{c}  \< 3,4 \> }$ & \includegraphics[scale=0.2]{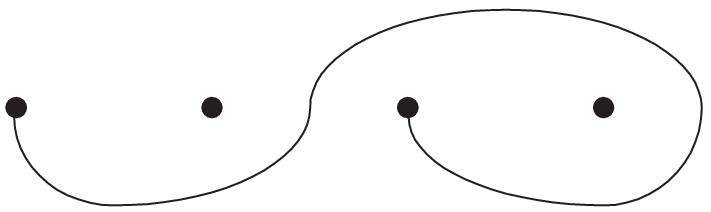} & ${\Delta _{a_1}^{\frac{1}{2}}  \< 2,3 \> }$ & \includegraphics[scale=0.2]{2ndGSX1.eps} \\
\hline
 & & & \\
 & & ${\Delta _{c}  \< 3,4 \> }$ & \includegraphics[scale=0.2]{2ndGSX1.eps} \\
\hline
\hline
\end{tabular}
\caption{Braid monodromy for the curve defined by the polynomial $(x^2-(y+\frac{3}{4})^2-1)((x-\frac{1}{2})^2-y^2-1)$}\label{BM 2nd}
\end{center}}

\clearpage

\end{document}